\documentclass[12pt, draftclsnofoot, onecolumn]{IEEEtran}
\usepackage[T1]{fontenc}   
\usepackage[utf8]{inputenc} 
\usepackage[T1]{fontenc}
\usepackage[utf8]{inputenc} 
\usepackage{url}          
\usepackage{booktabs}       
\usepackage{hyperref}
\usepackage{amsfonts}       
\usepackage{nicefrac}       
\usepackage{microtype}      
\usepackage{lipsum}		
\usepackage[final]{graphicx}
\usepackage{float}
\usepackage{cite}
\usepackage{doi}
\usepackage{wrapfig}
\usepackage{algorithm,algorithmic}
\usepackage{amsmath, amssymb, bbm, amsthm, dsfont}
\usepackage{cuted}
\usepackage{xcolor} 
\usepackage{enumitem}
\usepackage{diagbox}
\allowdisplaybreaks

\usepackage{subcaption}

\usepackage{aliascnt}
\usepackage{bm}
\usepackage{amsfonts}
\usepackage{mathrsfs}

\def\bx{\boldsymbol{x}}
\def\by{\boldsymbol{y}}

\def\bs{\boldsymbol{s}}

\def\bd{\boldsymbol{d}}
\def\bm{\boldsymbol{g}}

\def\bpi{\boldsymbol{\pi}}

\def\mF{\boldsymbol{\mathcal{F}}}

\DeclareMathOperator{\mhE}{{\scriptstyle\boldsymbol{\mathcal{E}}}}

\DeclareMathOperator{\mX}{{\scriptstyle\boldsymbol{\mathcal{X}}}}

\DeclareMathOperator{\mY}{{\scriptstyle\boldsymbol{\mathcal{Y}}}}

\DeclareMathOperator{\mZ}
{{\scriptstyle\boldsymbol{\mathcal{Z}}}}

\DeclareMathOperator{\mM}
{{\scriptstyle\boldsymbol{\mathcal{G}}}}

\DeclareMathOperator{\mD}
{{\scriptstyle\boldsymbol{\mathcal{D}}}}

\DeclareMathOperator{\mS}
{{\scriptstyle\boldsymbol{\mathcal{S}}}}

\def\mJ{\mathcal{J}}
\def\mQ{\mathcal{Q}}
\def\mT{\mathcal{T}}
\def\mU{\mathcal{U}}
\def\mW{\mathcal{W}}

\def\bxi{\boldsymbol{\xi}}
\def\bzeta{\boldsymbol{\zeta}}
\DeclareMathOperator{\mA}{{\mathcal{A}}}
\DeclareMathOperator{\mB}{{\mathcal{B}}}
\def\mC{\mathcal{C}}

\def\mP{\mathcal{P}}
\def\mR{\mathbb{R}}

\def\bxi{\boldsymbol{\xi}}

\def\mE{\mathbb{E}}

\def\cblue{\textcolor{blue}}

\newtheorem{Lemma}{Lemma}
\newtheorem{Theorem}{Theorem}

\newtheorem{Corollary}{Corollary}

\usepackage{booktabs}
\usepackage{tabularx}
\usepackage{makecell}
\usepackage[table]{xcolor}
\usepackage{threeparttable}
\usepackage{pifont} 
\definecolor{tableShade}{RGB}{245,248,252}
\newcolumntype{s}{>{\hsize=.0\hsize\centering\arraybackslash}X}
\newcolumntype{B}{>{\centering\arraybackslash}X}  
\newcolumntype{L}{>{\centering\arraybackslash}X} 

\title{\textbf{DAMA}: A Unified Accelerated Approach for Decentralized Nonconvex Minimax Optimization---Part II: Convergence and Performance Analyses}

\author{Haoyuan Cai,
        Sulaiman A. Alghunaim,
        and Ali H. Sayed
\thanks{Haoyuan Cai and Ali H. Sayed are with  the Institute of Electrical and Micro Engineering,  École Polytechnique Fédérale
de Lausanne, Switzerland (emails: \{haoyuan.cai, ali.sayed\}@epfl.ch).}
\thanks{Sulaiman A. Alghunaim is with Kuwait University, Kuwait (email:
sulaiman.alghunaim@ku.edu.kw).}
}

\begin{document}
\maketitle

\begin{abstract}
\textcolor{black}{In Part I of this work \cite{cai2025dama1},
we developed an accelerated algorithmic framework, \textbf{DAMA} (\textbf{D}ecentralized \textbf{A}ccelerated \textbf{M}inimax \textbf{A}pproach), for nonconvex Polyak–Łojasiewicz (PL) minimax optimization over decentralized multi-agent networks.
To further enhance convergence in online and offline scenarios, Part I of this work \cite{cai2025dama1} also proposed a novel accelerated gradient estimator, namely, \textbf{GRACE} (\textbf{GR}adient \textbf{AC}celeration \textbf{E}stimator), which unifies several momentum-based methods (e.g., STORM) and loopless variance-reduction techniques (e.g., PAGE, Loopless SARAH), thereby enabling accelerated gradient updates within \textbf{DAMA}.
Part I reported a unified performance bound for \textbf{DAMA} and refined guarantees for specific algorithmic instances, demonstrating the superior performance of several new variants on sparsely connected networks.
In this Part II, we focus on the convergence and performance bounds that substantiate the main results presented in Part I \cite{cai2025dama1}. In particular, we establish a unified performance bound for \textbf{DAMA} using the transformed recursion derived in Part I and subsequently refine this bound for its various special cases.}
\end{abstract}

\begin{IEEEkeywords}
Minimax optimization, unified decentralized strategies, probabilistic gradient estimator, acceleration technique
\end{IEEEkeywords}

\section{Introduction and Review of Part I \cite{cai2025dama1}}
\label{sec:introduction}
We provide a brief review of the main construction from Part I.
This work and its accompanying Part I \cite{cai2025dama1}
aim to solve the following distributed stochastic minimax problem:
\begin{subequations}
\begin{align}
&\min_{x \in \mathbb{R}^{d_1}}
\max_{y \in \mathbb{R}^{d_2}}
J(x, y)
= \frac{1}{K}\sum_{k=1}^K J_k(x,y),   
\label{main:sec1:problem}\\
&\text{where }J_k(x, y) = 
\begin{cases}
\mathbb{E}_{\bzeta_k}[Q_k(x, y; \bzeta_k)] & \textbf{(Online)}\\
\frac{1}{N_k} \sum_{s=1}^{N_k} Q_{k}(x, y; \bzeta_{k}(s)) & \textbf{(Offline)}
\end{cases}\label{main:sec1:problem_stochastic},
\end{align}
\end{subequations}
where the global risk (cost) function  $J(x,y)$ is defined as the average of $K$ local objectives $J_k(x,y)$, each privately owned by an agent $k$ in a fully decentralized network. We assume that $J(x,y)$ is $L_f$-smooth and nonconvex in $x \in \mR^{d_1}$
but 
{\em possibly} nonconcave in $y \in \mR^{d_2}$ while satisfying the
Polyak-Łojasiewicz (PL) condition
(see \cite{karimi2016linear}).
We focus on two important learning scenarios, where each agent $k$ can draw independent and identically distributed (i.i.d.) samples: $\bzeta_k$ in the online streaming setting or $\bzeta_k(s)$ in the finite-sum ( or empirical) setting, with $\bzeta_k(s)$ denoting the $s$-th local sample of agent $k$.
Problems of the type \eqref{main:sec1:problem}---\eqref{main:sec1:problem_stochastic} can be found in many applications,
such as 
generative adversarial networks (GANs) \cite{goodfellow2014gan}, reinforcement learning \cite{li2019robust}, domain adaptation \cite{acuna2021f}, and distributionally robust optimization \cite{sinha2017certifying}.
A growing body of works \cite{cai2024diffusion, xian2021faster, ghiasvand2025robust,gao2022decentralized, huang2023near, cai2025communication, liu2024decentralized,chen2024efficient,mancino2023variance} have been proposed to address  \eqref{main:sec1:problem}---\eqref{main:sec1:problem_stochastic} (see Part I \cite{cai2025dama1} for a detailed overview). However, they only
 explore limited synergies between decentralized learning strategies and gradient acceleration techniques. To bridge this gap, Part I of this work \cite{cai2025dama1} proposed a unified approach
    \textbf{DAMA} (\textbf{D}ecentralized \textbf{A}ccelerated \textbf{M}inimax \textbf{A}pproach) based on a new probabilistic gradient estimator 
    \textbf{GARCE} (\textbf{GR}adient \textbf{AC}celeration \textbf{E}stimator). 
\textbf{DAMA} subsumes existing bias-correction methods, e.g., gradient-tracking (GT), and introduces new variants for the minimax problem, e.g.,
exact diffusion (ED) \cite{yuan2018exact,yuan2018exact2}, and EXTRA \cite{shi2015extra}.
Furthermore, \textbf{GRACE} unifies the implementation of several momentum-based methods (e.g., STORM \cite{cutkosky2019momentum}) and loopless variance-reduction techniques (e.g., PAGE \cite{li2021page}, Loopless SARAH \cite{li2020convergence}) for constructing accelerated gradients to enhance \textbf{DAMA}.


We assume the reader is familiar with the notation introduced in part I \cite{cai2025dama1}. We collect the local estimates across all agents into the block variables
\begin{align}
    \mX_{i} &= \mbox{\rm col}\{\bx_{1,i}, \dots, \bx_{K,i}\} \in \mR^{Kd_1}, \quad& 
\mY _{i}=& \mbox{\rm col}\{\by_{1,i}, \dots, \by_{K,i}\} \in \mR^{Kd_2}.
\end{align}
To solve problem \eqref{main:sec1:problem}---\eqref{main:sec1:problem_stochastic}, Part I motivated a two-level primal-dual problem and proposed the following iterations $\forall i \ge 0$:
\begin{subequations}
\begin{align}
\mX_{i+1}
&= \mathcal{A}_x(\mathcal{C}_x \mX_i- \mu_{x}
\mM_{x,i}) - \mathcal{B}_x
\mD_{x,i} ,\label{main:X_update} \\
\mY_{i+1}
&= \mathcal{A}_y (\mathcal{C}_y \mY_i {\cblue + } \mu_{y}
\mM_{y,i}) - \mathcal{B}_y
\mD_{y,i}  ,\label{main:Y_update}\\
\mD_{x,i+1}
&= 
\mD_{x,i}
+\mB_x\mX_{i+1}\label{main:Ux_update} ,\\
\mD_{y,i+1}
&= 
\mD_{y,i}
+\mB_y\mY_{i+1},
\label{main:Uy_update}
\end{align}
\end{subequations}
where  $\{{\cal A}_x, {\cal A}_y, {\cal B}_x, {\cal B}_y, {\cal C}_x, {\cal C}_y\}$
are design matrices (see Table \ref{tab:matrix_choices} for possible choices). The block vectors $\mM_{x,i}, \mM_{y,i}$
serve as approximations of the local true gradients, and $\mD_{x,i} $ and $ \mD_{y,i}$ are the corresponding dual variables \cite{cai2025dama1}, with their block entries given by
\begin{subequations}
\begin{align}
\mM_{x,i} &= \mbox{col}\{\bm^x_{1,i}, \dots, \bm^x_{K,i}\} \ \in \mR^{Kd_1}, \quad& \mM_{y,i} =& \mbox{col}\{ \bm^y_{1,i}, \dots, \bm^y_{K,i}\} \in \mR^{Kd_2},\\
\mD_{x,i} &= \mbox{\rm col}\{\bd^x_{1,i}, \dots, \bd^x_{K,i}\} \in \mR^{Kd_1}, \quad& 
\mD_{y,i}=& \mbox{\rm col}\{\bd^y_{1,i}, \dots, \bd^y_{K,i}\} \in \mR^{Kd_2}.
\end{align}
\end{subequations} 

\begin{table*}[!htbp]
\centering
\caption{Matrix choices for different decentralized {\em minimax} learning strategies. Below, $W \in \mR^{K\times K}$ is a symmetric, doubly stochastic matrix.}
\label{tab:matrix_choices}
\renewcommand{\arraystretch}{1.2}
\setlength{\tabcolsep}{3pt}
\begin{tabular}{l|cccccc}
\toprule
\diagbox[width=2.5cm]{\textbf{Strategy}}{\textbf{Choices}}
& $\mA_x$ & $\mA_y$ & $\mB_x$ & $\mB_y$ & $\mC_x$ & $\mC_y$ \\ 
\toprule 
ED     & $W \otimes \mathrm{I}_{d_1}$ & $W \otimes \mathrm{I}_{d_2}$ & $(\mathrm{I}_{Kd_1} - W \otimes \mathrm{I}_{d_1})^{1/2}$ & $(\mathrm{I}_{Kd_2} -W \otimes \mathrm{I}_{d_2})^{1/2}$ & $\mathrm{I}_{Kd_1}$ & $\mathrm{I}_{Kd_2}$ \\ 
EXTRA  & $\mathrm{I}_{Kd_1}$ & $\mathrm{I}_{Kd_2}$ & $(\mathrm{I}_{Kd_1} - W \otimes \mathrm{I}_{d_1})^{1/2}$ & $(\mathrm{I}_{Kd_2} - W \otimes \mathrm{I}_{d_2})^{1/2}$ & $W \otimes \mathrm{I}_{d_1}$ & $W \otimes \mathrm{I}_{d_2}$ \\
ATC-GT & $(W \otimes \mathrm{I}_{d_1})^2$ & $(W \otimes \mathrm{I}_{d_2})^2$ & $\mathrm{I}_{Kd_1} - W \otimes \mathrm{I}_{d_1}$ & $\mathrm{I}_{Kd_2} - W \otimes \mathrm{I}_{d_2}$  & $\mathrm{I}_{Kd_1}$ & $\mathrm{I}_{Kd_2}$  \\ 
semi-ATC-GT & $W \otimes \mathrm{I}_{d_1}$ & $W \otimes \mathrm{I}_{d_2}$ & $\mathrm{I}_{Kd_1} - W \otimes \mathrm{I}_{d_1}$ & $\mathrm{I}_{Kd_2} - W \otimes \mathrm{I}_{d_2}$  & $W \otimes \mathrm{I}_{d_1}$ & $W \otimes \mathrm{I}_{d_2}$  \\ 
non-ATC-GT & $ \mathrm{I}_{Kd_1}$ & $\mathrm{I}_{Kd_2}$ & $\mathrm{I}_{Kd_1} - W \otimes \mathrm{I}_{d_1}$ & $\mathrm{I}_{Kd_2} - W \otimes \mathrm{I}_{d_2}$  & $(W \otimes \mathrm{I}_{d_1})^2$ & $(W \otimes \mathrm{I}_{d_2})^2$  \\
\bottomrule
\end{tabular}
\end{table*}


Part I of this work \cite{cai2025dama1} laid the groundwork for theoretical analysis.
Specifically, we transformed recursions \eqref{main:X_update}---\eqref{main:Uy_update} into the following equivalent form
\begin{subequations}
\begin{align}
\bx_{c,i+1} &= \bx_{c,i}-\frac{\mu_x}{K} \sum_{k=1}^{K} \bm^x_{k,i}\label{main:eq:x_finalrecursion},\\
\by_{c,i+1} &= \by_{c,i} \textcolor{blue}{ + } \frac{\mu_y}{K}\sum_{k=1}^{K}
\bm^y_{k,i}\label{main:eq:y_finalrecursion} ,\\
\mhE_{x,i+1} &= \mT_x \mhE_{x,i} 
-\frac{\mu_x}{\tau_x}\mQ^{-1}_x 
\begin{bmatrix}
0\\
\widehat{\Lambda}^{-1}_{b_x} \widehat{\Lambda}_{a_x} \widehat{\mU}^\top_x(\mM_{x,i} -\mM_{x,i+1})
\end{bmatrix} \label{main:eq:ex_finalrecursion} ,\\
\mhE_{y,i+1} &= \mT_y \mhE_{y,i} \textcolor{blue}{ + } \frac{\mu_y}{\tau_y}  \mQ^{-1}_y \begin{bmatrix}
0\\ 
\widehat{\Lambda}^{-1}_{b_y}\widehat{\Lambda}_{a_y} \widehat{\mU}^\top_y(\mM_{y,i} - \mM_{y,i+1})
\end{bmatrix}\label{main:eq:ey_finalrecursion},
\end{align}
\end{subequations}
where $\tau_x, \tau_y >0$ are arbitrary constants, $(\bx_{c,i}, \by_{c,i})$ is the network centroid, i.e., 
\begin{align}
   \bx_{c,i} \triangleq \frac{1}{K}\sum_{k=1}^K \bx_{k,i} \in \mR^{d_1},  
\by_{c,i} \triangleq \frac{1}{K}\sum_{k=1}^K \by_{k,i} \in \mR^{d_2}, 
\end{align}
and the coupled block error terms
$\mhE_{x,i}, \mhE_{y,i}$ are defined as 
\begin{align}
\mhE_{x,i} &\triangleq
   \frac{1}{\tau_x} \widehat{\mQ}^{-1}_x
\begin{bmatrix}
 \widehat{\mU}^\top_x \mX_i \\
 \widehat{\Lambda}_{b_x}^{-1} \widehat{\mU}^\top_x \mZ_{x,i}
\end{bmatrix} \in  \mathbb{R}^{2(K-1)d_1}, \ 
\mhE_{y,i}  \triangleq
   \frac{1}{\tau_y} \widehat{\mQ}^{-1}_y
\begin{bmatrix}
 \widehat{\mU}^\top_y \mY_i \\
 \widehat{\Lambda}_{b_y}^{-1} \widehat{\mU}^\top_y \mZ_{y,i}
\end{bmatrix}  \in  \mathbb{R}^{2(K-1)d_2}.
\label{main:tranform:def_couplederror}
\end{align}
Moreover, $\mZ_{x,i}$ and $\mZ_{y,i}$ are block auxiliary variables, i.e., 
\begin{align}
    \mZ_{x,i} &\triangleq \mu_x \mA_x\mM_{x,i} + \mB_x\mD_{x,i}-\mB^2_x\mX_{i} , \label{main:ZxExpression}\\
    \mZ_{y,i} &\triangleq {\cblue -}\mu_y \mA_y \mM_{y,i} + \mB_y \mD_{y,i}-\mB^2_{y} \mY_{i}.
    \label{main:ZyExpression}
\end{align}
and 
\begin{align}
&\mT_x, \mQ^{-1}_x \in \mathbb{R}^{2(K-1)d_1 \times 2(K-1)d_1},
\widehat{\Lambda}^{-1}_{b_x}, \widehat{\Lambda}_{a_x}
\in \mathbb{R}^{(K-1)d_1 \times (K-1)d_1}, \widehat{\mU}_x
\in \mathbb{R}^{Kd_1 \times (K-1)d_1}, \label{main:eigen_structure1}\\
&\mT_y, \mQ^{-1}_y \in \mathbb{R}^{2(K-1)d_2 \times 2(K-1)d_2},
\widehat{\Lambda}^{-1}_{b_y}, \widehat{\Lambda}_{a_y} 
\in \mathbb{R}^{(K-1)d_2 \times (K-1)d_2}, \widehat{\mU}_y
\in \mathbb{R}^{Kd_2 \times (K-1)d_2},
\label{main:eigen_structure2}
\end{align}
are {\em fixed} matrices that can be constructed from the eigen information of
$\{\mA_x, \mB_x, \mC_x\}$ and $\{\mA_y, \mB_y, \mC_y\}$, respectively.
In Part I \cite[Lemma 1]{cai2025dama1}, we provided the explicit construction for the matrices in \eqref{main:eigen_structure1}---\eqref{main:eigen_structure2}.

In this Part II, we aim to establish the theoretical analysis of \textbf{DAMA} to substantiate the results presented in Part I bassed on the transformed recursions \eqref{main:X_update}---\eqref{main:Uy_update}. For convenience, we introduce the notation below.

\color{black}
\subsection{Notation}
\label{sec:notations}
The purpose of this subsection is to provide a summary of the key notation used throughout the analysis.

The centroid (averaged) quantities at communication round $i$, are given by
\begin{align}
    \bx_{c,i} &\triangleq  \frac{1}{K}
    \sum_{k=1}^{K} \bx_{k,i} \in \mathbb{R}^{d_1},  &
    \by_{c,i} &\triangleq  \frac{1}{K}
    \sum_{k=1}^{K} \by_{k,i} \in \mathbb{R}^{d_2},  \\ 
   \bm^x_{c,i} &\triangleq \frac{1}{K}
   \sum_{k=1}^K 
   \bm^x_{k,i} \in \mathbb{R}^{d_1},  &
   \bm^y_{c,i} &\triangleq \frac{1}{K}
   \sum_{k=1}^K 
   \bm^y_{k,i} \in \mathbb{R}^{d_2}.
\end{align}

For network block vectors at communication round $i$, we denote
\begin{align}
\mX_{i} &\triangleq 
{\rm col}\{\bx_{1,i}, \dots, \bx_{K,i} 
\} \in \mathbb{R}^{Kd_1}, 
 \\
\mY_{i} &\triangleq 
{\rm col}\{\by_{1,i}, \dots, \by_{K,i} 
\} \in \mathbb{R}^{Kd_2},   \\
    \mX_{c,i} &\triangleq {\rm col}\{ \bx_{c,i}, \dots, \bx_{c,i}\}  \in \mathbb{R}^{Kd_1}, \\
    \mY_{c,i} &\triangleq {\rm col}\{ \by_{c,i}, \dots, \by_{c,i} \} \in \mathbb{R}^{Kd_2},  \\
\mM_{x,i} &\triangleq 
{\rm col}\{\bm^x_{1,i}, \dots, \bm^x_{K,i} 
\} \in \mathbb{R}^{Kd_1}, 
\\
\mM_{y,i} &\triangleq 
{\rm col}\{\bm^y_{1,i}, \dots, \bm^y_{K,i} 
\} \in \mathbb{R}^{Kd_2}, 
 \\
\mD_{x,i} &\triangleq
{\rm col}\{\bd^x_{1,i}, \dots, \bd^x_{K,i} 
\} \in \mathbb{R}^{Kd_1}, 
 \\ 
\mD_{y,i} &\triangleq
{\rm col}\{\bd^y_{1,i}, \dots, \bd^y_{K,i} 
\} \in \mathbb{R}^{Kd_2}. 
\end{align}
The network block (true) gradient vectors at communication round $i$
are defined as 
\begin{align}
\nabla_x \mJ(\mX_i, \mY_i)
 &\triangleq{\rm col} \{
\nabla_x J_1(\bx_{1,i}, \by_{1,i})
,\dots, 
\nabla_x J_K(\bx_{K,i}, \by_{K,i})
\} \in \mathbb{R}^{Kd_1},  \\
\nabla_y \mJ(\mX_i, \mY_i)
&\triangleq{\rm col} \{
\nabla_y J_1(\bx_{1,i}, \by_{1,i})
,\dots, 
\nabla_y J_K(\bx_{K,i}, \by_{K,i})
\} \in \mathbb{R}^{Kd_2},  \\
\nabla_x \mJ(\mX_{c,i}, \mY_{c,i})
&\triangleq{\rm col} \{
\nabla_x J_1(\bx_{c,i}, \by_{c,i})
,\dots, 
\nabla_x J_K(\bx_{c,i}, \by_{c,i})
\} \in \mathbb{R}^{Kd_1} ,\\
\nabla_y \mJ(\mX_{c,i}, \mY_{c,i})
&\triangleq{\rm col} \{
\nabla_y J_1(\bx_{c,i}, \by_{c,i})
,\dots, 
\nabla_y J_K(\bx_{c,i}, \by_{c,i})
\} \in \mathbb{R}^{Kd_2}. 
\end{align}
The network block error vector of the gradient estimator  and its centroid at communication round $i$ are denoted by
\begin{align}
\mS_{x,i} &\triangleq \mM_{x,i} - \nabla_x \mJ(\mX_i, \mY_i) \in \mathbb{R}^{Kd_1}, \\
\mS_{y,i} &\triangleq \mM_{y,i} - \nabla_y \mJ(\mX_i, \mY_i) \in \mathbb{R}^{Kd_2},  \\
 \bs^x_{c,i}  &\triangleq  \frac{1}{K}(\mathds{1}^\top_K \otimes \mathrm{I}_{d_1}) \mS_{x,i}\in \mathbb{R}^{d_1},  \\
   \bs^y_{c,i} &\triangleq  \frac{1}{K}(\mathds{1}^\top_K \otimes \mathrm{I}_{d_2}) \mS_{y,i}\in \mathbb{R}^{d_2}.  
\end{align}

Furthermore, the block combination matrices for the $x$- and the $y$-variables are defined by: 
\begin{align}
\mW_{x} &\triangleq  W \otimes \mathrm{I}_{d_1} \in \mathbb{R}^{K d_1 \times 
K d_1}, \quad \mW_{y} \triangleq W \otimes \mathrm{I}_{d_2} \in \mathbb{R}^{K d_2 \times 
Kd_2}.  
\end{align}
Consider the optimality value $P(x) = \max_y J(x,y)$ of the $y$-variable,
its optimality gap $P(x) - J(x,y)$
at the network centroid $(\bx_{c,i}, \by_{c,i})$
is defined as 
\begin{align}
\Delta_{c,i} \triangleq P(\bx_{c,i}) - J(\bx_{c,i}, \by_{c,i}).
\end{align}

\section{Theoretical analysis}
\label{sec:analysis}
In this section, we derive a performance bound for \textbf{DAMA}  in order to achieve an
$\varepsilon$-stationary point and then refine it for specific strategy instances. We first present several supporting lemmas and defer the technical details to the appendices.

\subsection{Key Lemmas}
We start by providing a bound on the consensus error of local models.
\begin{Lemma}[Consensus error]
\label{appendix:lemma:Lconsense}
Under Assumption 4 of Part I and setting 
$\tau_x = \sqrt{K}v_2, \tau_y = \sqrt{K}v_2$ for the scaled coupled error terms \eqref{main:Ux_update}–\eqref{main:Uy_update}, 
we can establish a bound for the weight
consensus error 
$\|\mX_{i} - \mX_{c,i}\|^2
+\|\mY_{i} - \mY_{c,i}\|^2$
as follows:
\begin{align}
\|\mX_i - \mX_{c,i}\|^2 
+\|\mY_i - \mY_{c,i}\|^2
\le 
Kv^2_1v^2_2\|\mhE_{x,i}\|^2
+Kv^2_1v^2_2\|\mhE_{y,i}\|^2,
\end{align}
where 
\begin{align}
    v^2_1 \triangleq \max\{ \| \widehat{\mQ}_x \|^2, \| \widehat{\mQ}_y \|^2\} , \
    v^2_2 \triangleq \max\{\| \widehat{\mQ}^{-1}_x \|^2, \| \widehat{\mQ}^{-1}_y \|^2\}.
\end{align}
\end{Lemma}
\begin{proof}
    See Appendix \ref{Appendix:proof:Lconsense}.
\end{proof}
The above lemma shows that the consensus error of the local models can be controlled by bounding the scaled coupled error in \eqref{main:Ux_update}–\eqref{main:Uy_update}. Since \textbf{DAMA} leverages the accelerated gradient estimator \textbf{GRACE}, we bound the local and averaged gradient estimation errors in the following two lemmas.
\begin{Lemma}
\label{appendix:subsec:onlinegradient}
Under Assumptions 2 and 3 of Part I, choosing smoothing factors 
$\beta_x=\beta_y=\beta\le 1, \bar{\beta}\triangleq p +\beta - p \beta \le 1$ and the Bernoulli parameter $p\le 1$ for the algorithm \textbf{DAMA},
we can establish a bound for the networked gradient estimation error as follows:
\begin{align}
&\frac{1}{T}
\sum_{i=0}^{T-1}(\mE\|\mS_{x,i}\|^2
+\mE\|\mS_{y,i}\|^2) \notag \\
&\le \frac{2}{\bar{\beta}T}
(\mE\|\mS_{x,0}\|^2
+
\mE\|\mS_{y,0}\|^2)
+ \frac{24KL^2_fv^2_1v^2_2}{b\bar{\beta}T}
\sum_{i=1}^{T}
(\mE\|\mhE_{x,i}\|^2 +\mE\|\mhE_{y,i}\|^2)
\notag \\
&\quad 
+\frac{12KL^2_fv^2_1v^2_2}{b\bar{\beta}T}
(\mE\|\mhE_{x,0}\|^2 +\mE\|\mhE_{y,0}\|^2)  + \frac{12KL^2_f}{b\bar{\beta}T}
\sum_{i=0}^{T-1}
(\mu^2_x\mE\|\bm^x_{c,i}\|^2 + \mu^2_y \mE\|\bm^y_{c,i}\|^2)
\notag \\
&\quad+\Big(
\frac{4K\beta^2\sigma^2}{b\bar{\beta}} + \frac{2pK\sigma^2}{B\bar{\beta}} \mathbb{I}_{\text{online}}\Big),
\end{align}
where 
$\mathbb{I}_{\text{online}} \in \{0, 1\}$ is the indicator of an online ($\mathbb{I}_{\text{online}} = 1$) or offline ($\mathbb{I}_{\text{online}} = 0$) setups.
\end{Lemma}
\begin{proof}
See Appendix \ref{appendix:lemma:onlinegradient}.
\end{proof}
\begin{Lemma}
\label{lemma:gradienterror_average}
Under Assumptions 2 and 3 of Part I, choosing $\beta_x=\beta_y=\beta \le 1, \bar{\beta} = p+\beta -p\bar{\beta} \le 1$ for the algorithm \textbf{DAMA},
 the averaged gradient estimation error can be bounded as follows:
\begin{align}
&\frac{1}{T}
\sum_{i=0}^{T-1}
(\mE\|\bs^x_{c,i}\|^2
+\mE\|\bs^y_{c,i}\|^2) \notag \\
&\le 
\frac{2}{T\bar{\beta}}
(\mE\|\bs^x_{c,0}\|^2
+
\mE\|\bs^y_{c,0}\|^2) 
+\frac{24
L^2_fv^2_1v^2_2}{bKT\bar{\beta}}\sum_{i=1}^{T}
(\mE\|\mhE_{x,i}\|^2 +\mE\|\mhE_{y,i}\|^2)+\frac{12
L^2_fv^2_1v^2_2}{bKT\bar{\beta}}
(\mE\|\mhE_{x,0}\|^2 \notag \\
&\quad +\mE\|\mhE_{y,0}\|^2) +
\frac{12L^2_f}{bKT\bar{\beta}}\sum_{i=0}^{T-1}
\mu^2_x\mE\|\bm^x_{c,i}\|^2 + 
\frac{12L^2_f}{bKT\bar{\beta}}\sum_{i=0}^{T-1}
\mu^2_y(1-(1-\bar{\beta})^{T-i}) \mE\|\bm^y_{c,i}\|^2 +
\frac{4\beta^2\sigma^2}{Kb\bar{\beta}} \notag \\
&\quad 
+ \frac{2p\sigma^2}{KB\bar{\beta}}
\mathbb{I}_{\text{online}},
\end{align}
where $\bar{\beta} = \beta+p-\beta p$.
\end{Lemma}
\begin{proof}
See Appendix \ref{appendix:gradienterror_average}.
\end{proof}
Lemmas~\ref{appendix:subsec:onlinegradient} and \ref{lemma:gradienterror_average} provide unified bounds on the gradient estimation error of \textbf{DAMA} in both the online stochastic and offline finite-sum settings. These results show that the error can be effectively reduced when the leading stochastic terms vanish and the variance $\sigma^2$
 is well controlled. The latter can be achieved by selecting appropriate values for the smoothing factor
$\beta$, the Bernoulli parameter 
$p$, the mini-batch size 
$b$, and the large-batch size 
$B$. Since the preceding lemmas depend on the scaled coupled error terms, we next establish a bound for these terms.
\begin{Lemma}
\label{lemma:coupled_error}
Under Assumptions 2---4 of Part I, 
choosing hyperparameters 
\begin{align}
\beta_x &= \beta_y  = \beta, \quad  p+\beta \le 1, \quad  \beta +bp \le b, \quad b \ge 1,\\
\mu_x &\le \mu_y \le \min\Big\{\frac{(1-\rho)\underline{\lambda_b}}{12L_fv_1v_2\lambda_a}, 
\frac{(1-\rho)\underline{\lambda_b}}{24L_fv_1v_2\lambda_a}\sqrt{\frac{b\bar{\beta}}{p+\beta^2}}
\Big\} ,
\end{align} 
for the algorithm \textbf{DAMA}, we get
\begin{align}
&\frac{1}{T}
\sum_{i=1}^{T}
(\mE\|\mhE_{x,i}\|^2
+\mE\|\mhE_{y,i}\|^2) \quad \Big( \text{or }
\frac{1}{T}
\sum_{i=0}^{T-1}
(\mE\|\mhE_{x,i}\|^2
+\mE\|\mhE_{y,i}\|^2)
\Big)\notag \\
&\le 
\frac{5(\mE\|\mhE_{x,0}\|^2 +\mE\|\mhE_{y,0}\|^2)}{T(1-\rho)} 
+\frac{288L^2_f\lambda^2_a\mu^2_y}{(1-\rho)^2\underline{\lambda^2_b}T}
\sum_{i=0}^{T-1}(\mu^2_x\mE\|\bm^x_{c,i}\|^2+\mu^2_y\mE\|\bm^y_{c,i}\|^2) + \frac{12\mu^2_y\lambda^2_a(p+\beta^2)}{(1-\rho)^2K\underline{\lambda^2_b}}\notag \\
&\quad 
\times
\Bigg(
\frac{2}{T\bar{\beta}}
(\mE\|\mS_{x,0}\|^2
+
\mE\|\mS_{y,0}\|^2)
+\frac{12KL^2_fv^2_1v^2_2}{bT\bar{\beta}}
(\mE\|\mhE_{x,0}\|^2 +\mE\|\mhE_{y,0}\|^2) 
\Bigg)
\notag \\
&\quad 
+ \frac{96\mu^2_y\lambda^2_a}{(1-\rho)^2\underline{\lambda^2_b}}
\Big(
\frac{p\sigma^2}{B} \mathbb{I}_{\text{online}}
+ \frac{\beta^2\sigma^2}{b}
\Big),
\end{align}
where 
\begin{align}
\rho &\triangleq \max \{\|\mT_x\|, \|\mT_y\|\} , 
    \lambda^2_a \triangleq \max \{ \|\widehat{\Lambda}_{a_x}\|^2, \|\widehat{\Lambda}_{a_y}\|^2\} , \  
    \frac{1}{\underline{\lambda^2_b}} \triangleq \max \Big\{ \|\widehat{\Lambda}^{-1}_{b_x}\|^2, \|\widehat{\Lambda}^{-1}_{b_y}\|^2\Big\}.
\end{align}
\end{Lemma}
\begin{proof}
See Appendix \ref{appendix:lemma:coupled_error}.
\end{proof}

Since $J(x,y)$ is assumed to be $L_f$-smooth and $\nu$-PL in the $y$-variable, the function $P(x) = \max_y J(x,y)$ is $L \triangleq (L_f + \frac{\kappa L_f}{2})$-smooth \cite{nouiehed2019solving} and well-defined, where $\kappa \triangleq L_f/\nu$. 
We provide the following lemma that describes the evolution of the argmax solution
 $y^o(x)= \arg \max_y J(x,y)$ over time.

\begin{Lemma}
\label{lemma:optimality_gap}
Under Assumptions 1---2 of Part I, choosing step sizes
$\mu_x \le \min\{\frac{\mu_y}{16\kappa^2}, \frac{1}{32L}\}, \mu_y \le \min\{\frac{1}{\nu}, \frac{1}{2L_f}\}$ for the algorithm
\textbf{DAMA}, we can bound the optimality gap $\Delta_{c,i}$  as follows
\begin{align}
\frac{1}{T}\sum_{i=0}^{T-1}\Delta_{c,i} 
&\le 
\frac{3}{T\nu\mu_y}
\Delta_{c,0}
+
\frac{\mu_x}{4\nu\mu_yT} 
\sum_{i=0}^{T-1}
\|\bm^x_{c,i}\|^2
-\frac{\mu_y}{4T}\sum_{i=1}^{T}
\sum_{j=0}^{i-1}\Big(
1- \frac{\nu\mu_y}{2}\Big)^{i-j-1}
\|\bm^y_{c,j}\|^2 \notag \\
&\quad +
\frac{4\kappa L_fv^2_1v^2_2}{T}
\sum_{i=0}^{T-1}
(\|\mhE_{x,i}\|^2
+\|\mhE_{y,i}\|^2)
+\frac{4}{T\nu}
\sum_{i=0}^{T-1}
\|\bs^y_{c,i}\|^2. 
\end{align}
\end{Lemma}
\begin{proof}
See Appendix \ref{subsec:optimality_gap}.
\end{proof}
The following lemma establishes the bound for the global maximum envelop function $P(x) = \max_y J(x,y)$.
\begin{Lemma}
\label{lemma:costvalue}
\label{appendix:lemma:costfunction}
    Under Assumptions 1---2 of Part I, and running the algorithm \textbf{DAMA}, 
    we can bound the cost value $P(\bx_{c,i+1})$
    as follows 
    \begin{align}
      P(\bx_{c,i+1}) &\le   P(\bx_{c,i})
-\frac{\mu_x}{2}
\|\nabla P(\bx_{c,i})\|^2
-\frac{\mu_x}{2}
(1-L\mu_x)
\|\bm^x_{c,i}\|^2
+\frac{20}{9}\mu_x \kappa L_f \Delta_{c,i}  \notag \\
&\quad + 10\mu_x L^2_fv^2_1 v^2_2 
(\|\mhE_{x,i}\|^2 + \|\mhE_{y,i}\|^2)
+ \mu_x\|\bs^x_{c,i}\|^2.
    \end{align}
\end{Lemma}
\begin{proof}
See Appendix \ref{appendix:proof:costfunction}.
\end{proof}
\subsection{Main result}
Under the supporting Lemmas \ref{appendix:subsec:onlinegradient}---\ref{lemma:costvalue}, we obtain the following main result:

\begin{Theorem}[\textbf{Main result}]
\label{main:theorem:main}
Under Assumptions 1---4 of Part I,
choosing appropriate hyperparameters for \textbf{DAMA},
it holds that  
\begin{align}
&\frac{1}{T}
\sum_{i=0}^{T-1}
\Big(\mE\|\nabla_x J(\bx_{c,i}, \by_{c,i})\|^2
+\mE\|\nabla_y J(\bx_{c,i}, \by_{c,i})\|^2\Big) \notag \\
&\le 
\mathcal{O}
\Big(
\underbrace{\frac{\mE G_{p,0}}{T\mu_x}
+  \frac{\kappa^2\mE\Delta_{c,0}}{T\mu_y} +\frac{\kappa^2(a^\prime+c^\prime\mu^2_y+f^\prime)\mu^2_y\zeta^2_0}{ T}}_{\text{Initial discrepancy}}
+\underbrace{\frac{\kappa^2b^\prime\mu^2_y\sigma^2}{ T} 
 + \kappa^2(d^\prime+e^\prime) \sigma^2}_{\text{Noisy bound}}
\Big).
\end{align}
where  $\mE G_{p,0}, \mE\Delta_{c,0}, \zeta^2_0$ account for the error arising from the initial discrepancy and
\begin{subequations}
\begin{align}
\rho &\triangleq \max \Big\{\|\mT_x\|, \|\mT_y\|\Big\}, \
\beta^\prime \triangleq  p+\beta^2, \ 
\bar{\beta} \triangleq p+\beta - p \beta, \\
\lambda_a &\triangleq \max\Big\{\lambda_{a_x}, \lambda_{b_x}\Big\}, \quad \frac{1}{\underline{\lambda^2_b}} \triangleq \max \Big\{
\frac{1}{\underline{\lambda^2_{b_x}}},\frac{1}{\underline{\lambda^2_{b_y}}}
\Big\}, \\
a' &\triangleq \frac{L^2_f}{bK\bar{\beta}(1-\rho)\underline{\lambda^2_b}}, \quad 
b^\prime \triangleq \frac{L^2_f\lambda^2_a \beta^\prime}{bb_0K\bar{\beta}^2(1-\rho)^2\underline{\lambda^2_b}}, \\
c^\prime &\triangleq \frac{L^4_f\lambda^2_a\beta^\prime}{b^2K\bar{\beta}^2(1-\rho)^2\underline{\lambda^4_b}}, \\ 
d^\prime &\triangleq \frac{L^2_f\lambda^2_a}{bK\bar{\beta}(1-\rho)^2\underline{\lambda^2_b}}
\Big( \frac{p}{B} \mathbb{I}_{\text{online}}+\frac{\beta^2}{b}\Big), \\
e^\prime &\triangleq \frac{1}{b_0\bar{\beta}KT}
+\frac{\beta^2}{Kb\bar{\beta}}
+\frac{p}{KB\bar{\beta}} \mathbb{I}_{\text{online}}, \quad 
f^\prime \triangleq \frac{L^2_f}{bK\bar{\beta}\underline{\lambda^2_b}}.
\end{align}
\end{subequations}
Here, $\mathbb{I}_{\text{online}} \in \{0,1\}$
is an indicator of online streaming setup ($\mathbb{I}_{\text{online}} =1$) or offline finite-sum setup ($\mathbb{I}_{\text{online}} =0$).
\end{Theorem}
\begin{proof}
See Appendix \ref{appendix:proofofmain}.
\end{proof}
Theorem \ref{main:theorem:main} presents the performance bound for the unified approach \textbf{DAMA}; further discussions and  comparison are provided in Part I \cite{cai2025dama1}. Refined versions of this bound for special cases appear in Appendices~\ref{appendix:corollary:storm}–\ref{appendix:corollary:L2S}.

\bibliographystyle{IEEEtran}
\bibliography{refs}

\appendix
\section{Technical Lemmas}
\subsection{Proof of Lemma \ref{appendix:lemma:Lconsense}}
\label{Appendix:proof:Lconsense}
For the weight consensus error $\|\mX_{i} -\mX_{c,i}\|^2+\|\mY_i-\mY_{c,i}\|^2$,  
it holds that 
\begin{align}
&\|\mX_i - \mX_{c,i}\|^2 
+\|\mY_i - \mY_{c,i}\|^2 \notag \\
&= 
\Big\|(\mathrm{I}_{K} \otimes \mathrm{I}_{d_1})\mX_i -  \mathds{1}_K \otimes \frac{1}{K} (\mathds{1}^\top_K \otimes \mathrm{I}_{d_1}) \mX_i \Big\|^2 + 
\Big\|(\mathrm{I}_{K} \otimes \mathrm{I}_{d_2})\mY_i - \mathds{1}_K \otimes  \frac{1}{K} (\mathds{1}^\top_K \otimes \mathrm{I}_{d_2}) \mY_i\Big\|^2 \notag \\
&=
\Big\|\Big(\Big(\mathrm{I}_{K} -  \frac{1}{K} \mathds{1}_K \otimes \mathds{1}^\top_K\Big) \otimes \mathrm{I}_{d_1}\Big) \mX_i \Big\|^2 + 
\Big\|\Big(\Big(\mathrm{I}_{K} -  \frac{1}{K} \mathds{1}_K \otimes \mathds{1}^\top_K\Big)\otimes \mathrm{I}_{d_2}\Big)\mY_i\Big\|^2 \notag \\
&\overset{(a)}{=}
\Big\|\Big(\Big(\mathrm{I}_{K} -  \frac{1}{K} \mathds{1}_K \mathds{1}^\top_K\Big) \otimes \mathrm{I}_{d_1}\Big) \mX_i \Big\|^2 + 
\Big\|\Big(\Big(\mathrm{I}_{K} -  \frac{1}{K} \mathds{1}_K\mathds{1}^\top_K\Big)\otimes \mathrm{I}_{d_2}\Big)\mY_i\Big\|^2 \notag \\
&=
\Big\|\Big(\widehat{U} \widehat{U}^\top\otimes \mathrm{I}_{d_1}\Big) \mX_i \Big\|^2 + 
\Big\|\Big(\widehat{U} \widehat{U}^\top\otimes I_{d_2}\Big)\mY_i\Big\|^2 \notag \\
&= 
\Big\|\widehat{\mU}_{x} \widehat{\mU}^\top_x \mX_i\Big\|^2
+\Big\|\widehat{\mU}_{y} \widehat{\mU}^\top_y \mY_i\Big\|^2 \notag \\
&\overset{(b)}{=} 
\Big\|\widehat{\mU}^\top_x \mX_i\|^2 
+\|\widehat{\mU}^\top_y \mY_i\|^2,
\end{align}
where $(a)$ is due to $\mathds{1}_K \otimes \mathds{1}^\top_K = \mathds{1}_K \mathds{1}^\top_K$, and $(b)$ follows from $\widehat{\mU}^\top_x \widehat{\mU}_x = \mathrm{I}_{(K-1)d_1}, \widehat{\mU}^\top_y \widehat{\mU}_y = \mathrm{I}_{(K-1)d_2}$.
Note that the above terms can be upper bounded by recalling the definitions for the scaled coupled error \eqref{main:Ux_update}–\eqref{main:Uy_update}, which gives
\begin{subequations}
\begin{align}
\|\tau_x \widehat{\mQ}_x \mhE_{x,i} \|^2
&= \Bigg\|\begin{bmatrix}
\widehat{\mU}^\top_x \mX_i \\
\widehat{\Lambda}^{-1}_{b_x} \widehat{\mU}^\top_x \mZ_{x,i} 
\end{bmatrix}
\Bigg\|^2
= \|\widehat{\mU}^\top_x \mX_i\|^2 
+\|\widehat{\Lambda}^{-1}_{b_x} \widehat{\mU}^\top_x \mZ_{x,i}\|^2 ,\\ 
\|\tau_y \widehat{\mQ}_y \mhE_{y,i} \|^2
&= \Bigg\|\begin{bmatrix}
\widehat{\mU}^\top_y \mY_i \\
\widehat{\Lambda}^{-1}_{b_y} \widehat{\mU}^\top_y \mZ_{y,i} 
\end{bmatrix}
\Bigg\|^2
= \|\widehat{\mU}^\top_y \mY_i\|^2 
+\|\widehat{\Lambda}^{-1}_{b_y} \widehat{\mU}^\top_y \mZ_{y,i}\|^2.
\end{align}
\end{subequations}
Hence, we can derive that
\begin{align}
&\|\mX_i - \mX_{c,i}\|^2 
+\|\mY_i - \mY_{c,i}\|^2  \notag \\
&=\| \widehat{\mU}^\top_x \mX_i \|^2
+\|\widehat{\mU}^\top_y \mY_i \|^2 \notag \\
&=
\|\tau_x \widehat{\mQ}_x \mhE_{x,i}\|^2
+\|\tau_y \widehat{\mQ}_y \mhE_{y,i}\|^2
- (\|\widehat{\Lambda}^{-1}_{b_x} \widehat{\mU}^\top_x \mZ_{x,i}\|^2
+\|\widehat{\Lambda}^{-1}_{b_y} \widehat{\mU}^\top_y \mZ_{y,i}\|^2) \notag \\
&\le 
\|\tau_x \widehat{\mQ}_x \mhE_{x,i}\|^2
+\|\tau_y \widehat{\mQ}_y \mhE_{y,i}\|^2 \notag \\
&\le 
\tau^2_x\|\widehat{\mQ}_x\|^2
\|\mhE_{x,i}\|^2
+
\tau^2_y\|\widehat{\mQ}_y\|^2
\|\mhE_{y,i}\|^2.
\end{align}
Denoting 
\begin{align}
v^2_{x,1} &\triangleq 
\|\widehat{\mQ}_x\|^2, \quad 
v^2_{x,2} \triangleq 
\|\widehat{\mQ}^{-1}_x\|^2, \quad v^2_{y,1} \triangleq 
\|\widehat{\mQ}_y\|^2, \quad 
v^2_{y,2} \triangleq 
\|\widehat{\mQ}^{-1}_y\|^2, \quad 
v^2_1 \triangleq \max\{
v^2_{x,1}, 
v^2_{y,1}
\},  \notag \\
v^2_2 \triangleq & \max\{
v^2_{x,2}, 
v^2_{y,2}
\}
\label{appendix:proof:lemma_notation1}
\end{align}
and setting $\tau_x = \sqrt{K} v_2, \tau_y = \sqrt{K} v_2$, we complete the proof.

\subsection{Proof of Lemma \ref{appendix:subsec:onlinegradient}}
\label{appendix:lemma:onlinegradient}
For brevity,
we focus on  bounding $\mE\|\mS_{x,i}\|^2$ and then apply a similar argument to bound $\mE\|\mS_{y,i}\|^2$. 
Recalling the \textbf{GRACE} strategy, we write its exact update at communication round $i$ based on the realization of the Bernoulli variable $\bpi_{k,i}$. Note that all agents generate
$\bpi_{k,i}$ with a shared random seed.
If $\bpi_{k,i} = 1$, $\mM_{x,i}$
is expressed as
\begin{align}
\mM_{x,i} &=  \Big\{\frac{1}{|\mathcal{E}_{k,i}|}
\sum_{\bxi_{k,i} \in \mathcal{E}_{k,i} }
\nabla_x Q_k(\bx_{k,i}, \by_{k,i}; \bxi_{k,i})\Big\}_{k=1}^{K} .
\end{align}
On the other hand, 
if $\bpi_{k,i} =0$, we have 
\begin{align}
&\mM_{x, i} = (1-\beta_x) \notag \\
&\times 
\Big(\mM_{x,i-1} -
\frac{1}{b}\Big\{\underset{{\bxi_{k,i} \in  \mathcal{E}^\prime_{k,i} }}{\sum}  \nabla_x Q_k(\bx_{k,i-1}, \by_{k,i-1}; \bxi_{k,i}) \Big\}_{k=1}^{K} \Big)  + \frac{1}{b}\Big\{\underset{{\bxi_{k,i} \in  \mathcal{E}^\prime_{k,i} }}{\sum}
\nabla_x Q_k(\bx_{k,i}, \by_{k,i}; \bxi_{k,i})\Big\}^{K}_{k=1}.
\end{align}
Using the law of total expectation and the fact that $\bpi_{k,i} = \bpi_{i} \sim  \operatorname{Bernoulli}(p)$, 
we have
\begin{align}
&\mE\|\mS_{x,i}\|^2\notag \\
=&\mE \Big[\mE_{\bpi_{k,i} \sim  \operatorname{Bernoulli}(p)}[\|\mS_{x,i}\|^2\mid \boldsymbol{\mathcal{F}}_i, \bpi_0, \dots, \bpi_{i-1}]\Big] \notag \\
=&
\mE
\Bigg[
p\Big\|\Big\{\frac{1}{|\mathcal{E}_{k,i}|} 
\sum_{\bxi_{k,i} \in \mathcal{E}_{k,i}}
\nabla_x Q_k(\bx_{k,i}, \by_{k,i}; \bxi_{k,i})\Big\}_{k=1}^{K}- \nabla_x \mJ(\mX_{i}, \mY_{i})\Big\|^2 +(1-p)
\Big\| 
(1-\beta_x)\Big(\mM_{x,i-1}
\notag \\
& -
\frac{1}{b}\Big\{\underset{{\bxi_{k,i} \in  \mathcal{E}^\prime_{k,i} }}{\sum}  \nabla_x Q_k(\bx_{k,i-1}, \by_{k,i-1}; \bxi_{k,i}) \Big\}_{k=1}^{K} \Big)  + \frac{1}{b}\Big\{\underset{{\bxi_{k,i} \in  \mathcal{E}^\prime_{k,i} }}{\sum}
\nabla_x Q_k(\bx_{k,i}, \by_{k,i}; \bxi_{k,i})\Big\}^{K}_{k=1}
\notag \\
&\quad- \nabla_x \mJ(\mX_{i}, \mY_{i})\Big\|^2
\Bigg] .
\end{align}
For convenience,
we denote 
\begin{align}
\boldsymbol{S}_1 &\triangleq \Big\| \Big\{ \frac{1}{|\mathcal{E}_{k,i}|}
\sum_{\bxi_{k,i} \in \mathcal{E}_{k,i}}
\nabla_x Q_k(\bx_{k,i}, \by_{k,i}; \bxi_{k,i})\Big\}_{k=1}^{K}- \nabla_x \mJ(\mX_{i}, \mY_{i})\Big\|^2,  \\
\boldsymbol{S}_2 &\triangleq
\Big\| 
(1-\beta_x)\Big(\mM_{x,i-1} -
\frac{1}{b}\Big\{\underset{{\bxi_{k,i} \in  \mathcal{E}^\prime_{k,i} }}{\sum}  \nabla_x Q_k(\bx_{k,i-1}, \by_{k,i-1}; \bxi_{k,i}) \Big\}_{k=1}^{K} \Big)   \notag \\
 &\quad + \frac{1}{b}\Big\{\underset{{\bxi_{k,i} \in  \mathcal{E}^\prime_{k,i} }}{\sum}
\nabla_x Q_k(\bx_{k,i}, \by_{k,i}; \bxi_{k,i})\Big\}^{K}_{k=1}
 - \nabla_x \mJ(\mX_{i}, \mY_{i})\Big\|^2.
\end{align}
In what follows, 
we bound $\boldsymbol{S}_1$
and $\boldsymbol{S}_2$, respectively.
In an offline finite-sum scenario, each agent $k$
utilizes a full batch of local samples, namely,
$\mathcal{E}_{k,i} = \{\bzeta_{k}(1), \bzeta_{k}(2),  \dots, \bzeta_{k}(N_k)\}$, which implies $\mE\boldsymbol{S}_1 = 0$; in the online stochastic setup, we have $|\mathcal{E}_{k,i}| = B <N_k$ as a full batch sample is not possible in a streaming data scenario. 
Thus, we can derive a unified bound for both scenarios
\begin{align}
\mE\boldsymbol{S}_1
&= \mE\Big\|\frac{1}{B} \Big\{
\sum_{\bxi_{k,i} \in \mathcal{E}_{k,i}}
\nabla_x Q_k(\bx_{k,i}, \by_{k,i}; \bxi_{k,i})\Big\}_{k=1}^{K}- \nabla_x \mJ(\mX_{i}, \mY_{i})\Big\|^2 \notag \\
&=
\sum_{k=1}^{K}\mE
\Big\| \frac{1}{B}\sum_{\bxi_{k,i} \in \mathcal{E}_{k,i}}
\nabla_x Q_k(\bx_{k,i}, \by_{k,i}; \bxi_{k,i})
-\nabla_x J_k(\bx_{k,i}, \by_{k,i})\Big\|^2 \notag \\
&=\sum_{k=1}^{K}\frac{1}{B^2}\mE
\Big\|\sum_{\bxi_{k,i} \in \mathcal{E}_{k,i}}
\nabla_x Q_k(\bx_{k,i}, \by_{k,i}; \bxi_{k,i})
-\nabla_x J_k(\bx_{k,i}, \by_{k,i})\Big\|^2\le 
\frac{K\sigma^2}{B} \mathbb{I}_{\text{online}},
\end{align}
where in the last inequality we used Assumption 3 of Part I, $\mathbb{I}_{\text{online}} \in \{0,1\}$ is an indicator that distinguishes the two setups: $\mathbb{I}_{\text{online}}=1$ for the online streaming regime and $\mathbb{I}_{\text{online}}=0$ for the offline finite-sum case.
We proceed to bound $\boldsymbol{S}_2$
as follows:
\begin{align}
&\mE\boldsymbol{S}_2 \notag \\
&=\mE
\Big\| 
(1-\beta_x)\Big(\mM_{x,i-1} -
\frac{1}{b}\Big\{\underset{{\bxi_{k,i} \in  \mathcal{E}^\prime_{k,i} }}{\sum}  \nabla_x Q_k(\bx_{k,i-1}, \by_{k,i-1}; \bxi_{k,i}) \Big\}_{k=1}^{K} \Big)  
\notag \\
&\quad + \frac{1}{b}\Big\{\underset{{\bxi_{k,i} \in  \mathcal{E}^\prime_{k,i} }}{\sum}
\nabla_x Q_k(\bx_{k,i}, \by_{k,i}; \bxi_{k,i})\Big\}^{K}_{k=1}- \nabla_x \mJ(\mX_{i}, \mY_{i})
\Big\|^2
\notag \\
&=
\mE
\Big\|
(1-\beta_x)\Big(\mM_{x,i-1} -
\frac{1}{b}\Big\{\underset{{\bxi_{k,i} \in  \mathcal{E}^\prime_{k,i} }}{\sum}  (\nabla_x Q_k(\bx_{k,i-1}, \by_{k,i-1}; \bxi_{k,i})
- \nabla_x Q_k(\bx_{k,i}, \by_{k,i}; \bxi_{k,i}))
\Big\}_{k=1}^{K} \Big) \notag \\
&\quad 
+
\frac{\beta_x}{b}\Big\{\underset{{\bxi_{k,i} \in  \mathcal{E}^\prime_{k,i} }}{\sum}
\nabla_x Q_k(\bx_{k,i}, \by_{k,i}; \bxi_{k,i})\Big\}^{K}_{k=1} -\nabla_x \mJ(\mX_{i}, \mY_{i})
\Big\|^2 \notag \\
&=
\mE
\Big\|
(1-\beta_x)
\Big(\mM_{x,i-1}
-\nabla_x \mJ(\mX_{i-1}, \mY_{i-1})\Big)
-(1-\beta_x)
\Big(
\frac{1}{b}\Big\{\underset{{\bxi_{k,i} \in  \mathcal{E}^\prime_{k,i} }}{\sum}\nabla_x Q_k(\bx_{k,i-1}, \by_{k,i-1}; \bxi_{k,i})
\notag \\
&\quad
- \nabla_x Q_k(\bx_{k,i}, \by_{k,i}; \bxi_{k,i}) 
\Big\}_{k=1}^{K}
-\nabla_x \mJ(\mX_{i-1}, \mY_{i-1})
+\nabla_x \mJ(\mX_{i}, \mY_{i})
\Big) \notag \\
&
\quad 
+ \beta_x \Big( \frac{1}{b}\Big\{\underset{{\bxi_{k,i} \in  \mathcal{E}^\prime_{k,i} }}{\sum}
\nabla_x Q_k(\bx_{k,i}, \by_{k,i}; \bxi_{k,i})\Big\}^{K}_{k=1} - \nabla_x \mJ(\mX_i, \mY_{i})
\Big)
\Big\|^2 \notag \\
&\overset{(a)}{=}\mE\Big\|(1-\beta_x)
\Big(\mM_{x,i-1}
-\nabla_x \mJ(\mX_{i-1}, \mY_{i-1})\Big)\Big\|^2
+ \mE\Big\|
-(1-\beta_x)
\Big(
\frac{1}{b}\Big\{\underset{{\bxi_{k,i} \in  \mathcal{E}^\prime_{k,i} }}{\sum}\nabla_x Q_k(\bx_{k,i-1}, \by_{k,i-1}; \bxi_{k,i})
\notag  
\\
&\quad
- \nabla_x Q_k(\bx_{k,i}, \by_{k,i}; \bxi_{k,i}) 
\Big\}_{k=1}^{K}
-\nabla_x \mJ(\mX_{i-1}, \mY_{i-1})
+\nabla_x \mJ(\mX_{i}, \mY_{i})
\Big) \notag \\
&
\quad 
+ \beta_x \Big( \frac{1}{b}\Big\{\underset{{\bxi_{k,i} \in  \mathcal{E}^\prime_{k,i} }}{\sum}
\nabla_x Q_k(\bx_{k,i}, \by_{k,i}; \bxi_{k,i})\Big\}^{K}_{k=1} - \nabla_x \mJ(\mX_i, \mY_{i})
\Big)
\Big\|^2 \notag \\
\newline
&\overset{(b)}{\le}(1-\beta_x)^2
\mE\|\mM_{x,i-1}
-\nabla_x \mJ(\mX_{i-1}, \mY_{i-1})\|^2 \notag \\
&\quad+ 
2(1-\beta_x)^2
\mE
\Big\|
\frac{1}{b}\Big\{\underset{{\bxi_{k,i} \in  \mathcal{E}^\prime_{k,i} }}{\sum}\nabla_x Q_k(\bx_{k,i-1}, \by_{k,i-1}; \bxi_{k,i})
\notag  
- \nabla_x Q_k(\bx_{k,i}, \by_{k,i}; \bxi_{k,i}) 
\Big\}_{k=1}^{K}
-\nabla_x \mJ(\mX_{i-1}, \mY_{i-1})
\\
&\quad+\nabla_x \mJ(\mX_{i}, \mY_{i})
\Big\|^2 +
2\beta^2_x
\mE\Big\|\frac{1}{b}\Big\{\underset{{\bxi_{k,i} \in  \mathcal{E}^\prime_{k,i} }}{\sum}
\nabla_x Q_k(\bx_{k,i}, \by_{k,i}; \bxi_{k,i})\Big\}^{K}_{k=1} - \nabla_x \mJ(\mX_i, \mY_{i})\Big\|^2 \notag \\
&\overset{(c)}{\le} 
(1-\beta_x)^2
\mE\|\mM_{x,i-1}
-\nabla_x \mJ(\mX_{i-1}, \mY_{i-1})\|^2+\frac{2(1-\beta_x)^2}{b^2}
\sum_{k=1}^K\underset{{\bxi_{k,i} \in  \mathcal{E}^\prime_{k,i} }}{\sum}
\mE
\Big\| 
 \nabla_x Q_k(\bx_{k,i-1}, \by_{k,i-1}; \bxi_{k,i})
\notag\\
&\quad 
-
\nabla_x Q_k(\bx_{k,i}, \by_{k,i}; \bxi_{k,i})-\nabla_x J_k(\bx_{k,i-1}, \by_{k,i-1})+\nabla_x J_k(\bx_{k,i}, \by_{k,i})
\Big\|^2
+ \frac{2K\beta^2_x\sigma^2}{b} \notag \\
&\overset{(d)}{\le} 
(1-\beta_x)^2
\mE\|\mM_{x,i-1}
-\nabla_x \mJ(\mX_{i-1}, \mY_{i-1})\|^2+\frac{2(1-\beta_x)^2}{b^2}
\sum_{k=1}^K\underset{{\bxi_{k,i} \in  \mathcal{E}^\prime_{k,i} }}{\sum}
\mE
\Big\| 
 \nabla_x Q_k(\bx_{k,i-1}, \by_{k,i-1}; \bxi_{k,i})
\notag\\
&\quad 
-
\nabla_x Q_k(\bx_{k,i}, \by_{k,i}; \bxi_{k,i})
\Big\|^2
+ \frac{2K\beta^2_x\sigma^2}{b} 
\notag \\
&\overset{(e)}{\le}
(1-\beta_x)^2
\mE\|\mM_{x,i-1}
-\nabla_x \mJ(\mX_{i-1}, \mY_{i-1})\|^2
+\frac{2(1-\beta_x)^2
L^2_f}{b}(\mE\|\mX_{i}-\mX_{i-1}\|^2+\mE\|\mY_i - \mY_{i-1}\|^2)
\notag \\
&\quad 
+ \frac{2K\beta^2_x\sigma^2}{b},
\label{appendix:proof:gradienterrorS_2}
\end{align}
where $(a)$ follows from Assumption 3 of Part I, $(b)$
follows from Jensen's inequality,  $(c)$
is because $\{\bxi_{k,i}\}$
is i.i.d over $k$ and $i$ (see Assumption 3 of Part I), $(d)$
follows from the inequality
$\mE\|\bxi - \mE\bxi\|^2 \le 
\mE\|\bxi\|^2$;
$(e)$ follows from the expected $L_f$-smoothness assumption.
For the weight increment
$\mE\|\mX_{i}-\mX_{i-1}\|^2 + \mE\|\mY_{i}-\mY_{i-1}\|^2$, we have 
\begin{align}
&\mE(\|\mX_{i} - \mX_{i-1}\|^2
+ \|\mY_{i} - \mY_{i-1}\|^2)\notag \\
&= 
\mE(\|\mX_{i} - \mX_{c,i} +\mX_{c,i} - \mX_{c,i-1} +
\mX_{c,i-1} - \mX_{i-1}\|^2 +
\|\mY_{i} - \mY_{c,i} +\mY_{c,i} - \mY_{c,i-1} +
\mY_{c,i-1} - \mY_{i-1}\|^2)
\notag \\
&\overset{(a)}{\le} 
3\mE(\|\mX_i - \mX_{c,i}\|^2+\|\mY_i - \mY_{c,i}\|^2)
+3K(\mu^2_x\mE\|\bm^x_{c,i-1}\|^2 + \mu^2_y\mE\|\bm^y_{c,i-1}\|^2) \notag \\
&\quad 
+3\mE(\|\mX_{i-1}-\mX_{c,i-1}\|^2 + \|\mY_{i-1}-\mY_{c,i-1}\|^2)\notag \\
&
\overset{(b)}{\le} 
3Kv^2_1v^2_2(\mE\|\mhE_{x,i}\|^2 +\mE\|\mhE_{y,i}\|^2)
+
3Kv^2_1v^2_2(\mE\|\mhE_{x,i-1}\|^2 +\mE\|\mhE_{y,i-1}\|^2) 
+ 3K(\mu^2_x\mE\|\bm^x_{c,i-1}\|^2 \notag \\
&\quad 
+ \mu^2_y\mE\|\bm^y_{c,i-1}\|^2),
\label{proof:weightincrement}
\end{align}
where $(a)$ follows from Jensen's inequality 
and recursions \eqref{main:X_update}---\eqref{main:Y_update},
and $(b)$ follows from Lemma 
\ref{appendix:lemma:Lconsense}.
Choosing $\beta_x \le 1$,
we have 
\begin{align}
&\mE\boldsymbol{S}_2
\notag \\
&\le 
(1-\beta_x)\mE\|\mM_{x,i-1}
-\nabla_x \mJ(\mX_{i-1}, \mY_{i-1})\|^2
+\frac{6
KL^2_fv^2_1v^2_2}{b}
(\mE\|\mhE_{x,i}\|^2 +\mE\|\mhE_{y,i}\|^2)
\notag \\
&\quad 
+\frac{6
KL^2_fv^2_1v^2_2}{b}
(\mE\|\mhE_{x,i-1}\|^2 +\mE\|\mhE_{y,i-1}\|^2)
+
\frac{6KL^2_f}{b}
(\mu^2_x\mE\|\bm^x_{c,i-1}\|^2 + \mu^2_y \mE\|\bm^y_{c,i-1}\|^2) +
\frac{2K\beta^2_x\sigma^2}{b} .
\end{align}
Summing up the 
results for $\mE\boldsymbol{S}_1$
and 
$\mE\boldsymbol{S}_2$,
we obtain
\begin{align}
&\mE\|\mS_{x,i}\|^2 \notag \\
&= p \mE\boldsymbol{S}_1
+(1-p)
\mE\boldsymbol{S}_2 \notag \\
&= \frac{pK\sigma^2}{B}
\mathbb{I}_{\text{online}}
+(1-p)(1-\beta_x)
\mE\|\mM_{x,i-1}- \nabla_x \mJ(\mX_{i-1}, \mY_{i-1})\|^2+\frac{6(1-p)
KL^2_fv^2_1v^2_2}{b}
(\mE\|\mhE_{x,i}\|^2 \notag \\
&\quad +\mE\|\mhE_{y,i}\|^2)
+\frac{6(1-p)
KL^2_fv^2_1v^2_2}{b}
(\mE\|\mhE_{x,i-1}\|^2 +\mE\|\mhE_{y,i-1}\|^2)
+
\frac{6(1-p)KL^2_f}{b}
(\mu^2_x\mE\|\bm^x_{c,i-1}\|^2 \notag \\
&\quad + \mu^2_y \mE\|\bm^y_{c,i-1}\|^2) +
\frac{2(1-p)K\beta^2_x\sigma^2}{b} \notag \\
&\overset{(a)}{\le} 
(1-p)(1-\beta_x)
\mE\|\mM_{x,i-1}- \nabla_x \mJ(\mX_{i-1}, \mY_{i-1})\|^2+\frac{6
KL^2_fv^2_1v^2_2}{b}
(\mE\|\mhE_{x,i}\|^2 +\mE\|\mhE_{y,i}\|^2)\notag \\
&\quad 
+\frac{6
KL^2_fv^2_1v^2_2}{b}
(\mE\|\mhE_{x,i-1}\|^2 +\mE\|\mhE_{y,i-1}\|^2)
+
\frac{6KL^2_f}{b}
(\mu^2_x\mE\|\bm^x_{c,i-1}\|^2 + \mu^2_y \mE\|\bm^y_{c,i-1}\|^2) +
\frac{2K\beta^2_x\sigma^2}{b} \notag \\
&\quad + \frac{pK\sigma^2}{B}
\mathbb{I}_{\text{online}} \notag \\
&= 
(1-p)(1-\beta_x)
\mE\|\mS_{x,i-1}\|^2+\frac{6
KL^2_fv^2_1v^2_2}{b}
(\mE\|\mhE_{x,i}\|^2 +\mE\|\mhE_{y,i}\|^2)\notag \\
&\quad 
+\frac{6
KL^2_fv^2_1v^2_2}{b}
(\mE\|\mhE_{x,i-1}\|^2 +\mE\|\mhE_{y,i-1}\|^2)
+
\frac{6KL^2_f}{b}
(\mu^2_x\mE\|\bm^x_{c,i-1}\|^2 + \mu^2_y \mE\|\bm^y_{c,i-1}\|^2) +
\frac{2K\beta^2_x\sigma^2}{b} \notag \\
&\quad + \frac{pK\sigma^2}{B}
\mathbb{I}_{\text{online}},
\end{align}
where $(a)$
follows from the fact that the Bernoulli parameter satisfies $p\le 1$.
Note that a similar argument applies to 
$\mE\|\mS_{y,i}\|^2$,
choosing $\beta_x =\beta_y = \beta \le 1$,
and summing up the results for $\mE\|\mS_{x,i}\|^2$ and $\mE\|\mS_{y,i}\|^2$,
we get 
\begin{align}
&\mE\|\mS_{x,i}\|^2
+\mE\|\mS_{y,i}\|^2 \notag \\
&\le 
(1-p)(1-\beta)
(\mE\|\mS_{x,i-1}\|^2
+
\mE\|\mS_{y,i-1}\|^2) 
+\frac{12
KL^2_fv^2_1v^2_2}{b}
(\mE\|\mhE_{x,i}\|^2 +\mE\|\mhE_{y,i}\|^2)\notag \\
&\quad+\frac{12
KL^2_fv^2_1v^2_2}{b}
(\mE\|\mhE_{x,i-1}\|^2  +\mE\|\mhE_{y,i-1}\|^2) +
\frac{12KL^2_f}{b}
(\mu^2_x\mE\|\bm^x_{c,i-1}\|^2 + \mu^2_y \mE\|\bm^y_{c,i-1}\|^2) \notag \\
&\quad+
\frac{4K\beta^2\sigma^2}{b} + \frac{2pK\sigma^2}{B}
\mathbb{I}_{\text{online}}.
\end{align}
Let 
$\bar{\beta} \triangleq 
p+\beta - p \beta$ and iterate the above inequality for $i= 1, 2, \dots,$
\begin{align}
&\mE\|\mS_{x,i}\|^2
+\mE\|\mS_{y,i}\|^2 \notag \\
&\le 
(1-\bar{\beta})^i
(\mE\|\mS_{x,0}\|^2
+
\mE\|\mS_{y,0}\|^2) +\frac{12KL^2_fv^2_1v^2_2}{b}
\sum_{j=1}^{i}(1-\bar{\beta})^{i-j}
(\mE\|\mhE_{x,j}\|^2 +\mE\|\mhE_{y,j}\|^2) \notag \\
&\quad +\frac{12KL^2_fv^2_1v^2_2}{b}
\sum_{j=0}^{i-1}(1-\bar{\beta})^{i-j-1}
(\mE\|\mhE_{x,j}\|^2 +\mE\|\mhE_{y,j}\|^2) +\frac{12KL^2_f}{b}
\sum_{j=0}^{i-1}(1-\bar{\beta})^{i-j-1}
(\mu^2_x\mE\|\bm^x_{c,j}\|^2
\notag \\
&\quad  + \mu^2_y \mE\|\bm^y_{c,j}\|^2) + \sum_{j=0}^{i-1}(1-\bar{\beta})^{i-j-1}\Big(
\frac{4K\beta^2\sigma^2}{b} + \frac{2pK\sigma^2}{B} \mathbb{I}_{\text{online}}\Big) .
\label{proof:boundSxsY}
\end{align}
Averaging the above inequality over 
$i=1,\dots, T$,  we get  
\begin{align}
&\frac{1}{T}
\sum_{i=1}^{T}(\mE\|\mS_{x,i}\|^2
+\mE\|\mS_{y,i}\|^2) \notag \\
&\le  
\frac{1}{T}
\sum_{i=1}^{T}(1-\bar{\beta})^i
(\mE\|\mS_{x,0}\|^2
+
\mE\|\mS_{y,0}\|^2) +
\frac{12KL^2_fv^2_1v^2_2}{b}\frac{1}{T}
\sum_{i=1}^{T}
\sum_{j=1}^{i}(1-\bar{\beta})^{i-j}
(\mE\|\mhE_{x,j}\|^2 +\mE\|\mhE_{y,j}\|^2) \notag \\
&\quad +
\frac{12KL^2_fv^2_1v^2_2}{b}\frac{1}{T}
\sum_{i=1}^{T}
\sum_{j=0}^{i-1}(1-\bar{\beta})^{i-j-1}
(\mE\|\mhE_{x,j}\|^2 +\mE\|\mhE_{y,j}\|^2)+\frac{12KL^2_f}{b}\frac{1}{T}
\sum_{i=1}^T\sum_{j=0}^{i-1}
(1-\bar{\beta})^{i-j-1}
\notag \\
&\quad \times (\mu^2_x\mE
\|\bm^x_{c,j}\|^2+\mu^2_y\mE\|\bm^y_{c,j}\|^2) + \frac{1}{T}
\sum_{i=1}^{T}\sum_{j=0}^{i-1}(1-\bar{\beta})^{i-j-1}\Big(
\frac{4K\beta^2\sigma^2}{b} + \frac{2pK\sigma^2}{B} \mathbb{I}_{\text{online}}\Big)
\notag \\
&\le \frac{1}{T\bar{\beta}}
(\mE\|\mS_{x,0}\|^2
+
\mE\|\mS_{y,0}\|^2)
+ \frac{12KL^2_fv^2_1v^2_2}{bT\bar{\beta}}
\sum_{i=1}^{T}
(\mE\|\mhE_{x,i}\|^2 +\mE\|\mhE_{y,i}\|^2) +\frac{12KL^2_fv^2_1v^2_2}{bT\bar{\beta}}
  \notag \\
&\quad \times 
\sum_{i=0}^{T-1}
(\mE\|\mhE_{x,i}\|^2 +\mE\|\mhE_{y,i}\|^2)+ \frac{12KL^2_f}{bT\bar{\beta}}
\sum_{i=0}^{T-1}
(\mu^2_x\mE\|\bm^x_{c,i}\|^2 + \mu^2_y \mE\|\bm^y_{c,i}\|^2)
\notag \\
&\quad 
+\Big(
\frac{4K\beta^2\sigma^2}{b\bar{\beta}} + \frac{2pK\sigma^2}{B\bar{\beta}} \mathbb{I}_{\text{online}}\Big), \\
&\le \frac{1}{T\bar{\beta}}
(\mE\|\mS_{x,0}\|^2
+
\mE\|\mS_{y,0}\|^2)
+ \frac{24KL^2_fv^2_1v^2_2}{bT\bar{\beta}}
\sum_{i=1}^{T}
(\mE\|\mhE_{x,i}\|^2 +\mE\|\mhE_{y,i}\|^2) +\frac{12KL^2_fv^2_1v^2_2}{bT\bar{\beta}}
(\mE\|\mhE_{x,0}\|^2   \notag \\
&\quad+\mE\|\mhE_{y,0}\|^2)+ \frac{12KL^2_f}{bT\bar{\beta}}
\sum_{i=0}^{T-1}
(\mu^2_x\mE\|\bm^x_{c,i}\|^2 + \mu^2_y \mE\|\bm^y_{c,i}\|^2)+\Big(
\frac{4K\beta^2\sigma^2}{b\bar{\beta}} + \frac{2pK\sigma^2}{B\bar{\beta}} \mathbb{I}_{\text{online}}\Big),
\end{align}
where the second inequality follows from 
\begin{subequations}
\begin{align}
\sum_{i=1}^{T}(1-\bar{\beta})^i &\le 
\sum_{i=0}^{\infty}
(1-\bar{\beta})^{i}
\le \frac{1}{\bar{\beta}}, \label{appendix:sequence:inequality1} \\ 
\sum_{i=1}^T\sum_{j=0}^{i-1} c^{i-j-1} w_j
&=\sum_{i=0}^{T-1}
\Big(\sum^{T-1-i}_{j=0} c^j\Big) w_i \leq \Big(\sum^{\infty}_{j=0} c^j\Big) \sum_{i=0}^{T-1}
 w_i = \frac{1}{1-c}\sum_{i=0}^{T-1} w_i,
\\ 
\sum_{i=1}^T\sum_{j=1}^{i} c^{i-j} w_j
&= 
\sum_{i=1}^{T}
\Big(\sum^{T-i}_{j=0} c^j\Big) w_i \le \Big(\sum^{\infty}_{j=0} c^j\Big)\sum_{i=1}^{T}
 w_i \le \frac{1}{1-c}\sum_{i=1}^{T} w_i ,\label{appendix:sequence:inequality2}
\end{align}
\end{subequations}
for any $c\in(0,1)$
and nonnegative sequence $\{w_j\}$.
Adding $\frac{1}{T\bar{\beta}}
(\mE\|\mS_{x,0}\|^2
+
\mE\|\mS_{y,0}\|^2)$ and choosing 
$ \bar{\beta}\le 1 \Longrightarrow \frac{1}{T\bar{\beta}} \ge \frac{1}{T}
$, we have 
\begin{align}
&\frac{1}{T}
\sum_{i=0}^{T-1}(\mE\|\mS_{x,i}\|^2
+\mE\|\mS_{y,i}\|^2) \notag \\
&\le \frac{2}{T\bar{\beta}}
(\mE\|\mS_{x,0}\|^2
+
\mE\|\mS_{y,0}\|^2)
+ \frac{24KL^2_fv^2_1v^2_2}{bT\bar{\beta}}
\sum_{i=1}^{T}
(\mE\|\mhE_{x,i}\|^2 +\mE\|\mhE_{y,i}\|^2)
 \notag \\
&\quad+\frac{12KL^2_fv^2_1v^2_2}{b\bar{\beta}T}
(\mE\|\mhE_{x,0}\|^2 +\mE\|\mhE_{y,0}\|^2) + \frac{12KL^2_f}{bT\bar{\beta}}
\sum_{i=0}^{T-1}
(\mu^2_x\mE\|\bm^x_{c,i}\|^2 + \mu^2_y \mE\|\bm^y_{c,i}\|^2)
\notag \\
&\quad+\Big(
\frac{4K\beta^2\sigma^2}{b\bar{\beta}} + \frac{2pK\sigma^2}{B\bar{\beta}} \mathbb{I}_{\text{online}}\Big).
\end{align}

\subsection{Proof of Lemma \ref{lemma:gradienterror_average}
} 
\label{appendix:gradienterror_average}
For brevity,
we focus on bounding  $\mE\|\bs^x_{c,i}\|^2$ and then apply a similar argument to bound $\mE\|\bs^y_{c,i}\|^2$. 
Note that we use the same random seed for generating the Bernoulli variable for all agents, and
if $\bpi_{k,i} = 1 \ \forall k$, the average of the gradient estimator $\bm^x_{c,i}$
can be expressed as follows
\begin{align}
\bm^x_{c,i} &= \frac{1}{K}\sum_{k=1}^{K} \frac{1}{|\mathcal{E}_{k,i}|}
\sum_{\bxi_{k,i} \in \mathcal{E}_{k,i} }
\nabla_x Q_k(\bx_{k,i}, \by_{k,i}; \bxi_{k,i}).
\end{align}
On the other hand, if $\bpi_{k,i} =0$, we get 
\begin{align}
\bm^x_{c, i} &= (1-\beta_x) \notag \\
&\times\Big(\bm^x_{c,i-1} - \frac{1}{K}\sum_{k=1}^{K}
\frac{1}{b}\underset{{\bxi_{k,i} \in  \mathcal{E}^\prime_{k,i} }}{\sum}  \nabla_x Q_k(\bx_{k,i-1}, \by_{k,i-1}; \bxi_{k,i}) \Big)  + \frac{1}{K}\sum_{k=1}^{K}\frac{1}{b}\underset{{\bxi_{k,i} \in  \mathcal{E}^\prime_{k,i} }}{\sum}
\nabla_x Q_k(\bx_{k,i}, \by_{k,i}; \bxi_{k,i}).
\end{align}
Because $\bpi_{k,i} \sim  \operatorname{Bernoulli}(p)$,
recalling the definition for $\bs^x_{c,i}$,
we have
\begin{align}
&\mE\|\bs^x_{c,i}\|^2\notag \\
=&\mE \Big[\mE_{\bpi_{k,i} \sim  \operatorname{Bernoulli}(p)}[\|\bs^x_{c,i}\|^2\mid \boldsymbol{\mathcal{F}}_i, \bpi_0, \dots, \bpi_{i-1}]\Big]  \\
=&
\mE
\Bigg[
p\Big\| \frac{1}{K}\sum_{k=1}^K\frac{1}{|\mathcal{E}_{k,i}|} 
\sum_{\bxi_{k,i} \in \mathcal{E}_{k,i}}
\nabla_x Q_k(\bx_{k,i}, \by_{k,i}; \bxi_{k,i})-  \frac{1}{K}
\sum_{k=1}^K \nabla_x J_k(\bx_{k,i},\by_{k,i})\Big\|^2 +(1-p)
\notag \\
& \Big\| 
(1-\beta_x)\Big(\bm^x_{c,i-1}-
\frac{1}{K}\sum_{k=1}^K\frac{1}{b}\underset{{\bxi_{k,i} \in  \mathcal{E}^\prime_{k,i} }}{\sum}  \nabla_x Q_k(\bx_{k,i-1}, \by_{k,i-1}; \bxi_{k,i})  \Big)  \notag \\
&\quad 
+ \frac{1}{K}
\sum_{k=1}^{K}\frac{1}{b}\underset{{\bxi_{k,i} \in  \mathcal{E}^\prime_{k,i} }}{\sum}
\nabla_x Q_k(\bx_{k,i}, \by_{k,i}; \bxi_{k,i})
 - \frac{1}{K}
\sum_{k=1}^K \nabla_x J_k(\bx_{k,i},\by_{k,i})\Big\|^2
\Bigg] .\notag
\end{align}
For simplicity,
we denote 
\begin{align}
\boldsymbol{S}_3 &\triangleq \Big\| \frac{1}{K}\sum_{k=1}^K\frac{1}{|\mathcal{E}_{k,i}|} 
\sum_{\bxi_{k,i} \in \mathcal{E}_{k,i}}
\nabla_x Q_k(\bx_{k,i}, \by_{k,i}; \bxi_{k,i})-  \frac{1}{K}
\sum_{k=1}^K \nabla_x J_k(\bx_{k,i},\by_{k,i})\Big\|^2, \notag \\
\boldsymbol{S}_4 &\triangleq
\Big\| 
(1-\beta_x)\Big(\bm^x_{c,i-1} -
\frac{1}{K}\sum_{k=1}^K\frac{1}{b}\underset{{\bxi_{k,i} \in  \mathcal{E}^\prime_{k,i} }}{\sum}  \nabla_x Q_k(\bx_{k,i-1}, \by_{k,i-1}; \bxi_{k,i})  \Big) 
\notag \\
& + \frac{1}{K}
\sum_{k=1}^{K}\frac{1}{b}\underset{{\bxi_{k,i} \in  \mathcal{E}^\prime_{k,i} }}{\sum}
\nabla_x Q_k(\bx_{k,i}, \by_{k,i}; \bxi_{k,i})
 - \frac{1}{K}
\sum_{k=1}^K \nabla_x J_k(\bx_{k,i},\by_{k,i})\Big\|^2.
\end{align}
In what follows, 
we bound $\boldsymbol{S}_3$
and $\boldsymbol{S}_4$, respectively.
If $\bpi_{k,i} =1$ in an offline finite-sum scenario, each agent computes local full-batch
gradient, which implies $\boldsymbol{S}_3 = 0$; if $\bpi_{k,i}=1$ in the online streaming data scenario,
we have  $|\mathcal{E}_{k,i}| = B <N_k$. It follows that
\begin{align}
\mE\boldsymbol{S}_3
&= \Big\| \frac{1}{K}\sum_{k=1}^K\frac{1}{B}
\sum_{\bxi_{k,i} \in \mathcal{E}_{k,i}}
\nabla_x Q_k(\bx_{k,i}, \by_{k,i}; \bxi_{k,i})-  \frac{1}{K}
\sum_{k=1}^K \nabla_x J_k(\bx_{k,i},\by_{k,i})\Big\|^2  \notag \\
&=\frac{1}{K^2B^2}\mE
\Big\|\sum_{k=1}^{K}\sum_{\bxi_{k,i} \in \mathcal{E}_{k,i}} (
\nabla_x Q_k(\bx_{k,i}, \by_{k,i}; \bxi_{k,i})
-\nabla_x J_k(\bx_{k,i}, \by_{k,i}))\Big\|^2\le 
\frac{\sigma^2}{KB} \mathbb{I}_{\text{online}}.
\end{align}
We can bound $\boldsymbol{S}_4$
as follows
\begin{align}
&\mE\boldsymbol{S}_4 \notag \\
&=\mE
\Big\| 
(1-\beta_x)\Big(\bm^x_{c,i-1} -
\frac{1}{K}\sum_{k=1}^K\frac{1}{b}\underset{{\bxi_{k,i} \in  \mathcal{E}^\prime_{k,i} }}{\sum}  \nabla_x Q_k(\bx_{k,i-1}, \by_{k,i-1}; \bxi_{k,i})  \Big) 
\notag \\
& + \frac{1}{K}
\sum_{k=1}^{K}\frac{1}{b}\underset{{\bxi_{k,i} \in  \mathcal{E}^\prime_{k,i} }}{\sum}
\nabla_x Q_k(\bx_{k,i}, \by_{k,i}; \bxi_{k,i})
 - \frac{1}{K}
\sum_{k=1}^K \nabla_x J_k(\bx_{k,i},\by_{k,i})\Big\|^2
\notag \\
&\overset{(a)}{\le}
(1-\beta_x)^2
\mE\|\bs^x_{c,i-1}\|^2 
+\frac{2(1-\beta_x)^2
L^2_f}{bK^2}(\mE\|\mX_{i}-\mX_{i-1}\|^2+\mE\|\mY_i - \mY_{i-1}\|^2)
+ \frac{2\beta^2_x\sigma^2}{Kb} 
\notag \\
&\overset{(b)}{\le} 
(1-\beta_x)\mE\|\bs^x_{c,i-1}\|^2
+\frac{6
L^2_fv^2_1v^2_2}{bK}
(\mE\|\mhE_{x,i}\|^2 +\mE\|\mhE_{y,i}\|^2)
\notag \\
&\quad 
+\frac{6
L^2_fv^2_1v^2_2}{bK}
(\mE\|\mhE_{x,i-1}\|^2 +\mE\|\mhE_{y,i-1}\|^2)
+
\frac{6L^2_f}{bK}
(\mu^2_x\mE\|\bm^x_{c,i-1}\|^2 + \mu^2_y \mE\|\bm^y_{c,i-1}\|^2) +
\frac{2\beta^2_x\sigma^2}{Kb},
\end{align}
where $(a)$ is bounded 
using an argument similar to 
equation \eqref{appendix:proof:gradienterrorS_2}, and $(b)$ follows from the inequality \eqref{proof:weightincrement}.
Putting together the
bounds for $\mE\boldsymbol{S}_3$
and 
$\mE\boldsymbol{S}_4$,
we obtain
\begin{align}
&\mE\|\bs^x_{c,i}\|^2 \notag \\
&= p \mE\boldsymbol{S}_3
+(1-p)
\mE\boldsymbol{S}_4 \notag \\
&= \frac{p\sigma^2}{KB}
\mathbb{I}_{\text{online}}
+(1-p)(1-\beta_x)
\mE\|\bs^x_{c,i-1}\|^2+\frac{6(1-p)
L^2_fv^2_1v^2_2}{bK}
(\mE\|\mhE_{x,i}\|^2 +\mE\|\mhE_{y,i}\|^2)\notag \\
&\quad 
+\frac{6(1-p)L^2_fv^2_1v^2_2}{bK}
(\mE\|\mhE_{x,i-1}\|^2 +\mE\|\mhE_{y,i-1}\|^2)
+
\frac{6(1-p)L^2_f}{bK}
(\mu^2_x\mE\|\bm^x_{c,i-1}\|^2 + \mu^2_y \mE\|\bm^y_{c,i-1}\|^2) \notag \\
&\quad 
+
\frac{2(1-p)\beta^2_x\sigma^2}{Kb} \notag \\
&\overset{(a)}{\le} 
(1-p)(1-\beta_x)
\mE\|\bs^x_{c,i-1}\|^2+\frac{6
L^2_fv^2_1v^2_2}{bK}
(\mE\|\mhE_{x,i}\|^2 +\mE\|\mhE_{y,i}\|^2)\notag \\
&\quad 
+\frac{6
L^2_fv^2_1v^2_2}{bK}
(\mE\|\mhE_{x,i-1}\|^2 +\mE\|\mhE_{y,i-1}\|^2)
+
\frac{6L^2_f}{bK}
(\mu^2_x\mE\|\bm^x_{c,i-1}\|^2 + \mu^2_y \mE\|\bm^y_{c,i-1}\|^2) +
\frac{2\beta^2_x\sigma^2}{Kb} \notag \\
&\quad 
+ \frac{p\sigma^2}{KB}
\mathbb{I}_{\text{online}},
\end{align}
where $(a)$
follows from the fact that the Bernoulli parameter $p$ satisfies $p\le 1$.
Note that a similar argument applies to 
$\mE\|\bs^y_{c,i}\|^2$,
choosing $\beta_x =\beta_y = \beta \le 1, b \ge 1$,
and summing up the results for $\mE\|\bs^x_{c,i}\|^2$ and $\mE\|\bs^y_{c,i}\|^2$,
we get 
\begin{align}
&\mE\|\bs^x_{c,i}\|^2
+\mE\|\bs^y_{c,i}\|^2 \notag \\
&\le 
(1-p)(1-\beta)
(\mE\|\bs^x_{c,i-1}\|^2
+
\mE\|\bs^y_{c,i-1}\|^2) 
+\frac{12
L^2_fv^2_1v^2_2}{bK}
(\mE\|\mhE_{x,i}\|^2 +\mE\|\mhE_{y,i}\|^2)+\frac{12
L^2_fv^2_1v^2_2}{bK}\notag \\
&\quad 
\times (\mE\|\mhE_{x,i-1}\|^2 +\mE\|\mhE_{y,i-1}\|^2) +
\frac{12L^2_f}{bK}
(\mu^2_x\mE\|\bm^x_{c,i-1}\|^2 + \mu^2_y \mE\|\bm^y_{c,i-1}\|^2) +
\frac{4\beta^2\sigma^2}{Kb} + \frac{2p\sigma^2}{KB}
\mathbb{I}_{\text{online}}.
\end{align}
Setting  $\bar{\beta} = p+\beta-p\beta \le  1$ and using an  argument similar to the proof of Lemma \ref{appendix:subsec:onlinegradient},
we can bound the average of the above inequality over $T$ communication rounds 
as follows:
\begin{align}
&\frac{1}{T}
\sum_{i=1}^{T}
(\mE\|\bs^x_{c,i}\|^2
+\mE\|\bs^y_{c,i}\|^2) \notag \\
&\le 
\frac{1}{T\bar{\beta}}
(\mE\|\bs^x_{c,0}\|^2
+
\mE\|\bs^y_{c,0}\|^2) 
+\frac{24
L^2_fv^2_1v^2_2}{bKT\bar{\beta}}\sum_{i=1}^{T}
(\mE\|\mhE_{x,i}\|^2 +\mE\|\mhE_{y,i}\|^2)\notag \\
&\quad 
+
\frac{12
L^2_fv^2_1v^2_2}{bKT\bar{\beta}}
(\mE\|\mhE_{x,0}\|^2 +\mE\|\mhE_{y,0}\|^2) 
+
\frac{12L_f^2}{bKT\bar{\beta}}
\sum_{i=0}^{T-1}
\mu^2_x\mE\|\bm^x_{c,i}\|^2 +
\frac{12L^2_f}{bKT}\sum_{i=1}^{T}
\sum_{j=0}^{i-1}(1-\bar{\beta})^{i-j-1} \notag\\
&\quad \times  \mu^2_y \mE\|\bm^y_{c,j}\|^2+
\frac{4\beta^2\sigma^2}{Kb\bar{\beta}} + \frac{2p\sigma^2}{KB\bar{\beta}}
\mathbb{I}_{\text{online}} \\
&\le 
\frac{1}{T\bar{\beta}}
(\mE\|\bs^x_{c,0}\|^2
+
\mE\|\bs^y_{c,0}\|^2) 
+\frac{24
L^2_fv^2_1v^2_2}{bKT\bar{\beta}}\sum_{i=1}^{T}
(\mE\|\mhE_{x,i}\|^2 +\mE\|\mhE_{y,i}\|^2)+
\frac{12
L^2_fv^2_1v^2_2}{bKT\bar{\beta}}
(\mE\|\mhE_{x,0}\|^2 \notag \\
&\quad +\mE\|\mhE_{y,0}\|^2) 
+
\frac{12L_f^2}{bKT\bar{\beta}}
\sum_{i=0}^{T-1}
\mu^2_x\mE\|\bm^x_{c,i}\|^2 +
\frac{12L^2_f}{bKT}\sum_{i=0}^{T-1}
\mu^2_y\mE\|\bm^y_{c,i}\|^2
\Big(
\sum_{j=0}^{T-1-i}
(1-\bar{\beta})^j
\Big)
\notag\\
&\quad +
\frac{4\beta^2\sigma^2}{Kb\bar{\beta}} + \frac{2p\sigma^2}{KB\bar{\beta}}
\mathbb{I}_{\text{online}} \notag \\
&\le 
\frac{1}{T\bar{\beta}}
(\mE\|\bs^x_{c,0}\|^2
+
\mE\|\bs^y_{c,0}\|^2) 
+\frac{24
L^2_fv^2_1v^2_2}{bKT\bar{\beta}}\sum_{i=1}^{T}
(\mE\|\mhE_{x,i}\|^2 +\mE\|\mhE_{y,i}\|^2)+
\frac{12
L^2_fv^2_1v^2_2}{bKT\bar{\beta}}
(\mE\|\mhE_{x,0}\|^2 \notag \\
&\quad 
 +\mE\|\mhE_{y,0}\|^2)+
\frac{12L_f^2}{bKT\bar{\beta}}
\sum_{i=0}^{T-1}
\mu^2_x\mE\|\bm^x_{c,i}\|^2 +
\frac{12L^2_f}{bKT\bar{\beta}}\sum_{i=0}^{T-1}
\mu^2_y
(1-(1-\bar{\beta})^{T-i})
\mE\|\bm^y_{c,i}\|^2
\Big)
\notag\\
&\quad +
\frac{4\beta^2\sigma^2}{Kb\bar{\beta}} + \frac{2p\sigma^2}{KB\bar{\beta}}
\mathbb{I}_{\text{online}}.
\end{align}
Adding 
$\frac{1}{T\bar{\beta}}
(\mE\|\bs^x_{c,0}\|^2
+
\mE\|\bs^y_{c,0}\|^2)$
on both sides and choosing 
$\frac{1}{T\bar{\beta}} \ge \frac{1}{T} \Longrightarrow \bar{\beta} \le 1$,
 we have 
\begin{align}
&\frac{1}{T}
\sum_{i=0}^{T-1}
(\mE\|\bs^x_{c,i}\|^2
+\mE\|\bs^y_{c,i}\|^2) \notag \\
&\le 
\frac{2}{T\bar{\beta}}
(\mE\|\bs^x_{c,0}\|^2
+
\mE\|\bs^y_{c,0}\|^2) 
+\frac{24
L^2_fv^2_1v^2_2}{bKT\bar{\beta}}\sum_{i=1}^{T}
(\mE\|\mhE_{x,i}\|^2 +\mE\|\mhE_{y,i}\|^2)+\frac{12
L^2_fv^2_1v^2_2}{bKT\bar{\beta}}
(\mE\|\mhE_{x,0}\|^2 \notag \\
&\quad +\mE\|\mhE_{y,0}\|^2) +
\frac{12L^2_f}{bKT\bar{\beta}}\sum_{i=0}^{T-1}
\mu^2_x\mE\|\bm^x_{c,i}\|^2 + 
\frac{12L^2_f}{bKT\bar{\beta}}\sum_{i=0}^{T-1}
\mu^2_y(1-(1-\bar{\beta})^{T-i}) \mE\|\bm^y_{c,i}\|^2 +
\frac{4\beta^2\sigma^2}{Kb\bar{\beta}} \notag \\
&\quad 
+ \frac{2p\sigma^2}{KB\bar{\beta}}
\mathbb{I}_{\text{online}}.
\end{align}

\subsection{Proof of Lemma \ref{lemma:coupled_error}}
\label{appendix:lemma:coupled_error}
Taking the squared $\ell_2$-norm on the expressions of \(\mhE_{x,i+1}, \mhE_{y,i+1}\) (see \eqref{main:Ux_update}
and \eqref{main:Uy_update}), we obtain
\begin{subequations}
\begin{align}
\|\mhE_{x,i+1}\|^2 &=\Bigg\| \mT_x \mhE_{x,i} 
-\frac{\mu_x}{\tau_x}\mQ^{-1}_x 
\begin{bmatrix}
0\\ 
\widehat{\Lambda}^{-1}_{b_x} \widehat{\Lambda}_{a_x} \widehat{\mU}^\top_x(\mM_{x,i} -\mM_{x,i+1})
\end{bmatrix} \Bigg\|^2 ,\\
\|\mhE_{y,i+1}\|^2 &= \Bigg\|\mT_y \mhE_{y,i} \textcolor{blue}{ + } \frac{\mu_y}{\tau_y}  \mQ^{-1}_y \begin{bmatrix}
0\\ 
\widehat{\Lambda}^{-1}_{b_y}\widehat{\Lambda}_{a_y} \widehat{\mU}^\top_y(\mM_{y,i} - \mM_{y,i+1})
\end{bmatrix}\Bigg\|^2.
\label{proof:ey:bound}
\end{align}
\end{subequations}
In what follows, we focus on bounding 
$\|\mhE_{y,i+1}\|^2$, and then apply a similar argument to bound
$\|\mhE_{x,i+1}\|^2$.
Let us introduce the following constants
\begin{align}
\rho &\triangleq \max \{\rho_x, \rho_y\}, \quad 
    \lambda^2_a \triangleq \max \{ \lambda^2_{a_x}, \lambda^2_{a_y}\} , \quad 
    \frac{1}{\underline{\lambda^2_b}} \triangleq \max \Big\{ \frac{1}{\underline{\lambda^2_{b_x}}}, \frac{1}{\underline{\lambda^2_{b_y}}}\Big\},\label{appendix:lemmacouple:notation}
\end{align}
where
\begin{subequations}
\begin{align}
\rho_x &\triangleq \|\mT_x\|, \ \rho_y \triangleq \|\mT_y\|, 
\lambda^2_{a_x} \triangleq \|\widehat{\Lambda}_{a_x}\|^2,  \ 
\lambda^2_{a_y}  \triangleq \|\widehat{\Lambda}_{a_y}\|^2,  \ 
\frac{1}{\underline{\lambda^2_{b_x}}} \triangleq
\|\widehat{\Lambda}^{-1}_{b_x}\|^2,   \ 
\frac{1}{\underline{\lambda^2_{b_y}}} \triangleq
\|\widehat{\Lambda}^{-1}_{b_y}\|^2 .
\end{align}
\end{subequations}
Applying Jensen's inequality
$\|a+b\|^2\le \frac{1}{t}\|a\|^2+\frac{1}{(1-t)}\|b\|^2$ to \eqref{proof:ey:bound},
and choosing 
$t = \rho_y \triangleq \|\mT_y\| <1$,
it follows that 
\begin{align}
&\|\mhE_{y,i+1}\|^2 \notag \\
&\le
\rho_y \|\mhE_{y,i} 
\|^2  +\frac{\mu^2_y}{(1-\rho_y)\tau^2_y}
\Bigg\| \mQ^{-1}_y
\begin{bmatrix}
0\\ 
\widehat{\Lambda}^{-1}_{b_y}\widehat{\Lambda}_{a_y} \widehat{\mU}^\top_y(\mM_{y,i+1} - \mM_{y,i})
\end{bmatrix}\Bigg\|^2 \notag \\
&\le 
\rho_y \|\mhE_{y,i} 
\|^2  +\frac{\mu^2_y}{(1-\rho_y)\tau^2_y}\|\mQ^{-1}_y\|^2
\Bigg\| 
\begin{bmatrix}
0\\ 
\widehat{\Lambda}^{-1}_{b_y}\widehat{\Lambda}_{a_y} \widehat{\mU}^\top_y(\mM_{y,i+1} - \mM_{y,i})
\end{bmatrix}\Bigg\|^2
\notag \\
&\le 
\rho_y \|\mhE_{y,i} 
\|^2  +\frac{\mu^2_y}{(1-\rho_y)\tau^2_y}\|\mQ^{-1}_y\|^2
\|\widehat{\Lambda}^{-1}_{b_y}\|^2\|\widehat{\Lambda}_{a_y} \|^2\|\widehat{\mU}^\top_y\|^2\|\mM_{y,i+1} - \mM_{y,i}\|^2,
\end{align}
where the last two inequalities follow from the submultiplicative property of the $\ell_2$-norm.
Taking 
expectation on both sides of the above relation and using the fact that $
\|\widehat{\mU}^\top_y\| \le 1$,
we get
\begin{align}
\mE\|\mhE_{y,i+1}\|^2 \le
\rho_y \mE\|\mhE_{y,i} 
\|^2  +
\frac{\mu^2_yv^2_{y,2}\lambda^2_{a_y}}{(1-\rho_y)\tau^2_y\underline{\lambda^2_{b_y}}}\mE\|\mM_{y,i+1}-\mM_{y,i}\|^2,
\end{align}
where $v^2_{y,2}$ is defined in 
\eqref{appendix:proof:lemma_notation1}.
We can now bound $\mE\|\mM_{y,i+1}-\mM_{y,i}\|^2$ by considering the expression of 
$\mM_{y,i+1}$ at communication round $i+1$. Note that $\mM_{y,i+1}$ is switching between two gradient estimators based on the realization of the Bernoulli variable $\bpi_{k,i+1}$. Specifically,
if the  Bernoulli variable $\bpi_{k,i+1} = 1$, $\mM_{y,i+1}$ 
is expressed as 
\begin{align}
\mM_{y,i+1} &=  \Big\{\frac{1}{|\mathcal{E}_{k,i+1}|}
\sum_{\bxi_{k,i+1} \in \mathcal{E}_{k,i+1}}
\nabla_y Q_k(\bx_{k,i+1}, \by_{k,i+1}; \bxi_{k,i+1})\Big\}_{k=1}^{K}.
\end{align} 
if $\bpi_{k,i+1} =0$,  $\mM_{y,i+1}$ is constructed by
\begin{align}
\mM_{y, i+1} = (1-\beta_y)\Big(\mM_{y,i} &-
\frac{1}{b}\Big\{\underset{{\bxi_{k,i+1} \in  \mathcal{E}^\prime_{k,i+1} }}{\sum}  \nabla_y Q_k(\bx_{k,i}, \by_{k,i}; \bxi_{k,i+1}) \Big\}_{k=1}^{K} \Big)  \notag \\
&+ \frac{1}{b}\Big\{\underset{{\bxi_{k,i+1} \in  \mathcal{E}^\prime_{k,i+1} }}{\sum}
\nabla_y Q_k(\bx_{k,i+1}, \by_{k,i+1}; \bxi_{k,i+1})\Big\}^{K}_{k=1}.
\end{align}
Using the law of total expectation and $\bpi_{k,i+1} = \bpi_{i+1} \sim  \operatorname{Bernoulli}(p)$,
we can bound the increment of the gradient estimator as follows 
\begin{align}
&\mE\|\mM_{y,i+1} -\mM_{y,i}\|^2 \notag \\
&=\mE \Big[\mE_{\bpi_{i+1} \sim \text{Bernoulli}(p)} [\|\mM_{y,i+1} - \mM_{y,i}\|^2 \mid \mF_{i+1}, \bpi_0, \dots, \bpi_{i}]\Big] \notag \\
&=\mE
\Bigg[ 
p\Big\| 
 \Big\{\frac{1}{|\mathcal{E}_{k,i+1}|}
\sum_{\bxi_{k,i+1} \in \mathcal{E}_{k,i+1}}
\nabla_y Q_k(\bx_{k,i+1}, \by_{k,i+1}; \bxi_{k,i+1})\Big\}_{k=1}^{K} - \mM_{y,i}
\Big\|^2 +(1-p)
\Big\| 
 (1-\beta_y)\Big(\mM_{y,i}\notag \\
&\quad 
 -
\frac{1}{b}\Big\{\underset{{\bxi_{k,i+1} \in  \mathcal{E}^\prime_{k,i+1} }}{\sum}  \nabla_y Q_k(\bx_{k,i}, \by_{k,i}; \bxi_{k,i+1}) \Big\}_{k=1}^{K} \Big) + \frac{1}{b}\Big\{\underset{{\bxi_{k,i+1} \in  \mathcal{E}^\prime_{k,i+1} }}{\sum}
\nabla_y Q_k(\bx_{k,i+1}, \by_{k,i+1}; \bxi_{k,i+1})\Big\}^{K}_{k=1}  
\notag \\
&\quad 
- \mM_{y,i}
\Big\|^2
\Bigg].
\end{align}
For clarity, we denote:
\begin{align}
\boldsymbol{T}_1 &\triangleq \Big\| 
 \Big\{\frac{1}{|\mathcal{E}_{k,i+1}|}
\sum_{\bxi_{k,i+1} \in \mathcal{E}_{k,i+1}}
\nabla_y Q_k(\bx_{k,i+1}, \by_{k,i+1}; \bxi_{k,i+1})\Big\}_{k=1}^{K} - \mM_{y,i}
\Big\|^2 \notag , \\
\boldsymbol{T}_2 &\triangleq \Big\| 
 (1-\beta_y)\Big(\mM_{y,i}  -
\frac{1}{b}\Big\{\underset{{\bxi_{k,i+1} \in  \mathcal{E}^\prime_{k,i+1} }}{\sum}  \nabla_y Q_k(\bx_{k,i}, \by_{k,i}; \bxi_{k,i+1}) \Big\}_{k=1}^{K} \Big) \notag \\
&\quad 
+ \frac{1}{b}\Big\{\underset{{\bxi_{k,i+1} \in  \mathcal{E}^\prime_{k,i+1} }}{\sum}
\nabla_y Q_k(\bx_{k,i+1}, \by_{k,i+1}; \bxi_{k,i+1})\Big\}^{K}_{k=1} 
- \mM_{y,i}
\Big\|^2.
\end{align}
To bound $\boldsymbol{T}_1$, we can add and subtract the true gradient as follows 
\begin{align}
\boldsymbol{T}_1 
&=
 \Big\| 
 \Big\{\frac{1}{|\mathcal{E}_{k,i+1}|}
\sum_{\bxi_{k,i+1} \in \mathcal{E}_{k,i+1}}
\nabla_y Q_k(\bx_{k,i+1}, \by_{k,i+1}; \bxi_{k,i+1})\Big\}_{k=1}^{K} - 
\nabla_y \mJ(\mX_{i+1}, \mY_{i+1})
\notag \\
&\quad+\nabla_y \mJ(\mX_{i+1},\mY_{i+1}) -\nabla_y \mJ(\mX_{i}, \mY_{i})
+\nabla_y \mJ(\mX_{i}, \mY_{i})-
\mM_{y,i}
\Big\|^2.
\end{align}
Applying Jensen's inequality to the above relation,
we get 
\begin{align}
&\boldsymbol{T}_1 \notag \\
&\overset{(a)}{\le} 
 3\Big\| 
 \Big\{\frac{1}{|\mathcal{E}_{k,i+1}|}
\sum_{\bxi_{k,i+1} \in \mathcal{E}_{k,i+1}}
\nabla_y Q_k(\bx_{k,i+1}, \by_{k,i+1}; \bxi_{k,i+1})\Big\}_{k=1}^{K} - 
\nabla_y \mJ(\mX_{i+1}, \mY_{i+1}) \Big\|^2
\notag \\
&\quad + 3
\Big\|\nabla_y \mJ(\mX_{i+1},\mY_{i+1})-\nabla_y \mJ(\mX_{i}, \mY_{i}) \Big\|^2 +
3\Big\|\mM_{y,i}- \nabla_y \mJ(\mX_{i}, \mY_{i})
\Big\|^2.
\end{align}
Using a full batch, the first term 
of $\boldsymbol{T}_1$ vanishes
when $\bpi_{i+1} = 1$ under a finite-sum scenario. On the other hand, when $\bpi_{i+1} =1$ under a online stochastic scenario, we have $|\mathcal{E}_{k,i+1}| = B < N_{k}$.
The unified bound for both scenarios is given by 
\begin{align}
&\mE\Big\|\frac{1}{B} \Big\{
\sum_{\bxi_{k,i+1} \in \mathcal{E}_{k,i+1}}
\nabla_y Q_k(\bx_{k,i+1}, \by_{k,i+1}; \bxi_{k,i+1})\Big\}_{k=1}^{K}- \nabla_y \mJ(\mX_{i+1}, \mY_{i+1})\Big\|^2 \notag \\
&=
\sum_{k=1}^{K}\mE
\Big\| \frac{1}{B}\sum_{\bxi_{k,i+1} \in \mathcal{E}_{k,i+1}} (
\nabla_y Q_k(\bx_{k,i+1}, \by_{k,i+1}; \bxi_{k,i+1})
-\nabla_y J_k(\bx_{k,i+1}, \by_{k,i+1}))\Big\|^2 \notag \\
&=\sum_{k=1}^{K}\frac{1}{B^2}\mE
\Big\|\sum_{\bxi_{k,i+1} \in \mathcal{E}_{k,i+1}}(
\nabla_y Q_k(\bx_{k,i+1}, \by_{k,i+1}; \bxi_{k,i+1})
-\nabla_y J_k(\bx_{k,i+1}, \by_{k,i+1}))\Big\|^2 
\notag \\
&\le 
\frac{K\sigma^2}{B} \mathbb{I}_{\text{online}}.
\end{align}
For the second term of $\boldsymbol{T}_1$,
we can bound it as follows 
\begin{align}
&\mE\Big\|\nabla_y \mJ(\mX_{i+1},\mY_{i+1}) -\nabla_y \mJ(\mX_{i}, \mY_{i}) \Big\|^2 \notag \\
&
\overset{(a)}{=}\mE\Big\|\nabla_y \mJ(\mX_{i+1},\mY_{i+1}) -
\nabla_y \mJ(\mX_{c,i+1},\mY_{c, i+1})+ 
\nabla_y \mJ(\mX_{c,i+1},\mY_{c, i+1}) -\nabla_y \mJ(\mX_{c,i},\mY_{c, i})
\notag \\
&\quad +\nabla_y \mJ(\mX_{c,i},\mY_{c, i})
-\nabla_y \mJ(\mX_{i}, \mY_{i}) \Big\|^2 \notag \\
&\overset{(b)}{\le} 
3
\mE\|\nabla_y \mJ(\mX_{i+1},\mY_{i+1}) -
\nabla_y \mJ(\mX_{c,i+1},\mY_{c, i+1})\|^2
+ 
3\mE\|\nabla_y \mJ(\mX_{c,i+1},\mY_{c, i+1}) -\nabla_y \mJ(\mX_{c,i},\mY_{c, i})\|^2
\notag \\
&\quad 
+3\mE\|\nabla_y \mJ(\mX_{i}, \mY_{i}) - \nabla_y \mJ(\mX_{c,i},\mY_{c, i})\|^2 \notag \\
&=
3\mE\|
{\rm col} \{\nabla_y J_k(\bx_{k,i+1}, \by_{k,i+1}) - \nabla_y J_k(\bx_{c,i+1}, \by_{c,i+1})\}_{k=1}^{K}
\|^2
+3\mE\|{\rm col} \{ \nabla_y J_k(\bx_{c,i+1}, \by_{c,i+1})\notag \\
&\quad-\nabla_y J_k(\bx_{c,i}, \by_{c,i})\}_{k=1}^K\|^2 
+3\mE\|
{\rm col} \{ 
\nabla_y J_k(\bx_{k,i}, \by_{k,i}) - \nabla_y J_k(\bx_{c,i}, \by_{c,i})
\}_{k=1}^{K}\|^2 \notag \\
&\overset{(c)}{\le} 
3L^2_f(\mE\|\mX_{i+1} - \mX_{c,i+1}\|^2+\mE\|\mY_{i+1} - \mY_{c,i+1}\|^2)
+3L^2_f
\sum_{k=1}^{K}(\mE\|\bx_{c,i+1}-\bx_{c,i}\|^2 \notag \\
&\quad +\mE\|\by_{c,i+1} - \by_{c,i}\|^2)
+3L^2_f(\mE\|\mX_i- \mX_{c,i}\|^2+\mE\|\mY_i -\mY_{c,i}\|^2) \notag \\
&\overset{(d)}{\le}  
3L^2_f(\mE\|\mX_{i+1} - \mX_{c,i+1}\|^2+\mE\|\mY_{i+1} - \mY_{c,i+1}\|^2)
+3KL^2_f
(\mu^2_x\mE\|\bm^x_{c,i}\|^2+\mu^2_y\mE\|\bm^y_{c,i}\|^2) \notag \\
&\quad 
+3L^2_f(\mE\|\mX_i- \mX_{c,i}\|^2+\mE\|\mY_i -\mY_{c,i}\|^2)
\notag \\
&\overset{(e)}{\le}  
3KL^2_f v^2_1v^2_2(\mE\|\mhE_{x,i+1}\|^2 + \mE\|\mhE_{y,i+1}\|^2)
+3KL^2_f
(\mu^2_x\mE\|\bm^x_{c,i}\|^2+\mu^2_y\mE\|\bm^y_{c,i}\|^2) \notag \\
&\quad 
+3KL^2_fv^2_1v^2_2(\mE\| \mhE_{x,i}\|^2 + \mE\|\mhE_{y,i}\|^2),
\end{align}
where $(a)$ follows from adding and subtracting the 
true gradients $\nabla_y \mJ(\mX_{c,i+1}, \mY_{c,i+1}),
\nabla_y \mJ(\mX_{c,i}, \mY_{c,i})$ into the expression;
$(b)$ follows from Jensen's inequality; $(c)$
follows from expected $L_f$-smooth assumption; $(d)$
follows from  the
recursions 
\eqref{main:X_update}---\eqref{main:Y_update}; $(e)$ follows from 
Lemma \ref{appendix:lemma:Lconsense}.
Recalling the notation for  $\mS_{y,i}$,
we can bound $\boldsymbol{T}_1$
as follows: 
\begin{align}
&\mE \boldsymbol{T}_1 \notag \\
&\le 
\frac{3K\sigma^2}{B} \mathbb{I}_{\text{online}}
+9KL^2_fv^2_1v^2_2 (\mE\|\mhE_{x,i+1}\|^2 +\mE\|\mhE_{y,i+1}\|^2)
+9KL^2_f
(\mu^2_x\mE\|\bm^x_{c,i}\|^2+\mu^2_y\mE\|\bm^y_{c,i}\|^2) \notag \\
&\quad 
+9KL^2_fv^2_1v^2_2(\mE\| \mhE_{x,i}\|^2 + \mE\|\mhE_{y,i}\|^2) + 3\mE\|\mS_{y,i}\|^2.
\end{align}
We now proceed to bound 
$\boldsymbol{T}_{2}$ as follows:
\begin{align}
\boldsymbol{T}_2 &= \Big\| 
 (1-\beta_y)\Big(\mM_{y,i}  -
\frac{1}{b}\Big\{\underset{{\bxi_{k,i+1} \in  \mathcal{E}^\prime_{k,i+1} }}{\sum}  \nabla_y Q_k(\bx_{k,i}, \by_{k,i}; \bxi_{k,i+1}) \Big\}_{k=1}^{K} \Big) \notag \\
&\quad 
+ \frac{1}{b}\Big\{\underset{{\bxi_{k,i+1} \in  \mathcal{E}^\prime_{k,i+1} }}{\sum}
\nabla_y Q_k(\bx_{k,i+1}, \by_{k,i+1}; \bxi_{k,i+1})\Big\}^{K}_{k=1} 
- \mM_{y,i}
\Big\|^2 \notag \\
&=
\Big\| 
 -\beta_y\mM_{y,i}  - 
\frac{(1-\beta_y) }{b}\Big\{\underset{{\bxi_{k,i+1} \in  \mathcal{E}^\prime_{k,i+1} }}{\sum}  \nabla_y Q_k(\bx_{k,i}, \by_{k,i}; \bxi_{k,i+1}) \Big\}_{k=1}^{K}  \notag \\
&\quad 
+ \frac{1}{b}\Big\{\underset{{\bxi_{k,i+1} \in  \mathcal{E}^\prime_{k,i+1} }}{\sum}
\nabla_y Q_k(\bx_{k,i+1}, \by_{k,i+1}; \bxi_{k,i+1})\Big\}^{K}_{k=1} 
\Big\|^2 \notag \\
& =
\Big\| 
 -\beta_y (\mM_{y,i} - \nabla_y \mJ(\mX_{i}, \mY_{i}))  
 +\frac{1}{b}\Big\{\underset{{\bxi_{k,i+1} \in  \mathcal{E}^\prime_{k,i+1} }}{\sum} \Big(
\nabla_y Q_k(\bx_{k,i+1}, \by_{k,i+1}; \bxi_{k,i+1}) \notag \\
&\quad
- \nabla_y Q_k(\bx_{k,i}, \by_{k,i}; \bxi_{k,i+1}) \Big)\Big\}^{K}_{k=1} + \beta_y \Big(
\frac{1}{b}\Big\{\underset{{\bxi_{k,i+1} \in  \mathcal{E}^\prime_{k,i+1} }}{\sum}
\nabla_y Q_k(\bx_{k,i}, \by_{k,i}; \bxi_{k,i+1})\Big\}^{K}_{k=1} \notag \\
&\quad 
-\nabla_y \mJ(\mX_i, \mY_{i})
\Big)
\Big\|^2  \notag \\
&\le 
3\beta^2_y
\|\mM_{y,i} - \nabla_y \mJ(\mX_{i}, \mY_{i})\|^2
+3 
\Big\|\frac{1}{b}\Big\{\underset{{\bxi_{k,i+1} \in  \mathcal{E}^\prime_{k,i+1} }}{\sum} \Big(
\nabla_y Q_k(\bx_{k,i+1}, \by_{k,i+1}; \bxi_{k,i+1}) \notag \\
&\quad
- \nabla_y Q_k(\bx_{k,i}, \by_{k,i}; \bxi_{k,i+1}) \Big)\Big\}^{K}_{k=1}\Big\|^2
+3 \beta^2_y
\Big\|\frac{1}{b}\Big\{\underset{{\bxi_{k,i+1} \in  \mathcal{E}^\prime_{k,i+1} }}{\sum}
\nabla_y Q_k(\bx_{k,i}, \by_{k,i}; \bxi_{k,i+1})\Big\}^{K}_{k=1} \notag \\
&\quad 
-\nabla_y \mJ(\mX_i, \mY_{i})\Big\|^2,
\end{align}
where the last inequality follows from Jensen's inequality.
Taking expectations on both sides, we can deduce that 
\begin{align}
&\mE \boldsymbol{T}_2 \notag \\
&\le 
3\beta^2_y\mE\|\mS_{y,i}\|^2
+3L^2_f(\mE\|\mX_{i+1}- \mX_{i}\|^2+\mE\|\mY_{i+1}-\mY_{i}\|^2)
+\frac{3K\beta^2_y\sigma^2}{b} \notag \\
&\overset{(a)}{\le} 
3\beta^2_y\mE\|\mS_{y,i}\|^2
+3L^2_f\Big(3\mE\|\mX_{i+1}- \mX_{c,i+1}\|^2
+3\mE\|\mX_{c, i+1}- \mX_{c,i}\|^2
+3\mE\|\mX_{i}- \mX_{c,i}\|^2
\notag \\
&\quad +3\mE\|\mY_{c,i+1}-\mY_{i+1}\|^2
+3\mE\|\mY_{c,i+1}-\mY_{c,i}\|^2
+3\mE\|\mY_{c,i}-\mY_{i}\|^2\Big)
+\frac{3K\beta^2_y\sigma^2}{b} \notag \\
&\le 
3\beta^2_y\mE\|\mS_{y,i}\|^2 +9KL^2_fv^2_1v^2_2
(\mE\|\mhE_{x,i+1}\|^2
+\mE\|\mhE_{y,i+1}\|^2
) + 9KL^2_fv^2_1v^2_2(\mE\|\mhE_{x,i}\|^2
+\mE\|\mhE_{y,i}\|^2
) \notag \\
&\quad 
+ 9KL^2_f(\mu^2_x \mE\|\bm^x_{c,i}\|^2 + \mu^2_y \mE\|\bm^y_{c,i}\|^2 ) + \frac{3K\beta^2_y\sigma^2}{b}.
\end{align}
Putting together the bound  for
$\boldsymbol{T}_1$
and $\boldsymbol{T}_2$, 
we get 
\begin{align}
\label{proof:bound_for_ey}
&\mE\|\mM_{y,i+1} - \mM_{y,i}\|^2 \notag \\
&= p \mE \boldsymbol{T}_1 + (1-p)
\mE \boldsymbol{T}_2
\notag \\
&\le 
p \Big( 
\frac{3K\sigma^2}{B} \mathbb{I}_{\text{online}}
+9KL^2_fv^2_1v^2_2 (\mE\|\mhE_{x,i+1}\|^2 + \mE\|\mhE_{y,i+1}\|^2)
+9KL^2_f
(\mu^2_x\mE\|\bm^x_{c,i}\|^2 \notag \\
&\quad 
+\mu^2_y\mE\|\bm^y_{c,i}\|^2)+9KL^2_fv^2_1v^2_2(\mE\| \mhE_{x,i}\|^2 + \mE\|\mhE_{y,i}\|^2) + 3\mE\|\mS_{y,i}\|^2
\Big)
+(1-p)
\Big( 
3\beta^2_y\mE\|\mS_{y,i}\|^2 
 \notag \\
&\quad +
9KL^2_fv^2_1v^2_2(\mE\|\mhE_{x,i+1}\|^2+\mE\|\mhE_{y,i+1}\|^2
) + 9KL^2_fv^2_1v^2_2(\mE\|\mhE_{x,i}\|^2+\mE\|\mhE_{y,i}\|^2
) \notag\\
&\quad + 9KL^2_f(\mu^2_x \mE\|\bm^x_{c,i}\|^2   + \mu^2_y \mE\|\bm^y_{c,i}\|^2 ) + \frac{3K\beta^2_y\sigma^2}{b}
\Big)
\notag \\
&\le 
9KL^2_fv^2_1v^2_2(\mE\|\mhE_{x,i+1}\|^2 + \mE\|\mhE_{y,i+1}\|^2)
+ 
9KL^2_fv^2_1v^2_2(\mE\|\mhE_{x,i}\|^2 + \mE\|\mhE_{y,i}\|^2) \notag \\
&\quad 
+ 9KL^2_f
(\mu^2_x\mE\|\bm^x_{c,i}\|^2 +\mu^2_y\mE\|\bm^y_{c,i}\|^2)
+ 
(3p +3\beta^2_y)
\mE\|\mS_{y,i}\|^2
+\frac{3Kp\sigma^2}{B}
\mathbb{I}_{\text{online}}
+\frac{3K\beta^2_y\sigma^2}{b},
\end{align}
where the last inequality follows from $p+(1-p) = 1$ and $p \le 1$.
Hence,
\(\mE\|\mhE_{y,i+1}\|^2\) can be bounded as
\begin{align}
&\mE\|\mhE_{y,i+1}\|^2
\notag \\
&\le 
\rho_y
\mE\|\mhE_{y,i}\|^2
+\frac{\mu^2_yv^2_{y,2}\lambda^2_{a_y}}{(1-\rho_y)\tau^2_y\underline{\lambda^2_{b_y}}}
\Big(
9KL^2_fv^2_1v^2_2(\mE\|\mhE_{x,i+1}\|^2 + \mE\|\mhE_{y,i+1}\|^2)
+ 
9KL^2_fv^2_1v^2_2(\mE\|\mhE_{x,i}\|^2  \notag \\
&\quad + \mE\|\mhE_{y,i}\|^2)
+ 9KL^2_f
(\mu^2_x\mE\|\bm^x_{c,i}\|^2 +\mu^2_y\mE\|\bm^y_{c,i}\|^2)
+ 
(3p +3\beta^2_y)
\mE\|\mS_{y,i}\|^2
+\frac{3Kp\sigma^2}{B}
\mathbb{I}_{\text{online}}
\notag \\
&\quad 
+\frac{3K\beta^2_y\sigma^2}{b}
\Big) .
\label{proof:ey_bound_final}
\end{align}
Similarly,
$\mE\|\mhE_{x,i+1}\|^2$
can be bounded as follows 
\begin{align}
&\mE\|\mhE_{x,i+1}\|^2
\notag \\
&\le 
\rho_x
\mE\|\mhE_{x,i}\|^2
+\frac{\mu^2_xv^2_{x,2}\lambda^2_{a_x}}{(1-\rho_x)\tau^2_x\underline{\lambda^2_{b_x}}}
\Big(
9KL^2_fv^2_1v^2_2(\mE\|\mhE_{x,i+1}\|^2 + \mE\|\mhE_{y,i+1}\|^2)
+ 
9KL^2_fv^2_1v^2_2(\mE\|\mhE_{x,i}\|^2  \notag \\
&\quad + \mE\|\mhE_{y,i}\|^2)
+ 9KL^2_f
(\mu^2_x\mE\|\bm^x_{c,i}\|^2 +\mu^2_y\mE\|\bm^y_{c,i}\|^2)
+ 
(3p +3\beta^2_x)
\mE\|\mS_{x,i}\|^2
+\frac{3Kp\sigma^2}{B}
\mathbb{I}_{\text{online}}
\notag \\
&\quad 
+\frac{3K\beta^2_x\sigma^2}{b}
\Big) .
\label{proof:ex_bound_final}
\end{align}
Recalling the notations 
$\rho, v^2_1, v^2_2, \lambda^2_a, \frac{1}{\underline{
\lambda^2_b}}$ (see \eqref{appendix:proof:lemma_notation1} and \eqref{appendix:lemmacouple:notation}) and $\tau^2_x = K v^2_2, \tau^2_y = Kv^2_2$,
it follows that:
\begin{subequations}
\begin{align}
\frac{\mu^2_xv^2_{x,2}\lambda^2_{a_x}}{(1-\rho_x)\tau^2_x\underline{\lambda^2_{b_x}}} &\le \frac{\mu^2_xv^2_{x,2}\lambda^2_a}{(1-\rho)Kv^2_{2}\underline{\lambda^2_{b}}} \le 
\frac{\mu^2_x\lambda^2_a}{(1-\rho)K\underline{\lambda^2_{b}}},\\
\frac{\mu^2_yv^2_{y,2}\lambda^2_{a_y}}{(1-\rho_y)\tau^2_y\underline{\lambda^2_{b_y}}} &\le \frac{\mu^2_yv^2_{y,2}\lambda^2_a}{(1-\rho)Kv^2_{2}\underline{\lambda^2_{b}}} \le 
\frac{\mu^2_y\lambda^2_a}{(1-\rho)K\underline{\lambda^2_{b}}} .
\end{align}
\end{subequations}
Setting $\beta_x = \beta_y = \beta$ and $\mu_x \le \mu_y$,
we can bound 
$\mE\|\mhE_{x,i+1}\|^2 + \mE\|\mhE_{y,i+1}\|^2$
as follows
\begin{align}
&\mE\|\mhE_{x,i+1}\|^2 + \mE\|\mhE_{y,i+1}\|^2
\notag \\
&\le 
\Big( 
\rho + \frac{9L^2_f(\mu^2_x+\mu^2_y)v^2_1v^2_2\lambda^2_a}{
(1-\rho)\underline{\lambda^2_{b}}
}
\Big)(\mE\|\mhE_{x,i}\|^2 + \mE\|\mhE_{y,i}\|^2) 
+ 
\frac{9L^2_f(\mu^2_x+\mu^2_y)v^2_1v^2_2\lambda^2_a}{
(1-\rho)\underline{\lambda^2_{b}}
}(\mE\|\mhE_{x,i+1}\|^2 \notag \\
&\quad + \mE\|\mhE_{y,i+1}\|^2)
+\frac{9L^2_f\lambda^2_a (\mu^2_x +\mu^2_y)}{(1-\rho)\underline{\lambda^2_{b}}}
(\mu^2_x \mE\|\bm^x_{c,i}\|^2 + \mu^2_y\mE\|\bm^y_{c,i}\|^2)
+\frac{\mu^2_y\lambda^2_a(3p+3\beta^2)}{(1-\rho)K\underline{\lambda^2_b}}
(\mE\|\mS_{x,i}\|^2 
\notag \\
&\quad 
+ \mE\|\mS_{y,i}\|^2)+
\frac{(\mu^2_x+\mu^2_y)\lambda^2_a}{(1-\rho)\underline{\lambda^2_b}}
\Big(
\frac{3p\sigma^2}{B} \mathbb{I}_{\text{online}} + \frac{3\beta^2\sigma^2}{b}\Big) .\label{proof:ey_final_second}
\end{align}
Choosing 
\begin{align}
&\frac{9L^2_f(\mu^2_x+\mu^2_y)v^2_1v^2_2\lambda^2_a}{(1-\rho) \underline{\lambda^2_b}}
\le \frac{1-\rho}{2} \notag \Longrightarrow
18L^2_f\mu^2_yv^2_1v^2_2\lambda^2_a \le \frac{(1-\rho)^2\underline{\lambda^2_b}}{2}
\Longrightarrow 
\mu_y \le \frac{(1-\rho)\underline{\lambda_b}}{6L_fv_1v_2\lambda_a}.
\end{align}
Hence,
we have 
\begin{align}
&\mE\|\mhE_{x,i+1}\|^2 + \mE\|\mhE_{y,i+1}\|^2
\notag \\
&\le 
\frac{1+\rho}{2}
(\mE\|\mhE_{x,i}\|^2 + \mE\|\mhE_{y,i}\|^2)
+ 
\frac{9L^2_f(\mu^2_x+\mu^2_y)v^2_1v^2_2\lambda^2_a}{(1-\rho)\underline{\lambda^2_b}}
(\mE\|\mhE_{x,i+1}\|^2
+\mE\|\mhE_{y,i+1}\|^2)
\notag \\ &\quad  +\frac{9L^2_f\lambda^2_a (\mu^2_x +\mu^2_y)}{(1-\rho)\underline{\lambda^2_{b}}}
(\mu^2_x \mE\|\bm^x_{c,i}\|^2 + \mu^2_y\mE\|\bm^y_{c,i}\|^2)
+\frac{\mu^2_y\lambda^2_a(3p+3\beta^2)}{(1-\rho)K\underline{\lambda^2_b}}
(\mE\|\mS_{x,i}\|^2 + \mE\|\mS_{y,i}\|^2)
\notag \\
&\quad 
+
\frac{(\mu^2_x+\mu^2_y)\lambda^2_a}{(1-\rho)\underline{\lambda^2_b}}
\Big(
\frac{3p\sigma^2}{B} \mathbb{I}_{\text{online}} + \frac{3\beta^2\sigma^2}{b}\Big) .
\end{align}
Applying the above inequality iteratively, it holds for any $i =0, 1, \dots$
\begin{align}
&\mE\|\mhE_{x,i+1}\|^2 + \mE\|\mhE_{y,i+1}\|^2
\notag \\
&\le 
\Big( 
\frac{1+\rho}{2}
\Big)^{i+1}
(\mE
\|\mhE_{x,0}\|^2
+
\mE
\|\mhE_{y,0}\|^2
)
+\frac{9L^2_f(\mu^2_x+\mu^2_y)v^2_1v^2_2\lambda^2_a}{(1-\rho)\underline{\lambda^2_b}}
\sum_{j=1}^{i+1}
\Big(
\frac{1+\rho}{2}
\Big)^{i-j+1}
(\mE\|\mhE_{x,j}\|^2  \notag \\
&\quad
+\mE\|\mhE_{y,j}\|^2)+ \frac{9L^2_f\lambda^2_a(\mu^2_x+\mu^2_y)}{(1-\rho)\underline{\lambda^2_b}}
\sum_{j=0}^{i}
\Big(\frac{1+\rho}{2}\Big)^{i-j}
(\mu^2_x \mE\|\bm^x_{c,j}\|^2+\mu^2_y\mE\|\bm^y_{c,j}\|^2)
+\frac{\mu^2_y\lambda^2_a(3p+3\beta^2)}{(1-\rho)K\underline{\lambda^2_b}}
\notag \\
&\quad \times \sum_{j=0}^{i}\Big(\frac{1+\rho}{2}\Big)^{i-j} (\mE\|\mS{_{x,j}}\|^2+\mE\|\mS_{y,j}\|^2)
+\frac{(\mu^2_x+\mu^2_y)\lambda^2_a}{(1-\rho)\underline{\lambda^2_b}}
\Big(
\frac{3p \sigma^2}{B}
\mathbb{I}_{\text{online}}
+\frac{3\beta^2\sigma^2}{b}
\Big)\notag\\
&\quad \times \sum_{j=0}^{i}\Big(\frac{1+\rho}{2}\Big)^{i-j}.
\end{align}
Averaging the above inequality over 
$i = 0, \dots, T-1$, it holds that 
\begin{align}
&\frac{1}{T}\sum_{i=0}^{T-1}\mE\|\mhE_{x,i+1}\|^2 + \mE\|\mhE_{y,i+1}\|^2
\notag \\
&\le \frac{1}{T}
\sum_{i=0}^{T-1}
\Big( 
\frac{1+\rho}{2}
\Big)^{i+1}
(\mE
\|\mhE_{x,0}\|^2
+
\mE
\|\mhE_{y,0}\|^2
)
+\frac{9L^2_f(\mu^2_x+\mu^2_y)v^2_1v^2_2\lambda^2_a}{(1-\rho)\underline{\lambda^2_b}T}
\sum_{i=0}^{T-1}\sum_{j=1}^{i+1}
\notag \\
&\quad 
\times \Big(
\frac{1+\rho}{2}
\Big)^{i-j+1}
(\mE\|\mhE_{x,j}\|^2 +\mE\|\mhE_{y,j}\|^2) 
+ \frac{9L^2_f\lambda^2_a(\mu^2_x+\mu^2_y)}{(1-\rho)\underline{\lambda^2_b}T}
\sum_{i=0}^{T-1}
\sum_{j=0}^{i}
\Big(\frac{1+\rho}{2}\Big)^{i-j}
\notag \\
&\quad \times (\mu^2_x \mE\|\bm^x_{c,j}\|^2+\mu^2_y\mE\|\bm^y_{c,j}\|^2)
+\frac{\mu^2_y\lambda^2_a(3p+3\beta^2)}{(1-\rho)K\underline{\lambda^2_b}T}
\sum_{i=0}^{T-1}
\sum_{j=0}^{i}\Big(\frac{1+\rho}{2}\Big)^{i-j}(\mE\|\mS{_{x,j}}\|^2+
\notag \\
&\quad  \mE\|\mS_{y,j}\|^2)
+\frac{(\mu^2_x+\mu^2_y)\lambda^2_a}{(1-\rho)\underline{\lambda^2_b}T}
\Big(
\frac{3p \sigma^2}{B}
\mathbb{I}_{\text{online}}
+\frac{3\beta^2\sigma^2}{b}
\Big)
\sum_{i=0}^{T-1}\sum_{j=0}^{i}\Big(\frac{1+\rho}{2}\Big)^{i-j}.
\end{align}
Using \eqref{appendix:sequence:inequality1}---\eqref{appendix:sequence:inequality2}, we can derive
\begin{subequations}
\begin{align}
\sum_{i=0}^{T-1}
\Big(\frac{1+\rho}{2}\Big)^{i+1}
&\le \frac{2}{1-\rho} ,\\
\sum_{i=0}^{T-1}
\sum_{j=1}^{i+1}
\Big( \frac{1+\rho}{2}\Big)^{i-j+1}
(\mE\|\mhE_{x,j}\|^2 + \mE\|\mhE_{y,j}\|^2)
&\le 
\frac{2}{1-\rho}\sum_{i=1}^{T}
(\mE\|\mhE_{x,i}\|^2 + \mE\|\mhE_{y,i}\|^2),  \\
\sum_{i=0}^{T-1}\sum_{j=0}^i
\Big(\frac{1+\rho}{2}\Big)^{i-j}
(\mu^2_x\mE\|\bm^x_{c,j}\|^2 +\mu^2_y\mE\|\bm^y_{c,j}\|^2)
&\le \frac{2}{1-\rho}
\sum_{i=0}^{T-1}
(\mu^2_x\mE\|\bm^x_{c,i}\|^2 +\mu^2_y\mE\|\bm^y_{c,i}\|^2),  \\
\sum_{i=0}^{T-1}
\sum_{j=0}^{i}
\Big(\frac{1+\rho}{2}\Big)^{i-j}
(\mE\|\mS_{x,j}\|^2 + 
(\mE\|\mS_{x,j}\|^2 ) &\le  
\frac{2}{1-\rho}
\sum_{i=0}^{T-1}
(\mE\|\mS_{x,i}\|^2 +  
(\mE\|\mS_{y,i}\|^2 ),  \\
\sum_{i=0}^{T-1}\sum_{j=0}^i\Big(\frac{1+\rho}{2}\Big)^{i-j} &\le 
\frac{2T}{1-\rho}.
\end{align}
\end{subequations}
It follows that
\begin{align}
&\frac{1}{T}
\sum_{i=1}^{T}
\mE\|\mhE_{x,i}\|^2
+\mE\|\mhE_{y,i}\|^2 \notag \\
&\le 
\frac{2(\mE\|\mhE_{x,0}\|^2 +\mE\|\mhE_{y,0}\|^2)}{T(1-\rho)}
+ \frac{18L^2_f(\mu^2_x+\mu^2_y)v^2_1v^2_2\lambda^2_a}{(1-\rho)^2\underline{\lambda^2_b}T}
\sum_{i=1}^{T}
(\mE\|\mhE_{x,i}\|^2+\mE\|\mhE_{y,i}\|^2)
\notag \\
&\quad 
+\frac{18L^2_f\lambda^2_a(\mu^2_x+\mu^2_y)}{(1-\rho)^2\underline{\lambda^2_b}T}
\sum_{i=0}^{T-1}(\mu^2_x\mE\|\bm^x_{c,i}\|^2+\mu^2_y\mE\|\bm^y_{c,i}\|^2)+ \frac{\mu^2_y\lambda^2_a(6p+6\beta^2)}{(1-\rho)^2K\underline{\lambda^2_b}T}
\sum_{i=0}^{T-1}(\mE\|\mS_{x,i}\|^2 
\notag \\
&\quad +\mE\|\mS_{y,i}\|^2)
+ \frac{(\mu^2_x+\mu^2_y)\lambda^2_a}{(1-\rho)^2\underline{\lambda^2_b}}
\Big(
\frac{6p\sigma^2}{B} \mathbb{I}_{\text{online}}
+ \frac{6\beta^2\sigma^2}{b}
\Big).
\end{align}
Invoking 
Lemma \ref{appendix:subsec:onlinegradient},
we get 
\begin{align}
&\frac{1}{T}
\sum_{i=1}^{T}
\mE\|\mhE_{x,i}\|^2
+\mE\|\mhE_{y,i}\|^2 \notag \\
&\le 
\frac{2(\mE\|\mhE_{x,0}\|^2 +\mE\|\mhE_{y,0}\|^2)}{T(1-\rho)}
+ \frac{18L^2_f(\mu^2_x+\mu^2_y)v^2_1v^2_2\lambda^2_a}{(1-\rho)^2\underline{\lambda^2_b}T}
\sum_{i=1}^{T}
(\mE\|\mhE_{x,i}\|^2+\mE\|\mhE_{y,i}\|^2)
\notag \\
&\quad 
+\frac{18L^2_f\lambda^2_a(\mu^2_x+\mu^2_y)}{(1-\rho)^2\underline{\lambda^2_b}T}
\sum_{i=0}^{T-1}(\mu^2_x\mE\|\bm^x_{c,i}\|^2+\mu^2_y\mE\|\bm^y_{c,i}\|^2)+ \frac{\mu^2_y\lambda^2_a(6p+6\beta^2)}{(1-\rho)^2K\underline{\lambda^2_b}}
\Bigg(
\frac{2}{T\bar{\beta}}
(\mE\|\mS_{x,0}\|^2
\notag \\
&\quad+
\mE\|\mS_{y,0}\|^2)+ \frac{24KL^2_fv^2_1v^2_2}{bT\bar{\beta}}
\sum_{i=1}^{T}
(\mE\|\mhE_{x,i}\|^2 +\mE\|\mhE_{y,i}\|^2)
+\frac{12KL^2_fv^2_1v^2_2}{b\bar{\beta}T}
(\mE\|\mhE_{x,0}\|^2 
\notag \\
&\quad+\mE\|\mhE_{y,0}\|^2)  + \frac{12KL^2_f}{bT\bar{\beta}}
\sum_{i=0}^{T-1}
(\mu^2_x\mE\|\bm^x_{c,i}\|^2 + \mu^2_y \mE\|\bm^y_{c,i}\|^2)
+\Big(
\frac{4K\beta^2\sigma^2}{b\bar{\beta}} + \frac{2pK\sigma^2}{B\bar{\beta}} \mathbb{I}_{\text{online}}\Big)
\Bigg)
\notag \\
&\quad 
+ \frac{(\mu^2_x+\mu^2_y)\lambda^2_a}{(1-\rho)^2\underline{\lambda^2_b}}
\Big(
\frac{6p\sigma^2}{B} \mathbb{I}_{\text{online}}
+ \frac{6\beta^2\sigma^2}{b}
\Big).
\end{align}
Recalling that $\mu_x \le \mu_y$ and choosing
\begin{align}
&\frac{18L^2_f(\mu^2_x+\mu^2_y)v^2_1v^2_2\lambda^2_a}{(1-\rho)^2\underline{\lambda^2_b}} \le \frac{1}{4}  \Longrightarrow
\frac{36L^2_f\mu^2_yv^2_1v^2_2\lambda^2_a}{(1-\rho)^2\underline{\lambda^2_b}} \le \frac{1}{4} \Longrightarrow
\mu_y \le 
\frac{(1-\rho)\underline{\lambda_b}}{12L_fv_1v_2\lambda_a},
\\
&\frac{\mu^2_y\lambda^2_a(6p+6\beta
^2)}{(1-\rho
)^2K\underline{\lambda^2_b}}
\times  \frac{24KL^2_fv^2_1v^2_2}{b\bar{\beta}} \le \frac{1}{4} \Longrightarrow
\mu_y \le \frac{(1-\rho)\underline{\lambda
_b}}{24L_fv_1v_2\lambda_a}
\sqrt{\frac{b\bar{\beta}}{p+\beta^2}} ,
\end{align}
\begin{align}
\frac{\mu^2_y\lambda^2_a(6p+6\beta
^2)}{(1-\rho
)^2K\underline{\lambda^2_b}} \times 
\frac{12KL^2_f}{b\bar{\beta}} \le \frac{72L^2_f\lambda^2_a\mu^2_y}{(1-\rho)^2\underline{\lambda^2_b}}
&\Longrightarrow p+\beta^2 \le b\bar{\beta} \notag \\
&\Longrightarrow 
p+\beta^2 \le b(p+\beta - p \beta) \notag \\
&
\Longrightarrow
0 \le (b-1) p + \beta(b(1-p)-\beta) \notag\\
&\Longrightarrow
b\ge 1, \beta + bp \le b,
\label{proof:consensus_condition1}
\end{align}
and 
\begin{align}
&\frac{\mu^2_y\lambda^2_a(6p+6\beta
^2)}{(1-\rho
)^2K\underline{\lambda^2_b}} 
\times \frac{4K\beta^2\sigma^2}{b\bar{\beta}}
\le \frac{24\mu^2_y\beta^2\lambda^2_a\sigma^2}{(1-\rho)^2\underline{\lambda^2_b}b} \Longrightarrow p+\beta^2 \le \bar{\beta} \Longrightarrow p+\beta \le 1 ,\\
&\frac{\mu^2_y\lambda^2_a(6p+6\beta
^2)}{(1-\rho
)^2K\underline{\lambda^2_b}}
\times \frac{2pK\sigma^2}{B\bar{\beta}} \le \frac{12\mu^2_yp\lambda^2_a\sigma^2}{(1-\rho)^2\underline{\lambda^2_b}B} \Longrightarrow p+\beta^2 \le \bar{\beta} \Longrightarrow p +\beta \le 1 .
\end{align}
Hence,
we get
\begin{align}
&\frac{1}{T}
\sum_{i=1}^{T}
\mE\|\mhE_{x,i}\|^2
+\mE\|\mhE_{y,i}\|^2 \notag \\
&\le 
\frac{2(\mE\|\mhE_{x,0}\|^2 +\mE\|\mhE_{y,0}\|^2)}{T(1-\rho)}
+ \frac{1}{2T}
\sum_{i=1}^{T}
(\mE\|\mhE_{x,i}\|^2+\mE\|\mhE_{y,i}\|^2) 
+\frac{144L^2_f\lambda^2_a\mu^2_y}{(1-\rho)^2\underline{\lambda^2_b}T}
\sum_{i=0}^{T-1}(\mu^2_x\mE\|\bm^x_{c,i}\|^2\notag \\
&\quad +\mu^2_y\mE\|\bm^y_{c,i}\|^2)
+ \frac{\mu^2_y\lambda^2_a(6p+6\beta^2)}{(1-\rho)^2K\underline{\lambda^2_b}}
\Bigg(
\frac{2}{T\bar{\beta}}
(\mE\|\mS_{x,0}\|^2
+
\mE\|\mS_{y,0}\|^2)+\frac{12KL^2_fv^2_1v^2_2}{bT\bar{\beta}}
(\mE\|\mhE_{x,0}\|^2 
\notag \\
&\quad  +\mE\|\mhE_{y,0}\|^2) 
\Bigg)+ \frac{48\mu^2_y\lambda^2_a}{(1-\rho)^2\underline{\lambda^2_b}}
\Big(
\frac{p\sigma^2}{B} \mathbb{I}_{\text{online}}
+ \frac{\beta^2\sigma^2}{b}
\Big).
\end{align}
It follows that
\begin{align}
&\frac{1}{T}
\sum_{i=1}^{T}
\mE\|\mhE_{x,i}\|^2
+\mE\|\mhE_{y,i}\|^2 \notag \\
&\le 
\frac{4(\mE\|\mhE_{x,0}\|^2 +\mE\|\mhE_{y,0}\|^2)}{T(1-\rho)} 
+\frac{288L^2_f\lambda^2_a\mu^2_y}{(1-\rho)^2\underline{\lambda^2_b}T}
\sum_{i=0}^{T-1}(\mu^2_x\mE\|\bm^x_{c,i}\|^2+\mu^2_y\mE\|\bm^y_{c,i}\|^2) + \frac{12\mu^2_y\lambda^2_a(p+\beta^2)}{(1-\rho)^2K\underline{\lambda^2_b}}\notag \\
&\quad 
\times
\Bigg(
\frac{2}{T\bar{\beta}}
(\mE\|\mS_{x,0}\|^2
+
\mE\|\mS_{y,0}\|^2)
+\frac{12KL^2_fv^2_1v^2_2}{bT\bar{\beta}}
(\mE\|\mhE_{x,0}\|^2 +\mE\|\mhE_{y,0}\|^2) 
\Bigg)
\notag \\
&\quad 
+ \frac{96\mu^2_y\lambda^2_a}{(1-\rho)^2\underline{\lambda^2_b}}
\Big(
\frac{p\sigma^2}{B} \mathbb{I}_{\text{online}}
+ \frac{\beta^2\sigma^2}{b}
\Big).
\end{align}
Adding $\frac{\mE\|\mhE_{x,0}\|^2 +\mE\|\mhE_{y,0}\|^2}{T(1-\rho)}$ into the above relation and using  $\frac{1}{(1-\rho)T} \ge \frac{1}{T}$, we obtain 
\begin{align}
&\frac{1}{T}
\sum_{i=0}^{T-1}
(\mE\|\mhE_{x,i}\|^2 +\mE\|\mhE_{y,i}\|^2)
\notag \\
&\le 
\frac{1}{T}
\sum_{i=1}^{T}
(\mE\|\mhE_{x,i}\|^2 +\mE\|\mhE_{y,i}\|^2)
+ \frac{\mE\|\mhE_{x,0}\|^2 +\mE\|\mhE_{y,0}\|^2}{T(1-\rho)} \notag \\
&\le 
\frac{5(\mE\|\mhE_{x,0}\|^2 +\mE\|\mhE_{y,0}\|^2)}{T(1-\rho)} 
+\frac{288L^2_f\lambda^2_a\mu^2_y}{(1-\rho)^2\underline{\lambda^2_b}T}
\sum_{i=0}^{T-1}(\mu^2_x\mE\|\bm^x_{c,i}\|^2+\mu^2_y\mE\|\bm^y_{c,i}\|^2) +\frac{12\mu^2_y\lambda^2_a(p+\beta^2)}{(1-\rho)^2K\underline{\lambda^2_b}}\notag \\
&\quad 
\times 
\Bigg(
\frac{2}{T\bar{\beta}}
(\mE\|\mS_{x,0}\|^2
+
\mE\|\mS_{y,0}\|^2)
+\frac{12KL^2_fv^2_1v^2_2}{bT\bar{\beta}}
(\mE\|\mhE_{x,0}\|^2 +\mE\|\mhE_{y,0}\|^2) 
\Bigg)
\notag \\
&\quad 
+ \frac{96\mu^2_y\lambda^2_a}{(1-\rho)^2\underline{\lambda^2_b}}
\Big(
\frac{p\sigma^2}{B} \mathbb{I}_{\text{online}}
+ \frac{\beta^2\sigma^2}{b}
\Big).
\end{align}

\subsection{Proof of Lemma \ref{lemma:optimality_gap}}
\label{subsec:optimality_gap}

Recall that the following optimality gap at the network centroid
\begin{align}
    \Delta_{c, i} \triangleq P(\bx_{c,i}) - J(\bx_{c,i},\by_{c,i}).
\end{align}
Since 
$-J(x, \cdot)$ is $L_f$-smooth, we have 
\begin{align}
&-J(\bx_{c,i+1}, \by_{c,i+1}) \notag \\
&\le 
-J(\bx_{c,i+1}, \by_{c,i})
- \mu_y\langle 
\nabla_y J(\bx_{c,i+1}, \by_{c,i}),
\bm^y_{c,i}\rangle
+ \frac{L_f\mu^2_y}{2}
\|\bm^y_{c,i}\|^2
\notag \\
&=
-J(\bx_{c,i+1}, \by_{c,i})
-
\frac{\mu_y}{2}\|\nabla_y J(\bx_{c,i+1}, \by_{c,i})\|^2
-
\frac{\mu_y}{2}
\|\bm^y_{c,i}
\|^2
+\frac{\mu_y}{2}
\|\nabla_y J(\bx_{c,i+1}, \by_{c,i})
-\bm^y_{c,i}\|^2 \notag \\
&\quad 
+ 
\frac{L_f\mu^2_y}{2}
\|\bm^y_{c,i}\|^2 \notag \\
&\overset{(a)}{\le}  
-J(\bx_{c,i+1}, \by_{c,i})
-\mu_y\nu[P(\bx_{c,i+1}) - J(\bx_{c,i+1}, \by_{c,i})]
-\frac{\mu_y}{2}
\|\bm^y_{c,i}\|^2
+\mu_y
\|\nabla_y J(\bx_{c,i+1},\by_{c,i}) \notag \\
&\quad -\nabla_y J(\bx_{c,i}, \by_{c,i})\|^2 + 
\mu_y\|\nabla_y J(\bx_{c,i}, \by_{c,i}) - \bm^y_{c,i}\|^2
+ \frac{L_f\mu^2_y}{2}
\|\bm^y_{c,i}\|^2,
\end{align}
where the $(a)$ follows from the $\nu$-PL condition
\(
\|\nabla_y J(\bx_{c,i+1}, \by_{c,i})\|^2 \ge 2\nu (P(\bx_{c,i+1}) - J(\bx_{c,i+1}, \by_{c,i}))
\) and Jensen's inequality.
Adding \(P(\bx_{c,i+1})\)
into the above expression and taking the expectation, we get 
\begin{align}
&\mE\Delta_{c,i+1}
\notag\\
&\overset{(a)}{\le} \mE\Big[ 
(1-\mu_y\nu)[P(\bx_{c,i+1}) - J(\bx_{c,i+1}, \by_{c,i})]
-\frac{\mu_y}{2}
\|\bm^y_{c,i}\|^2
+ \mu_y\mu^2_xL^2_f\|\bm^x_{c,i}\|^2 +2\mu_y\Big\|\nabla_y J(\bx_{c,i}, 
 \by_{c,i})\notag \\
&\quad 
 - \frac{1}{K}
 \sum_{k=1}^K \nabla_y J_k(\bx_{k,i}, \by_{k,i})\Big\|^2
+2\mu_y\Big\|\frac{1}{K}\sum_{k=1}^K \nabla_y J_k(\bx_{k,i}, \by_{k,i}) - \bm^y_{c,i}\Big\|^2
 +\frac{L_f\mu^2_y}{2}
 \|\bm^y_{c,i}\|^2 \Big] \notag \\
 &
 \overset{(b)}{\le} \mE\Big[
(1-\mu_y\nu)[P(\bx_{c,i+1}) - P(\bx_{c,i})+ P(\bx_{c,i})
-J(\bx_{c,i}, \by_{c,i})
+ J(\bx_{c,i}, \by_{c,i})
-  J(\bx_{c,i+1}, \by_{c,i})]
\notag \\
&\quad + \mu_y\mu^2_xL^2_f\|\bm^x_{c,i}\|^2
+2\mu_y\Big\|\nabla_y J(\bx_{c,i}, 
 \by_{c,i}) - \frac{1}{K}
 \sum_{k=1}^K \nabla_y J_k(\bx_{k,i}, \by_{k,i})\Big\|^2
 \notag \\
 &\quad+2\mu_y\Big\|\frac{1}{K}\sum_{k=1}^K \nabla_y J_k(\bx_{k,i}, \by_{k,i}) - \bm^y_{c,i}\Big\|^2 -\frac{\mu_y}{2}(1-L_f\mu_y)
 \|\bm^y_{c,i}\|^2 \Big]
 \notag \\
 &\overset{(c)}{\le} \mE\Big[
(1-\mu_y\nu)[P(\bx_{c,i+1}) - P(\bx_{c,i})+ P(\bx_{c,i})
-J(\bx_{c,i}, \by_{c,i})
+ J(\bx_{c,i}, \by_{c,i})
-  J(\bx_{c,i+1}, \by_{c,i})]
\notag \\
&\quad 
+ \mu_y\mu^2_xL^2_f\|\bm^x_{c,i}\|^2+\frac{2\mu_y L^2_f}{K}
(\|\mX_i -\mX_{c,i}\|^2 + \|\mY_i -\mY_{c,i}\|^2)
 +2\mu_y\Big\|\frac{1}{K}\sum_{k=1}^K \nabla_y J_k(\bx_{k,i}, \by_{k,i}) - \bm^y_{c,i}\Big\|^2
\notag \\
&\quad 
-\frac{\mu_y}{2}(1-L_f\mu_y)
 \|\bm^y_{c,i}\|^2 \Big],
\end{align}
where $(a)$ and $(c)$ are  due to the expected $L_f$-smooth assumption and Jensen's inequality;
$(b)$ follows by adding and subtracting 
$P(\bx_{c,i}), J(\bx_{c,i}, \by_{c,i})$.
Using the $L_f$-smoothness property for $-J(\cdot,y)$
and $L$-smooth property 
for 
$P(\cdot)$,
we get 
\begin{align}
&P(\bx_{c,i+1})  - P(\bx_{c,i})+J(\bx_{c,i}, \by_{c,i}) - J(\bx_{c,i+1}, \by_{c,i})
\notag \\
&\le   \langle \nabla P(\bx_{c,i}) - \nabla_x J(\bx_{c,i}, \by_{c,i}), \bx_{c,i+1} - \bx_{c,i}\rangle +\frac{(L+L_f)\mu^2_x}{2}\|\bm^x_{c,i}\|^2.
\end{align}
To bound the cross-term, we use the inequality $\langle a, b\rangle \le \frac{1}{2t}\|a\|^2 + \frac{t}{2}\|b\|^2$  and choose the parameter \(t = \frac{1}{8\mu_x}\) to get
\begin{align}
&\langle 
\nabla P(\bx_{c,i}) - \nabla_x J(\bx_{c,i}, \by_{c,i}), \bx_{c,i+1} -\bx_{c,i} \rangle \notag \\
&\le \frac{1}{16\mu_x}
\|\bx_{c,i+1} - \bx_{c,i}\|^2
+4\mu_x
\|\nabla P(\bx_{c,i}) - \nabla_x J(\bx_{c,i}, \by_{c,i})\|^2
\notag \\
&\le 
\frac{\mu_x}{16}
\|\bm^x_{c,i}\|^2
+4\mu_xL^2_f 
\|\by_{c,i} - \by^o(\bx_{c,i})\|^2
\notag \\
&\le 
\frac{\mu_x}{16}
\|\bm^x_{c,i}\|^2
+ 8 L_f\kappa\mu_x\Delta_{c,i},
\end{align}
where the last inequality follows from the quadratic growth property of the $\nu$-PL function
\cite{karimi2016linear}
\begin{align}
\Delta_{c,i} = P(\bx_{c,i})- J(\bx_{c,i}, \by_{c,i}) \ge \frac{\nu}{2}
\|\by_{c,i}-\by^o(\bx_{c,i})\|^2.
\end{align}
Therefore, we have
\begin{align}
&P(\bx_{c,i+1}) -P(\bx_{c,i})+J(\bx_{c,i}, \by_{c,i}) - J(\bx_{c,i+1}, \by_{c,i}) \notag \\
&\le 
8L_f\kappa\mu_x\Delta_{c,i}
+\frac{\mu_x}{16}
\|\bm^x_{c,i}\|^2
+\frac{(L+L_f)\mu^2_x}{2}
\|\bm^x_{c,i}\|^2.
\end{align}
Putting the results together,
we have
\begin{align}
&\Delta_{c,i+1} \notag \\
&\le 
(1-\mu_y\nu) \Big( \Delta_{c,i}+
8L_f\kappa\mu_x\Delta_{c,i}
+\frac{\mu_x}{16}
\|\bm^x_{c,i}\|^2
+\frac{(L+L_f)\mu^2_x}{2}
\|\bm^x_{c,i}\|^2\Big)
+\mu_y\mu^2_xL^2_f\|\bm^x_{c,i}\|^2 \notag \\
&\quad 
+\frac{2\mu_yL^2_f}{K}(\|\mX_{i} -\mX_{c,i}\|^2 +\|\mY_{i} - \mY_{c,i}\|^2)
 +2\mu_y\Big\|\frac{1}{K}\sum_{k=1}^K \nabla_y J_k(\bx_{k,i}, \by_{k,i}) - \bm^y_{c,i}\Big\|^2
\notag \\
&\quad 
-\frac{\mu_y}{2}
(1-L_f\mu_y)
\|\bm^y_{c,i}\|^2
.\end{align}
Invoking Lemma \ref{appendix:lemma:Lconsense} 
\begin{align}
\label{proof:consensus_to_e}
&\|\mX_i - \mX_{c,i}\|^2
+\|\mY_i -\mY_{c,i}\|^2 
\le
Kv^2_1v^2_2
(\|\mhE_{x,i}\|^2 +\|\mhE_{y,i}\|^2),
\end{align}
we  get 
\begin{align}
&\Delta_{c,i+1} \notag \\
&\le 
(1-\mu_y\nu) \Big( 
(1+8L_f\kappa\mu_x)\Delta_{c,i}
+\frac{\mu_x}{16}
\|\bm^x_{c,i}\|^2
+\frac{(L+L_f)\mu^2_x}{2}
\|\bm^x_{c,i}\|^2\Big)
+\mu_y\mu^2_xL^2_f\|\bm^x_{c,i}\|^2 \notag \\
&\quad 
+2\mu_yL^2_fv^2_1v^2_2(\|\mhE_{x,i}\|^2 +\|\mhE_{y,i}\|^2)
 +2\mu_y\Big\|\frac{1}{K}\sum_{k=1}^K \nabla_y J_k(\bx_{k,i}, \by_{k,i}) - \bm^y_{c,i}\Big\|^2
\notag \\
&\quad 
-\frac{\mu_y}{2}
(1-L_f\mu_y)
\|\bm^y_{c,i}\|^2.
\end{align}
To establish a descent relation,
we choose the following step size conditions:
\begin{subequations}
\begin{align}
1-\mu_y \nu& \ge 0 \Longrightarrow \mu_y \le \frac{1}{\nu},  \quad &
8\kappa L_f \mu_x &\le \frac{\nu\mu_y}{2} 
\Longrightarrow \mu_x \le 
\frac{\mu_y}{16\kappa^2} , \\
1-L_f\mu_y &\ge \frac{1}{2} \Longrightarrow  \mu_y \le \frac{1}{2L_f}, \quad &
2L\mu^2_x &\le \frac{\mu_x}{16}
\Longleftarrow
\mu_x \le \frac{1}{32L}.
\end{align}
\end{subequations}
We then have 
\begin{subequations}
\begin{align}
(1-\mu_y\nu)(1+8L_f\kappa\mu_x)
&\le (1-\mu_y\nu)(1+ \frac{\mu_y\nu}{2})
\le 1- \frac{\mu_y\nu}{2},
\\
-\frac{\mu_y}{2}(1-L_f\mu_y) 
&\le -\frac{\mu_y}{4}, \\
\frac{L+L_f}{2} &\le  L, \\
\frac{\mu_x}{16}
+ \frac{(L+L_f)\mu^2_x}{2} + \mu_y\mu^2_xL^2_f  &\le \frac{\mu_x}{16} 
+ L\mu^2_x
+ \frac{\mu^2_x L_f}{2} \le  \frac{\mu_x}{16} + 2L\mu^2_x \le \frac{\mu_x}{8}.
\end{align}
\end{subequations}
Recalling the notation 
$\bs^y_{c,i}$, we can conclude that
\begin{align}
\Delta_{c,i+1}
\le 
(1-\frac{\nu\mu_y}{2})
\Delta_{c,i}
+\frac{\mu_x}{8}
\|\bm^x_{c,i}\|^2-\frac{\mu_y}{4}
\|\bm^y_{c,i}\|^2
+2\mu_yL^2_fv^2_1v^2_2(\|\mhE_{x,i}\|^2 +\|\mhE_{y,i}\|^2)
+2\mu_y\|\bs^y_{c,i}\|^2.
\end{align}
Applying the above inequality iteratively, it holds for any $i =1, 2, \dots$
\begin{align}
\Delta_{c,i}
&\le 
\Big(
1- \frac{\nu\mu_y}{2}\Big)^{i}
\Delta_{c,0}
+
\frac{\mu_x}{8}
\sum_{j=0}^{i-1}
\Big(
1- \frac{\nu\mu_y}{2}\Big)^{i-j-1}
\|\bm^x_{c,j}\|^2
-\frac{\mu_y}{4}
\sum_{j=0}^{i-1}\Big(
1- \frac{\nu\mu_y}{2}\Big)^{i-j-1}
\|\bm^y_{c,j}\|^2
\notag \\
&\quad
+2\mu_yL^2_fv^2_1v^2_2
\sum_{j=0}^{i-1}
\Big(
1- \frac{\nu\mu_y}{2}\Big)^{i-j-1}
(\|\mhE_{x,j}\|^2+
\|\mhE_{y,j}\|^2)+
2\mu_y
\sum_{j=0}^{i-1}
\Big(
1- \frac{\nu\mu_y}{2}\Big)^{i-j-1}
\|\bs^y_{c,j}\|^2. 
\end{align}
Averaging above inequality over $i= 1, \dots, T$, we get 
\begin{align}
&\frac{1}{T}\sum_{i=1}^{T}\Delta_{c,i} \notag \\
&\le 
\frac{1}{T}\sum_{i=1}^{T}\Big(
1- \frac{\nu\mu_y}{2}\Big)^{i}
\Delta_{c,0}
+
\frac{\mu_x}{8}\frac{1}{T}\sum_{i=1}^{T}
\sum_{j=0}^{i-1}
\Big(
1- \frac{\nu\mu_y}{2}\Big)^{i-j-1}
\|\bm^x_{c,j}\|^2 \notag \\
&\quad 
-\frac{\mu_y}{4}\frac{1}{T}\sum_{i=1}^{T}
\sum_{j=0}^{i-1}\Big(
1- \frac{\nu\mu_y}{2}\Big)^{i-j-1}
\|\bm^y_{c,j}\|^2+
2\mu_yL^2_fv^2_1v^2_2\frac{1}{T}\sum_{i=1}^{T}
\sum_{j=0}^{i-1}
\Big(
1- \frac{\nu\mu_y}{2}\Big)^{i-j-1}
(\|\mhE_{x,j}\|^2\notag \\
&\quad+
\|\mhE_{y,j}\|^2)
+
\frac{2\mu_y}{T}\sum_{i=1}^{T}
\sum_{j=0}^{i-1}
\Big(
1- \frac{\nu\mu_y}{2}\Big)^{i-j-1}
\|\bs^y_{c,j}\|^2 \notag \\
&\le 
\frac{2}{T\nu\mu_y}
\Delta_{c,0}
+
\frac{\mu_x}{4\nu\mu_y T}
\sum_{i=0}^{T-1}
\|\bm^x_{c,i}\|^2
-\frac{\mu_y}{4T}\sum_{i=1}^{T}
\sum_{j=0}^{i-1}\Big(
1- \frac{\nu\mu_y}{2}\Big)^{i-j-1}
\|\bm^y_{c,j}\|^2 \notag \\
&\quad +
\frac{4\kappa L_fv^2_1v^2_2}{T}
\sum_{i=0}^{T-1}
(\|\mhE_{x,i}\|^2
+\|\mhE_{y,i}\|^2)
+\frac{4}{T\nu}
\sum_{i=0}^{T-1}
\|\bs^y_{c,i}\|^2.
\end{align}
Adding \(\frac{\Delta_{c,0}}{T\nu\mu_y}\) on both sides of the above relation and choosing $\mu_y\nu \le 1 \Longrightarrow \mu_y \le \frac{1}{\nu}$, we can deduce that 
\begin{align}
\frac{1}{T}
\sum_{i=0}^{T-1}
\Delta_{c,i}
&\le 
\frac{3}{T\nu\mu_y}
\Delta_{c,0}
+
\frac{\mu_x}{4\nu\mu_y T}
\sum_{i=0}^{T-1}
\|\bm^x_{c,i}\|^2
-\frac{\mu_y}{4T}\sum_{i=1}^{T}
\sum_{j=0}^{i-1}\Big(
1- \frac{\nu\mu_y}{2}\Big)^{i-j-1}
\|\bm^y_{c,j}\|^2  \notag \\
&\quad+
\frac{4\kappa L_fv^2_1v^2_2}{T}
\sum_{i=0}^{T-1}
(\|\mhE_{x,i}\|^2
+\|\mhE_{y,i}\|^2)
+\frac{4}{T\nu}
\sum_{i=0}^{T-1}
\|\bs^y_{c,i}\|^2.
\end{align}

\subsection{Proof of Lemma \ref{appendix:lemma:costfunction}}
\label{appendix:proof:costfunction}
 Since \(P(\cdot)\) is \(L\)-smooth \cite{nouiehed2019solving}, we have
\begin{align}
&P(\bx_{c,i+1}) \notag \\
&\le 
P(\bx_{c,i})
- \mu_x \langle 
\nabla P(\bx_{c,i}), \bm^x_{c,i}
\rangle + \frac{L\mu^2_x}{2}\|\bm^x_{c,i}\|^2 \notag \\
&\le 
P(\bx_{c,i})
-\frac{\mu_x}{2}
\|\nabla P(\bx_{c,i})\|^2
-\frac{\mu_x}{2}
\|\bm^x_{c,i}\|^2
+\frac{\mu_x}{2}
\|\nabla P(\bx_{c,i}) - \bm^x_{c,i}\|^2
+\frac{L\mu^2_x}{2}
\|\bm^x_{c,i}\|^2
\notag \\
&\overset{(a)}{\le} 
P(\bx_{c,i})
-\frac{\mu_x}{2}
\|\nabla P(\bx_{c,i})\|^2
-\frac{\mu_x}{2} (1-L\mu_x)
\|\bm^x_{c,i}\|^2  +
\mu_x
\Big\|\nabla P(\bx_{c,i}) - \frac{1}{K}\sum_{k=1}^K \nabla_x J_k(\bx_{k,i},\by_{k,i})\Big\|^2
\notag \\
&\quad +\mu_x \Big\|\frac{1}{K}\sum_{k=1}^{K}
\bm^x_{k,i} - \frac{1}{K}\sum_{k=1}^{K}\nabla_x J_k(\bx_{k,i}, \by_{k,i})\Big\|^2
\notag\\
&\overset{(b)}{\le} 
P(\bx_{c,i})
-\frac{\mu_x}{2}
\|\nabla P(\bx_{c,i})\|^2
-\frac{\mu_x}{2}
(1-L\mu_x)
\|\bm^x_{c,i}\|^2
+
\frac{10\mu_x}{9}\Big\|\nabla P(\bx_{c,i}) - \frac{1}{K}
\sum_{k=1}^K\nabla_x J_k(\bx_{c,i},\by_{c,i})\Big\|^2
\notag \\
&\quad 
+10\mu_x\Big\|\frac{1}{K}
\sum_{k=1}^K (\nabla_x J_k(\bx_{c,i},\by_{c,i}) - \nabla_x J_k(\bx_{k,i},\by_{k,i}))\Big\|^2 +\mu_x\|\bs^x_{c,i}\|^2\notag \\
&\overset{(c)}{\le} 
P(\bx_{c,i})
-\frac{\mu_x}{2}
\|\nabla P(\bx_{c,i})\|^2
-\frac{\mu_x}{2}
(1-L\mu_x)
\|\bm^x_{c,i}\|^2
+\frac{10}{9}\mu_xL^2_f\|\by_{c,i}-\by^o(\bx_{c,i})\|^2
\notag \\
&\quad +\frac{10\mu_xL^2_f}{K}(\|\mX_i-\mX_{c,i}\|^2+\|\mY_i -\mY_{c,i}\|^2) + \mu_x\|\bs^x_{c,i}\|^2\notag \\
&\overset{(d)}{\le} P(\bx_{c,i})
-\frac{\mu_x}{2}
\|\nabla P(\bx_{c,i})\|^2
-\frac{\mu_x}{2}
(1-L\mu_x)
\|\bm^x_{c,i}\|^2
+\frac{10}{9}\mu_xL^2_f\|\by_{c,i} - \by^o(\bx_{c,i})\|^2 \notag \\
&\quad 
+ 10\mu_xL^2_f v^2_1 v^2_2 
(\|\mhE_{x,i}\|^2 + \|\mhE_{y,i}\|^2) + \mu_x\|\bs^x_{c,i}\|^2 \notag \\
&\overset{(e)}{\le} P(\bx_{c,i})
-\frac{\mu_x}{2}
\|\nabla P(\bx_{c,i})\|^2
-\frac{\mu_x}{2}
(1-L\mu_x)
\|\bm^x_{c,i}\|^2
+\frac{20}{9}\mu_x \kappa L_f \Delta_{c,i} \notag \\
&\quad 
+ 10\mu_xL^2_f v^2_1 v^2_2 
(\|\mhE_{x,i}\|^2 + \|\mhE_{y,i}\|^2) + \mu_x\|\bs^x_{c,i}\|^2,
\end{align}
where $(a)$ and $(b)$ follow from Jensen's inequality, i.e., $\|a+b\|^2 \le \frac{\|a\|^2}{t} +\frac{\|b\|^2}{1-t}$ for some  $t \in (0, 1)$;
$(c)$ follows from $L_f$-smooth assumption;$(d)$ follows from 
Lemma \ref{appendix:lemma:Lconsense}; $(e)$ follows from the quadratic growth of the $\nu$-PL function \cite{karimi2016linear}.

\subsection{Proof of the main result}
\label{appendix:proofofmain}
\begin{proof}
Invoking Lemma 
\eqref{lemma:costvalue},
we deduce that 
\begin{align}
&\|\nabla P(\bx_{c,i})\|^2
\notag\\
&\le 
\frac{2(P(\bx_{c,i})-P(\bx_{c,i+1}))}{\mu_x}
-(1-L\mu_x)
\|\bm^x_{c,i}\|^2
+\frac{40}{9}\kappa L_f \Delta_{c,i}
+20v^2_1v^2_2L^2_f
(\|\mhE_{x,i}\|^2 +\|\mhE_{y,i}\|^2)
\notag\\
&\quad
+2\|\bs^x_{c,i}\|^2.
\end{align}
Taking expectation and 
adding 
$1.5\kappa L_f \mE\Delta_{c,i}$
on both sides of the above inequality, we get 
\begin{align}
&\mE\|\nabla P(\bx_{c,i})\|^2
+1.5\kappa L_f\mE\Delta_{c,i} \notag\\
&\le 
\frac{2\mE(P(\bx_{c,i})-P(\bx_{c,i+1}))}{\mu_x}
-(1-L\mu_x)
\mE\|\bm^x_{c,i}\|^2
+6\kappa L_f \mE\Delta_{c,i}
+20v^2_1v^2_2L^2_f
(\mE\|\mhE_{x,i}\|^2 \notag \\
&\quad 
+\mE\|\mhE_{y,i}\|^2)
+2\mE\|\bs^x_{c,i}\|^2.
\end{align}
Averaging the above 
inequality
for $i = 0, \dots, T-1$
and telescoping the first term, we get
\begin{align}
&\frac{1}{T}
\sum_{i=0}^{T-1}
(\mE\|\nabla P(\bx_{c,i})\|^2
+1.5\kappa L_f\mE\Delta_{c,i}) \notag \\
&\le 
\frac{2\mE(P(\bx_{c,0})-P^\star)}{T\mu_x}
-\frac{(1-L\mu_x)}{T}
\sum_{i=0}^{T-1}
\mE\|\bm^x_{c,i}\|^2
+ \frac{6\kappa L_f}{T}
\sum_{i=0}^{T-1}\mE\Delta_{c,i}+\frac{20v^2_1v^2_2L^2_f}{T}
\sum_{i=0}^{T-1}
(\mE\|\mhE_{x,i}\|^2
\notag \\
&\quad 
 +\mE\|\mhE_{y,i}\|^2)+\frac{2}{T}
\sum_{i=0}^{T-1}\mE\|\bs^x_{c,i}\|^2.
\end{align}
Invoking Lemma \ref{lemma:optimality_gap}, we have 
\begin{align}
&\frac{1}{T}
\sum_{i=0}^{T-1}
(\mE\|\nabla P(\bx_{c,i})\|^2
+1.5\kappa L_f\mE\Delta_{c,i}) \notag \\
&\le 
\frac{2\mE(P(\bx_{c,0})-P^\star)}{T\mu_x}
-\frac{(1-L\mu_x)}{T}
\sum_{i=0}^{T-1}
\mE\|\bm^x_{c,i}\|^2
+
\Bigg( 
\frac{18\kappa^2}{T\mu_y}
\mE\Delta_{c,0}
+ \frac{3\kappa^2\mu_x}{2\mu_y T}
\sum_{i=0}^{T-1}
\mE\|\bm^x_{c,i}\|^2
\notag \\
&\quad 
-\frac{3\kappa L_f \mu_y}{2T}
\sum_{i=1}^{T}
\sum_{j=0}^{i-1}\Big(
1- \frac{\nu\mu_y}{2}\Big)^{i-j-1}
\mE\|\bm^y_{c,j}\|^2 
+ \frac{24\kappa^2L^2_fv^2_1v^2_2}{T}\sum_{i=0}^{T-1}
(\mE\|\mhE_{x,i}\|^2
+\mE\|\mhE_{y,i}\|^2)
\notag\\
&\quad 
+\frac{24\kappa^2}{T}\sum_{i=0}^{T-1}
\mE\|\bs^y_{c,i}\|^2
\Bigg)+\frac{20v^2_1v^2_2L^2_f}{T}
\sum_{i=0}^{T-1}
(\mE\|\mhE_{x,i}\|^2 +\mE\|\mhE_{y,i}\|^2)
+\frac{2}{T}
\sum_{i=0}^{T-1}\mE\|\bs^x_{c,i}\|^2.
\end{align}
Using the fact that $20 < 20\kappa^2$ and $2 < 24\kappa^2$,  the last few terms can be absorbed as follows:
\begin{align}
&\frac{1}{T}
\sum_{i=0}^{T-1}
(\mE\|\nabla P(\bx_{c,i})\|^2
+1.5\kappa L_f\mE\Delta_{c,i}) \notag \\
&\le 
\frac{2\mE(P(\bx_{c,0}-P^\star))}{T\mu_x}
-\frac{(1-L\mu_x)}{T}
\sum_{i=0}^{T-1}
\mE\|\bm^x_{c,i}\|^2
+
\Bigg( 
\frac{18\kappa^2}{T\mu_y}
\mE\Delta_{c,0}
+ \frac{3\kappa^2\mu_x}{2\mu_yT}
\sum_{i=0}^{T-1}
\mE\|\bm^x_{c,i}\|^2 
\notag \\
&\quad 
-\frac{3\kappa L_f \mu_y}{2T}
\sum_{i=1}^{T}
\sum_{j=0}^{i-1}\Big(
1- \frac{\nu\mu_y}{2}\Big)^{i-j-1}
\mE\|\bm^y_{c,j}\|^2  \Bigg)+\frac{44\kappa^2L^2_fv^2_1v^2_2}{T}
\sum_{i=0}^{T-1}
(\mE\|\mhE_{x,i}\|^2 +\mE\|\mhE_{y,i}\|^2)
\notag\\
&\quad 
+\frac{24\kappa^2}{T}
\sum_{i=0}^{T-1}(\mE\|\bs^x_{c,i}\|^2+
\mE\|\bs^y_{c,i}\|^2).
\end{align}
Invoking Lemma
\ref{lemma:gradienterror_average}, we get 
\begin{align}
&\frac{1}{T}
\sum_{i=0}^{T-1}
(\mE\|\nabla P(\bx_{c,i})\|^2
+1.5\kappa L_f\mE\Delta_{c,i}) \notag \\
&\le 
\frac{2\mE(P(\bx_{c,0})-P^\star)}{T\mu_x}
-\frac{(1-L\mu_x)}{T}
\sum_{i=0}^{T-1}
\mE\|\bm^x_{c,i}\|^2
+
\Bigg( 
\frac{18\kappa^2}{T\mu_y}
\mE\Delta_{c,0}
+ \frac{3\kappa^2\mu_x}{2\mu_yT}
\sum_{i=0}^{T-1}
\mE\|\bm^x_{c,i}\|^2
\notag \\
&\quad 
-\frac{3\kappa L_f \mu_y}{2T}
\sum_{i=1}^{T}
\sum_{j=0}^{i-1}\Big(
1- \frac{\nu\mu_y}{2}\Big)^{i-j-1}
\mE\|\bm^y_{c,j}\|^2  \Bigg)+\frac{44\kappa^2L^2_fv^2_1v^2_2}{T}
\sum_{i=0}^{T-1}
(\mE\|\mhE_{x,i}\|^2 +\mE\|\mhE_{y,i}\|^2) \notag \\
&\quad +
\Bigg( 
\frac{48\kappa^2}{T\bar{\beta}}
(\mE\|\bs^x_{c,0}\|^2
\notag\\
&\quad +
\mE\|\bs^y_{c,0}\|^2)
+\frac{576\kappa^2L^2_fv^2_1v^2_2}{bKT\bar{\beta}}\sum_{i=1}^{T}
(\mE\|\mhE_{x,i}\|^2 +\mE\|\mhE_{y,i}\|^2) +\frac{288\kappa^2 L^2_fv^2_1v^2_2}{bKT\bar{\beta}} (\mE\|\mhE_{x,0}\|^2+\mE\|\mhE_{y,0}\|^2)
\notag \\
&\quad 
+
\frac{288\kappa^2 L^2_f}{bKT\bar{\beta}}
\sum_{i=0}^{T-1}
\mu^2_x\mE\|\bm^x_{c,i}\|^2
+\frac{288\kappa^2 L^2_f}{bKT\bar{\beta}}
\sum_{i=0}^{T-1}\mu^2_{y}(1-(1-\bar{\beta})^{T-i})\mE\|\bm^y_{c,i}\|^2+ \frac{96\kappa^2\beta^2\sigma^2}{Kb\bar{\beta}} \notag \\
&\quad 
+\frac{48\kappa^2p\sigma^2}{KB\bar{\beta}}\mathbb{I}_{\text{online}}
\Bigg).
\end{align}
Regrouping terms,
we get 
\begin{align}
&\frac{1}{T}
\sum_{i=0}^{T-1}
(\mE\|\nabla P(\bx_{c,i})\|^2
+1.5\kappa L_f\mE\Delta_{c,i}) \notag \\
&\le 
\frac{2\mE(P(\bx_{c,0})-P^\star)}{T\mu_x}-\Bigg(1-L\mu_x - \frac{3\kappa^2\mu_x}{2\mu_y}- \frac{288\kappa^2L^2_f\mu^2_x}{bK\bar{\beta}}\Bigg) \frac{1}{T}
\sum_{i=0}^{T-1}
\mE\|\bm^x_{c,i}\|^2
\notag \\
&\quad 
+
\Bigg( \frac{288\kappa^2 L^2_f}{bKT\bar{\beta}}
\sum_{i=0}^{T-1}\mu^2_{y}(1-(1-\bar{\beta})^{T-i})\mE\|\bm^y_{c,i}\|^2
-\frac{3\kappa L_f \mu_y}{2T}
\sum_{i=1}^{T}
\sum_{j=0}^{i-1}\Big(
1- \frac{\nu\mu_y}{2}\Big)^{i-j-1}
\mE\|\bm^y_{c,j}\|^2  \Bigg) \notag \\
&\quad
+\frac{576\kappa^2L^2_fv^2_1v^2_2}{bK\bar{\beta}T}
\sum_{i=1}^{T}
(\mE\|\mhE_{x,i}\|^2 +\mE\|\mhE_{y,i}\|^2)+
\frac{44\kappa^2L^2_fv^2_1v^2_2}{T}
\sum_{i=0}^{T-1}
(\mE\|\mhE_{x,i}\|^2 +\mE\|\mhE_{y,i}\|^2)
\notag\\
&\quad 
+
\Bigg(  \frac{18\kappa^2}{T\mu_y}
\mE\Delta_{c,0} +
\frac{48\kappa^2}{T\bar{\beta}}
(\mE\|\bs^x_{c,0}\|^2+
\mE\|\bs^y_{c,0}\|^2)
+\frac{288\kappa^2 L^2_fv^2_1v^2_2}{bKT\bar{\beta}} (\mE\|\mhE_{x,0}\|^2+\mE\|\mhE_{y,0}\|^2) 
\notag \\
&\quad 
+ \frac{96\kappa^2\beta^2\sigma^2}{Kb\bar{\beta}}+\frac{48\kappa^2p\sigma^2}{KB\bar{\beta}}\mathbb{I}_{\text{online}}
\Bigg).
\end{align}
Choosing 
\begin{align}
44\kappa^2L^2_fv^2_1v^2_2 \le \frac{44\kappa^2L^2_fv^2_1v^2_2}{bK\bar{\beta}} \Longrightarrow 
 b\bar{\beta} \le \frac{1}{K}.
\end{align}
Invoking Lemma \ref{lemma:coupled_error}, we get 
\begin{align}
&\frac{1}{T}
\sum_{i=0}^{T-1}
(\mE\|\nabla P(\bx_{c,i})\|^2
+1.5\kappa L_f\mE\Delta_{c,i}) \notag \\
&\le 
\frac{2\mE(P(\bx_{c,0})-P^\star)}{T\mu_x}-\Bigg(1-L\mu_x - \frac{3\kappa^2\mu_x}{2\mu_y}- \frac{288\kappa^2L^2_f\mu^2_x}{bK\bar{\beta}} \Bigg) \frac{1}{T}
\sum_{i=0}^{T-1}
\mE\|\bm^x_{c,i}\|^2
\notag \\
&\quad 
+
\Bigg( \frac{288\kappa^2 L^2_f}{bKT\bar{\beta}}
\sum_{i=0}^{T-1}\mu^2_{y}(1-(1-\bar{\beta})^{T-i})\mE\|\bm^y_{c,i}\|^2-\frac{3\kappa L_f \mu_y}{2T}
\sum_{i=1}^{T}
\sum_{j=0}^{i-1}\Big(
1- \frac{\nu\mu_y}{2}\Big)^{i-j-1}
\mE\|\bm^y_{c,j}\|^2  \Bigg) 
\notag \\
&\quad 
+ \frac{620\kappa^2 L^2_f v^2_1v^2_2}{bK\bar{\beta}}
\Bigg(\frac{5(\mE\|\mhE_{x,0}\|^2 +\mE\|\mhE_{y,0}\|^2)}{T(1-\rho)} 
+\frac{288L^2_f\lambda^2_a\mu^2_y}{(1-\rho)^2\underline{\lambda^2_b}T}
\sum_{i=0}^{T-1}(\mu^2_x\mE\|\bm^x_{c,i}\|^2+\mu^2_y\mE\|\bm^y_{c,i}\|^2) \notag \\
&\quad 
+\frac{12\mu^2_y\lambda^2_a(p+\beta^2)}{(1-\rho)^2K\underline{\lambda^2_b}}
\Big(
\frac{2}{T\bar{\beta}}
(\mE\|\mS_{x,0}\|^2
+
\mE\|\mS_{y,0}\|^2)
+\frac{12KL^2_fv^2_1v^2_2}{bT\bar{\beta}}
(\mE\|\mhE_{x,0}\|^2 +\mE\|\mhE_{y,0}\|^2) 
\Big)
\notag \\
&\quad 
+ \frac{96\mu^2_y\lambda^2_a}{(1-\rho)^2\underline{\lambda^2_b}}
\Big(
\frac{p\sigma^2}{B} \mathbb{I}_{\text{online}}
+ \frac{\beta^2\sigma^2}{b}
\Big) \Bigg)+
\Bigg(  \frac{18\kappa^2}{T\mu_y}
\mE\Delta_{c,0} +\frac{48\kappa^2}{T\bar{\beta}}
(\mE\|\bs^x_{c,0}\|^2+
\mE\|\bs^y_{c,0}\|^2) 
\notag\\
&\quad 
+\frac{288\kappa^2 L^2_fv^2_1v^2_2}{bKT\bar{\beta}} (\mE\|\mhE_{x,0}\|^2+\mE\|\mhE_{y,0}\|^2) + \frac{96\kappa^2\beta^2\sigma^2}{Kb\bar{\beta}} +\frac{48\kappa^2p\sigma^2}{KB\bar{\beta}}\mathbb{I}_{\text{online}}
\Bigg).
\end{align}
Rearranging the terms, we get 
\begin{align}
&\frac{1}{T}
\sum_{i=0}^{T-1}
(\mE\|\nabla P(\bx_{c,i})\|^2
+1.5\kappa L_f\mE\Delta_{c,i}) \notag \\
&\le 
\frac{2\mE(P(\bx_{c,0})-P^\star)}{T\mu_x}-\Bigg(1-L\mu_x - \frac{3\kappa^2\mu_x}{2\mu_y} -\frac{288\kappa^2L^2_f\mu^2_x}{bK\bar{\beta}} - \frac{620\times 288 \kappa^2 L^4_fv^2_1v^2_2\lambda^2_a\mu^2_x\mu^2_y}{bK\bar{\beta}(1-\rho)^2\underline{\lambda^2_b}}\Bigg) \notag \\
&\quad \times \frac{1}{T}
\sum_{i=0}^{T-1}
\mE\|\bm^x_{c,i}\|^2+
\Bigg( \frac{288\kappa^2 L^2_f}{bKT\bar{\beta}}
\sum_{i=0}^{T-1}\mu^2_{y}(1-(1-\bar{\beta})^{T-i})\mE\|\bm^y_{c,i}\|^2
+\frac{620\times 288 \kappa^2 L^4_fv^2_1v^2_2\lambda^2_a\mu^4_y}{bK\bar{\beta}T(1-\rho)^2\underline{\lambda^2_b}}
\notag \\
&\quad \times \sum_{i=0}^{T-1}\mE\|\bm^y_{c,i}\|^2
-\frac{3\kappa L_f \mu_y}{2T}
\sum_{i=1}^{T}
\sum_{j=0}^{i-1}\Big(
1- \frac{\nu\mu_y}{2}\Big)^{i-j-1}
\mE\|\bm^y_{c,j}\|^2  \Bigg) + \frac{620\kappa^2 L^2_f v^2_1v^2_2}{bK\bar{\beta}}
\notag \\
&\quad \times \Bigg(\frac{5(\mE\|\mhE_{x,0}\|^2 +\mE\|\mhE_{y,0}\|^2)}{T(1-\rho)} 
+\frac{12\mu^2_y\lambda^2_a(p+\beta^2)}{(1-\rho)^2K\underline{\lambda^2_b}}\Big(
\frac{2}{T\bar{\beta}}
(\mE\|\mS_{x,0}\|^2
+
\mE\|\mS_{y,0}\|^2)
\notag \\
&\quad +\frac{12KL^2_fv^2_1v^2_2}{bT\bar{\beta}}
(\mE\|\mhE_{x,0}\|^2 +\mE\|\mhE_{y,0}\|^2) \Big)
+ \frac{96\mu^2_y\lambda^2_a}{(1-\rho)^2\underline{\lambda^2_b}}
\Big(
\frac{p\sigma^2}{B} \mathbb{I}_{\text{online}}
+ \frac{\beta^2\sigma^2}{b}
\Big) \Bigg)+
\Bigg( \frac{18\kappa^2}{T\mu_y}
\mE\Delta_{c,0} 
\notag\\
&\quad +
\frac{48\kappa^2}{T\bar{\beta}}
(\mE\|\bs^x_{c,0}\|^2+
\mE\|\bs^y_{c,0}\|^2) +\frac{288\kappa^2 L^2_fv^2_1v^2_2}{bKT\bar{\beta}} (\mE\|\mhE_{x,0}\|^2+\mE\|\mhE_{y,0}\|^2) 
\notag\\
&\quad + \frac{96\kappa^2\beta^2\sigma^2}{Kb\bar{\beta}}+\frac{48\kappa^2p\sigma^2}{KB\bar{\beta}}\mathbb{I}_{\text{online}}
\Bigg).
\end{align}
To eliminate the terms associated with $\frac{1}{T}\sum_{i=0}^{T-1} \mE\|\bm^x_{c,i}\|^2$, we choose $\mu_x \le \mu_y$ and  
\begin{align}
&L\mu_x \le \frac{1}{6} \Longrightarrow \mu_x \le \frac{1}{6L}, \quad  \frac{3\kappa^2\mu_x}{2\mu_y} \le \frac{1}{6} \Longrightarrow \frac{\mu_x}{\mu_y}
\le \frac{1}{9 \kappa^2},   \\
&
\frac{288\kappa^2 L^2_f \mu^2_x }{bK\bar{\beta}} \le \frac{1}{6} 
\Longrightarrow
\mu_x  \le \frac{\sqrt{Kb\bar{\beta}}}{24\sqrt{3}\kappa L_f}
, \notag \\
&
\frac{620\times 288 \kappa^2 L^4_fv^2_1v^2_2\lambda^2_a \mu^2_x\mu^2_y}{bK\bar{\beta}(1-\rho)^2\underline{\lambda^2_b}}
\le\frac{288\kappa^2L^2_f\mu^2_x}{bK\bar{\beta}}\Longrightarrow
\mu_y \le \frac{(1-\rho)\underline{\lambda_b}}{\sqrt{620}L_f\lambda_av_1v_2}.
\end{align}
Under these conditions, we have
\begin{align}
&\frac{1}{T}
\sum_{i=0}^{T-1}
(\mE\|\nabla P(\bx_{c,i})\|^2
+1.5\kappa L_f\mE\Delta_{c,i}) \notag \\
&\le 
\frac{2\mE(P(\bx_{c,0})-P^\star)}{T\mu_x} +
\Bigg( \frac{288\kappa^2 L^2_f}{bKT\bar{\beta}}
\sum_{i=0}^{T-1}\mu^2_{y}(1-(1-\bar{\beta})^{T-i})\mE\|\bm^y_{c,i}\|^2
+\frac{620\times 288 \kappa^2 L^4_fv^2_1v^2_2\lambda^2_a\mu^4_y}{bK\bar{\beta}T(1-\rho)^2\underline{\lambda^2_b}}
\notag \\
&\quad \times \sum_{i=0}^{T-1}\mE\|\bm^y_{c,i}\|^2
-\frac{3\kappa L_f \mu_y}{2T}
\sum_{i=1}^{T}
\sum_{j=0}^{i-1}\Big(
1- \frac{\nu\mu_y}{2}\Big)^{i-j-1}
\mE\|\bm^y_{c,j}\|^2  \Bigg) + \frac{620\kappa^2 L^2_f v^2_1v^2_2}{bK\bar{\beta}}
\notag \\
&\quad \times \Bigg(\frac{5(\mE\|\mhE_{x,0}\|^2 +\mE\|\mhE_{y,0}\|^2)}{T(1-\rho)} 
+\frac{12\mu^2_y\lambda^2_a(p+\beta^2)}{(1-\rho)^2K\underline{\lambda^2_b}}\Big(
\frac{2}{T\bar{\beta}}
(\mE\|\mS_{x,0}\|^2
+
\mE\|\mS_{y,0}\|^2)+\frac{12KL^2_fv^2_1v^2_2}{bT\bar{\beta}}
\notag \\
&\quad \times (\mE\|\mhE_{x,0}\|^2 +\mE\|\mhE_{y,0}\|^2) \Big)
+ \frac{96\mu^2_y\lambda^2_a}{(1-\rho)^2\underline{\lambda^2_b}}
\Big(
\frac{p\sigma^2}{B} \mathbb{I}_{\text{online}}
+ \frac{\beta^2\sigma^2}{b}
\Big) \Bigg)+
\Bigg( \frac{18\kappa^2}{T\mu_y}
\mE\Delta_{c,0} 
\notag\\
&\quad +
\frac{48\kappa^2}{T\bar{\beta}}
(\mE\|\bs^x_{c,0}\|^2+
\mE\|\bs^y_{c,0}\|^2) +\frac{288\kappa^2 L^2_fv^2_1v^2_2}{bKT\bar{\beta}} (\mE\|\mhE_{x,0}\|^2+\mE\|\mhE_{y,0}\|^2) 
\notag \\
&\quad 
+ \frac{96\kappa^2\beta^2\sigma^2}{Kb\bar{\beta}}+\frac{48\kappa^2p\sigma^2}{KB\bar{\beta}}\mathbb{I}_{\text{online}}
\Bigg).
\end{align}
Note that 
\begin{align}
&\frac{\kappa L_f\mu_y}{T}
\sum_{i=1}^{T}
\sum_{j=0}^{i-1}
\Big(1- \frac{\nu\mu_y}{2}\Big)^{i-j-1}
\mE\|\bm^y_{c,j}\|^2 \notag \\
&= \frac{\kappa L_f\mu_y}{ T}
\sum_{i=0}^{T-1}
\mE\|\bm^y_{c,i}\|^2
\Big(\sum_{j=0}^{T-1-i} \Big( 1-\frac{\nu\nu_y}{2}\Big)^j\Big) \notag \\
&= \frac{\kappa L_f\mu_y}{ T}
\sum_{i=0}^{T-1}
\mE\|\bm^y_{c,i}\|^2
\frac{1- \Big( 1- \frac{\nu \mu_y}{2}
\Big)^{T-i}}{1- \Big( 1- \frac{\nu \mu_y}{2}
\Big)} \notag \\
&=
\frac{2\kappa^2}{T}
\sum_{i=0}^{T-1} \Big( 1- \Big(1- \frac{\nu \mu_y}{2}\Big)^{T-i}\Big)
\mE\|\bm^y_{c,i}\|^2  \notag \\
&\ge 
\frac{2\kappa^2}{T}
\sum_{i=0}^{T-1}
\Big( 1- \Big(1- \frac{\nu \mu_y}{2}\Big)\Big)
\mE\|\bm^y_{c,i}\|^2 = \frac{\kappa^2 \nu \mu_y}{T}
\sum_{i=0}^{T-1}\mE\|\bm^y_{c,i}\|^2 \quad (\mu_y \le \frac{1}{\nu}) \label{appendix:mainproof:my}.
\end{align}
Splitting $\frac{3\kappa L_f \mu_y}{2T}
\sum_{i=1}^{T}
\sum_{j=0}^{i-1}\Big(
1- \frac{\nu\mu_y}{2}\Big)^{i-j-1}
\mE\|\bm^y_{c,j}\|^2$ into two parts and using the last two results of \eqref{appendix:mainproof:my}, we can deduce that
\begin{align}
&\Bigg( \frac{288\kappa^2 L^2_f}{bKT\bar{\beta}}
\sum_{i=0}^{T-1}\mu^2_{y}(1-(1-\bar{\beta})^{T-i})\mE\|\bm^y_{c,i}\|^2
+\frac{620\times 288 \kappa^2 L^4_fv^2_1v^2_2\lambda^2_a\mu^4_y}{bK\bar{\beta}T(1-\rho)^2\underline{\lambda}^2_b}
\sum_{i=0}^{T-1}\mE\|\bm^y_{c,i}\|^2
\notag \\
&\quad 
-\frac{3\kappa L_f \mu_y}{2T}
\sum_{i=1}^{T}
\sum_{j=0}^{i-1}\Big(
1- \frac{\nu\mu_y}{2}\Big)^{i-j-1}
\mE\|\bm^y_{c,j}\|^2  \Bigg) \notag \\
&\le 
\Bigg( \frac{288\kappa^2 L^2_f \mu^2_y}{bKT\bar{\beta}}
\sum_{i=0}^{T-1}(1-(1-\bar{\beta})^{T-i})\mE\|\bm^y_{c,i}\|^2
-\frac{2\kappa^2}{T}
\sum_{i=0}^{T-1}
\Big(1-\Big(
1- \frac{\nu\mu_y}{2}\Big)^{T-i}\Big)
\mE\|\bm^y_{c,i}\|^2 
\notag \\
&\quad 
+\frac{620\times 288 \kappa^2 L^4_fv^2_1v^2_2\lambda^2_a\mu^4_y}{bK\bar{\beta}T(1-\rho)^2\underline{\lambda}^2_b}
\sum_{i=0}^{T-1}\mE\|\bm^y_{c,i}\|^2
- \frac{\kappa^2 \nu \mu_y}{2T}\sum_{i=0}^{T-1}
\mE\|\bm^y_{c,i}\|^2\Bigg) .
\end{align} 
Choosing 
\begin{align}
&\bar{\beta} \le 1, \mu_y \le \frac{2}{\nu}, \quad (1-(1-\bar{\beta})^{T-i})
 \le \Big( 1- \Big(
 1- \frac{\nu\mu_y}{2}\Big)^{T-i}\Big) \Longrightarrow  \bar{\beta} \le \frac{\nu \mu_y}{2}  ,\\
 & \frac{288\kappa^2 L^2_f \mu^2_y}{bK\bar{\beta}} \le 2 \kappa^2 \Longrightarrow \frac{144L^2_f \mu^2_y}{K} \le b\bar{\beta}  \Longrightarrow
 \mu_y \le \frac{\sqrt{Kb\bar{\beta}}}{12L_f},\\
 & \frac{620\times 288 \kappa^2 L^4_f v^2_1v^2_2\lambda^2_a \mu^4_y}{bK\bar{\beta}(1-\rho)^2\underline{\lambda^2_b}} \le \frac{\kappa^2\nu \mu_y}{2} \Longleftarrow 
 \mu_y \le  \frac{(1-\rho)^{\frac{2}{3}}\underline{\lambda^{\frac{2}{3}}_b} (bK\bar{\beta})^{\frac{1}{3}}}{90L_f \kappa^{\frac{1}{3}}v^{\frac{2}{3}}_1 v^{\frac{2}{3}}_2 \lambda^{\frac{2}{3}}_a}.
\end{align}
Putting all results together, we obtain 
\begin{align}
&\frac{1}{T}
\sum_{i=0}^{T-1}
(\mE\|\nabla P(\bx_{c,i})\|^2
+1.5\kappa L_f\mE\Delta_{c,i}) \notag \\
&\le 
\frac{2\mE(P(\bx_{c,0})-P^\star)}{T\mu_x}  + \frac{620\kappa^2 L^2_f v^2_1v^2_2}{bK\bar{\beta}}
\Bigg(\frac{5(\mE\|\mhE_{x,0}\|^2 +\mE\|\mhE_{y,0}\|^2)}{T(1-\rho)} 
\notag \\
&\quad 
+\frac{12\mu^2_y\lambda^2_a(p+\beta^2)}{(1-\rho)^2K\underline{\lambda^2_b}}\Big(
\frac{2}{T\bar{\beta}}
(\mE\|\mS_{x,0}\|^2
+
\mE\|\mS_{y,0}\|^2)+\frac{12KL^2_fv^2_1v^2_2}{bT\bar{\beta}}
(\mE\|\mhE_{x,0}\|^2 +\mE\|\mhE_{y,0}\|^2) 
\Big)\notag \\
&\quad 
+ \frac{96\mu^2_y\lambda^2_a}{(1-\rho)^2\underline{\lambda^2_b}}
\Big(
\frac{p\sigma^2}{B} \mathbb{I}_{\text{online}}
+ \frac{\beta^2\sigma^2}{b}
\Big) \Bigg)+
\Bigg( \frac{18\kappa^2}{T\mu_y}
\mE\Delta_{c,0} +
\frac{48\kappa^2}{T\bar{\beta}}
(\mE\|\bs^x_{c,0}\|^2+
\mE\|\bs^y_{c,0}\|^2) 
\notag\\
&\quad +\frac{288\kappa^2 L^2_fv^2_1v^2_2}{bKT\bar{\beta}} (\mE\|\mhE_{x,0}\|^2+\mE\|\mhE_{y,0}\|^2) 
+ \frac{96\kappa^2\beta^2\sigma^2}{Kb\bar{\beta}}+\frac{48\kappa^2p\sigma^2}{KB\bar{\beta}}\mathbb{I}_{\text{online}}
\Bigg).
\label{appendix:proof:finalboundbefore}
\end{align}
Note that
\begin{align}
&\frac{1}{T}
\sum_{i=0}^{T-1}
(\mE\|\nabla_x J(\bx_{c,i}, \by_{c,i})\|^2
+\mE\|\nabla_y J(\bx_{c,i}, \by_{c,i})\|^2) \notag \\
& =
\frac{1}{T}
\sum_{i=0}^{T-1}
\Big(\mE\|\nabla_x J(\bx_{c,i}, \by_{c,i}) - \nabla_x J(\bx_{c,i}, \by^o(\bx_{c,i})) +\nabla_x J(\bx_{c,i}, \by^o(\bx_{c,i}))\|^2
+\mE\|\nabla_y J(\bx_{c,i}, \by_{c,i}) \notag \\
&\quad - 
\nabla_y J(\bx_{c,i}, \by^o(\bx_{c,i}))\|^2\Big) \notag \\
&\overset{(a)}{\le} 
\frac{1}{T}
\sum_{i=0}^{T-1}
\Big( 
4\mE\|\nabla P(\bx_{c,i})\|^2
+ \frac{7}{3}L^2_f\mE\|\by_{c,i} - \by^o(\bx_{c,i})\|^2
\Big) \notag \\
&\overset{(b)}{\le} 
\frac{1}{T}
\sum_{i=0}^{T-1}
\Big(
4\mE\|\nabla P(\bx_{c,i})\|^2
+ \frac{14}{3} 
\kappa L_f\mE\Delta_{c,i}
\Big) \notag \\
&\overset{(c)}{\le}
\frac{8\mE(P(\bx_{c,0})-P^\star)}{T\mu_x}  + \frac{2480\kappa^2 L^2_f v^2_1v^2_2}{bK\bar{\beta}}
\Bigg(\frac{5(\mE\|\mhE_{x,0}\|^2 +\mE\|\mhE_{y,0}\|^2)}{T(1-\rho)} 
\notag \\
&\quad 
+\frac{12\mu^2_y\lambda^2_a(p+\beta^2)}{(1-\rho)^2K\underline{\lambda^2_b}}\Big(
\frac{2}{T\bar{\beta}}
(\mE\|\mS_{x,0}\|^2
+
\mE\|\mS_{y,0}\|^2)+\frac{12KL^2_fv^2_1v^2_2}{bT\bar{\beta}}
(\mE\|\mhE_{x,0}\|^2 +\mE\|\mhE_{y,0}\|^2) 
\Big)\notag \\
&\quad 
+ \frac{96\mu^2_y\lambda^2_a}{(1-\rho)^2\underline{\lambda^2_b}}
\Big(
\frac{p\sigma^2}{B} \mathbb{I}_{\text{online}}
+ \frac{\beta^2\sigma^2}{b}
\Big) \Bigg)+
\Bigg( \frac{72\kappa^2}{T\mu_y}
\mE\Delta_{c,0} +
\frac{192\kappa^2}{T\bar{\beta}}
(\mE\|\bs^x_{c,0}\|^2+
\mE\|\bs^y_{c,0}\|^2) 
\notag\\
&\quad +\frac{1152\kappa^2 L^2_fv^2_1v^2_2}{bKT\bar{\beta}} (\mE\|\mhE_{x,0}\|^2+\mE\|\mhE_{y,0}\|^2) 
+ \frac{384\kappa^2\beta^2\sigma^2}{Kb\bar{\beta}}+\frac{192\kappa^2p\sigma^2}{KB\bar{\beta}}\mathbb{I}_{\text{online}}
\Bigg),
\end{align}
where $(a)$ follows from the expected $L_f$-smooth assumption and Jensen's inequality
$\|a+b\|^2 \le \frac{1}{t}\|a\|^2 + \frac{1}{1-t}\|b\|^2$, where $t= \frac{1}{4}$; $(b)$ follows from the quadratic growth property of a $\nu$-PL function \cite{karimi2016linear}; $(c)$ follows from 
\eqref{appendix:proof:finalboundbefore}.
In what follows, we summarize the hyperparameter conditions and refine the performance bound.

$\bullet$ \textbf{The bound of 
$\mE\|\mS_{x,0}\|^2, \mE\|\mS_{y,0}\|^2, \mE\|\bs^x_{c,0}\|^2, \mE\|\bs^y_{c,0}\|^2$}

Using Assumption 3 of Part I, it is not difficult to prove that 
\begin{align}
\mE\|\mS_{x,0}\|^2 &\le \frac{K\sigma^2}{b_0}, \quad    
\mE\|\mS_{y,0}\|^2 \le \frac{K\sigma^2}{b_0}, \\
\mE\|\bs^x_{c,0}\|^2 &\le \frac{\sigma^2}{b_0K}, \quad 
\mE\|\bs^y_{c,0}\|^2 \le \frac{\sigma^2}{b_0K}.
\end{align}

$\bullet$ 
\textbf{The bound of}
$\mE\|\mhE_{x,0}\|^2+\mE\|\mhE_{y,0}\|^2$

Recalling the definition for 
$\mhE_{x,0}, 
\mhE_{y,0}$, we get 
\begin{align}
&\mE\|\mhE_{x,0}\|^2
+ \mE\|\mhE_{y,0}\|^2
\notag \\
&=\mE
\Big\|
\frac{1}{\tau_x}
\widehat{\mQ}^{-1}_x
\begin{bmatrix}
\widehat{\mU}^\top_x \mX_0 \\
\widehat{\Lambda}^{-1}_{b_x}
\widehat{\mU}^\top_x
\mZ_{x,0}
\end{bmatrix}
\Big\|^2
+\mE
\Big\|
\frac{1}{\tau_y}
\widehat{\mQ}^{-1}_y
\begin{bmatrix}
\widehat{\mU}^\top_y \mY_0 \\
\widehat{\Lambda}^{-1}_{b_y}
\widehat{\mU}^\top_y
\mZ_{y,0}
\end{bmatrix}
\Big\|^2
\notag \\
&\overset{(a)}{\le}
\frac{1}{K\underline{\lambda^2_{b_x}}}\mE\|\widehat{\mU}^\top_x\mZ_{x,0}\|^2
+\frac{1}{K\underline{\lambda^2_{b_y}}}\mE\|\widehat{\mU
}^\top_y\mZ_{y,0}\|^2 \notag \\
&\le 
\frac{
\mE\|\widehat{\mU
}^\top_x\mZ_{x,0}\|^2
+\mE\|\widehat{\mU
}^\top_y\mZ_{y,0}\|^2
}{K\underline{\lambda^2_b}},
\end{align}
where $(a)$ 
 holds by setting identical local models, i.e.,
$\mX_{0} = \mathds{1}_{K} \otimes 
x_{0}, \mY_{0} = \mathds{1}_{K} \otimes 
y_{0}$ (for some $x_0 \in \mathbb{R}^{d_1}, y_0 \in \mathbb{R}^{d_2}$), and then the initial consensus error vanishes, i.e., $ \|\widehat{\mU}^\top_{x}\mX_{0}\|^2 = 0, \|\widehat{\mU}^\top_{y}\mY_{0}\|^2 = 0$.
Furthermore, we have 
\begin{align}
&\mE\|\widehat{\mU
}^\top_x\mZ_{x,0}\|^2
+\mE\|\widehat{\mU
}^\top_y\mZ_{y,0}\|^2
\notag \\
&= 
\mE\mZ^\top_{x,0}
\widehat{\mU
}_x
\widehat{\mU
}^\top_x
\widehat{\mU
}_x
\widehat{\mU
}^\top_x\mZ_{x,0}
+\mE\mZ^\top_{y, 0}
\widehat{\mU
}_y
\widehat{\mU
}^\top_y
\widehat{\mU
}_y
\widehat{\mU
}^\top_y\mZ_{y,0}
\notag \\
&= 
\mE\|\widehat{\mU
}_x
\widehat{\mU
}^\top_x\mZ_{x, 0}\|^2
+
\mE\|\widehat{\mU
}_y
\widehat{\mU
}^\top_y\mZ_{y, 0}\|^2
\notag \\
& =
\mE\Big\|\Big(\mathrm{I}_{K}
-\frac{1}{K} \mathds{1}_K \mathds{1}^\top_K  \Big)\otimes \mathrm{I}_{d_1}\mZ_{x,0}\Big\|^2
+
\mE\Big\|\Big(\mathrm{I}_{K}
-\frac{1}{K} \mathds{1}_{K} \mathds{1}^\top_{K}  \Big)\otimes \mathrm{I}_{d_2}\mZ_{y, 0}\Big\|^2
\notag \\
&=
\mE\Big\|
\Big(\mathrm{I}_{K}
-\frac{1}{K} \mathds{1}_{K} \mathds{1}^\top_{K}  \Big)\otimes \mathrm{I}_{d_1}\Big(\mu_x \mA_{x} 
\mM_{x,0}
+\mB_{x}
\mD_{x,0}
-\mB^2_{x}\mX_{0}\Big)
\Big\|^2
\notag \\
&\quad +
\mE\Big\|\Big(\mathrm{I}_{K}
-\frac{1}{K} \mathds{1}_{K} \mathds{1}^\top_{K}  \Big)\otimes \mathrm{I}_{d_2}\Big(-\mu_y \mA_{y} 
\mM_{y,0}
+\mB_{y}
\mD_{y,0}
-\mB^2_{y}\mY_{0}
\Big)\Big\|^2 \notag \\
&\overset{(a)}{=}
\mu^2_x 
\mE\Big\|
\Big(\mA_{x}  - 
\frac{1}{K}
\mathds{1}_{K}
\mathds{1}^\top_{K}
\otimes \mathrm{I}_{d_1}\Big)\mM_{x,0}
\Big\|^2
+
\mu^2_y 
\mE\Big\|
\Big(\mA_{y}  - 
\frac{1}{K}
\mathds{1}_{K}
\mathds{1}^\top_{K}
\otimes \mathrm{I}_{d_2}\Big)\mM_{y,0}
\Big\|^2
\notag \\
& \le 
K\mu^2_y \zeta^2_{0}  \quad (\mu_x \le \mu_y),
\end{align}
where 
$(a)$ follows from
Assumption 4 of Part I and the initialization
$
\mX_{0} = \mathds{1}_{K} \otimes 
x_{0}, \mY_{0} = \mathds{1}_{K} \otimes 
y_{0}, 
\mD
_{x,0} = 
\mathds{1}_{K} \otimes 
d_{x, 0}, \mD
_{y,0} = 
\mathds{1}_{K} \otimes 
d_{y, 0}
$
(for some $x_0, d_{x,0} \in \mathbb{R}^{d_1},y_0, d_{y,0} \in \mathbb{R}^{d_2}$).
Here,
$\zeta^2_0$ is defined as the initial discrepancy of the gradient estimator, namely
\begin{align}
\zeta^2_{0} 
\triangleq \frac{1}{K}
\mE\Big\|
\Big(\mA_{x}  - 
\frac{1}{K}
\mathds{1}_{K}
\mathds{1}^\top_{K}
\otimes \mathrm{I}_{d_1}\Big)\mM_{x,0}
\Big\|^2 +\frac{1}{K}
\mE\Big\|
\Big(\mA_{y}  - 
\frac{1}{K}
\mathds{1}_{K}
\mathds{1}^\top_{K}
\otimes \mathrm{I}_{d_2}\Big)\mM_{y,0}
\Big\|^2.
\end{align}
Hence, 
we have 
\begin{align}
&\mE
\|\mhE_{x,0}\|^2 +\mE
\|\mhE_{y,0}\|^2\le 
\frac{\mu^2_y \zeta^2_0}{\underline{\lambda^2_b}}.
\end{align}

Using the above results and denoting the initial gap $G_{p,0} \triangleq P(\bx_{c,0}) - P^\star$, 
we can get the final performance bound as follows 
\begin{align}
&\frac{1}{T}
\sum_{i=0}^{T-1}
(\mE\|\nabla_x J(\bx_{c,i}, \by_{c,i})\|^2
+\mE\|\nabla_y J(\bx_{c,i}, \by_{c,i})\|^2) \notag \\
&\le 
\frac{8\mE G_{p,0}}{T\mu_x}
+
\frac{12400\kappa^2L^2_fv^2_1v^2_2\zeta^2_0\mu^2_y}{bK\bar{\beta}T(1-\rho)\underline{\lambda^2_b}}
+\frac{2480\times48\kappa^2L^2_fv^2_1v^2_2\lambda^2_a\mu^2_y(p+\beta^2)\sigma^2}{bb_0K\bar{\beta}^2(1-\rho)^2\underline{\lambda^2_b}T} \notag \\
&
\quad +
\frac{2480\times 144 \kappa^2L^4_fv^4_1v^4_2\lambda^2_a\zeta^2_0\mu^4_y(p+\beta^2)}{b^2K\bar{\beta}^2(1-\rho)^2\underline{\lambda^4_b}T}
+ \frac{2480\times 96\kappa^2L^2_fv^2_1v^2_2\lambda^2_a\mu^2_y}{bK\bar{\beta}(1-\rho)^2\underline{\lambda^2_b}}
\Big(
\frac{p\sigma^2}{B}\mathbb{I}_{\text{online}} +\frac{\beta^2\sigma^2}{b}
\Big) \notag \\
&\quad 
+ \frac{72\kappa^2\mE\Delta_{c,0}}{T\mu_y}
+\frac{384\kappa^2\sigma^2}{b_0\bar{\beta}K T}
+ \frac{1152\kappa^2L^2_fv^2_1v^2_2\zeta^2_0\mu^2_y}{bK\bar{\beta}T \underline{\lambda^2_b}}
+\frac{384\kappa^2\beta^2\sigma^2}{Kb\bar{\beta}}
+\frac{192\kappa^2p\sigma^2}{KB\bar{\beta}}\mathbb{I}_{\text{online}}.
\end{align}

$\bullet$ \textbf{Summary of hyperparameter conditions}

Putting together all the hyperparameter conditions derived
from 
Lemma \ref{appendix:lemma:Lconsense}---\ref{appendix:lemma:costfunction}
and 
the proof of the main results, 
for the step sizes $\mu_x, \mu_y$, we get 
\begin{subequations}
\begin{align}
\mu_x &\le \min \Big\{
\frac{1}{32L}, \frac{\mu_y}{16\kappa^2}, \frac{\sqrt{Kb\bar{\beta}}}{24\sqrt{3}\kappa L_f} \Big\} ,  \label{appendix:stepsize:summarya}\\
\mu_y &\le 
\min \Big\{ 
\frac{1}{\nu}, \frac{1}{2L_f}, \frac{\sqrt{Kb\bar{\beta}}}{12L_f}, \frac{(1-\rho)\underline{\lambda_b}}{\sqrt{620}L_f v_1v_2 \lambda_a}, \frac{(1-\rho)\underline{\lambda_b}}{12L_f v_1v_2\lambda_a},\frac{(1-\rho)\underline{\lambda_b}}{24L_f v_1v_2\lambda_a}\sqrt{\frac{b \bar{\beta}}{{p+\beta^2}}},\frac{(1-\rho)^{\frac{2}{3}}\underline{\lambda^{\frac{2}{3}}_b} (bK\bar{\beta})^{\frac{1}{3}}}{90L_f \kappa^{\frac{1}{3}}v^{\frac{2}{3}}_1 v^{\frac{2}{3}}_2 \lambda^{\frac{2}{3}}_a}
\Big\},
\label{appendix:stepsize:summaryb}
\end{align}
for the smoothing factor $\beta_x =\beta_y =\beta$
and $\bar{\beta} = p+\beta - p\beta$, and minibatch size $b$, we have 
\begin{align}
\bar{\beta} \le \frac{\nu\mu_y}{2},\  b\bar{\beta} \le \frac{1}{K} , \ p+\beta \le 1, \ \beta+bp \le b , \ b\ge 1, \ \bar{\beta} \le 1, \ \beta_x =\beta_y =\beta \le 1.
\label{appendix:stepsize:summaryc}
\end{align}
While the hyperparameter conditions may seem convoluted at first, they can be satisfied in a straightforward manner by choosing parameters in a prescribed order with respect to the communication round 
$T$, provided 
$T$ is sufficiently large.
\end{subequations}
\end{proof}

\section{Refined performance bound for specific strategies}
\label{appendix:sec:special_results}

We establish the performance bound for a specific strategy. As a preliminary step, we summarize results from \cite{alghunaim2022unified} on the transition matrix to identify value of
$\lambda^2_a$ and $\underline{\lambda^2_b}$, included here for completeness.
\begin{Lemma}[Restatement of {\cite[Lemma 1]{alghunaim2022unified}}]
\label{appendix:lemma:combinationmatrix}
Given the matrices $\mA_{x}, \mB_{x}, \mC_{x}, \mA_{y}, \mB_{y}, \mC_{y}$ and $\mW_x = W \otimes \mathrm{I}_{d_1}, \mW_y = W \otimes \mathrm{I}_{d_2}$
that satisfy Assumption 
4 of Part I and assuming $W$ is positive semi-definite,
then 
\begin{align}
\mP_x &\triangleq \begin{bmatrix}
\widehat{\Lambda}_{a_{x}}\widehat{\Lambda}_{c_{x}} -\widehat{\Lambda}^2_{b_x} & - \widehat{\Lambda}_{b_x} \\
\widehat{\Lambda}_{b_x}  & \mathrm{I}_{(K-1)d_1} 
\end{bmatrix}\in \mathbb{R}^{2(K-1)d_1\times 2(K-1)d_1}, \notag \\
\mP_y &\triangleq \begin{bmatrix}
\widehat{\Lambda}_{a_{y}}\widehat{\Lambda}_{c_{y}} -\widehat{\Lambda}^2_{b_y}& - \widehat{\Lambda}_{b_y} \\
\widehat{\Lambda}_{b_y}  & \mathrm{I}_{(K-1)d_2}
\end{bmatrix} \in \mathbb{R}^{2(K-1)d_2\times 2(K-1)d_2},
\end{align}
admits the similarity transformation 
$\mP_x = \widehat{\mQ}_x \mT_x \widehat{\mQ}^{-1}_x, 
\mP_y = \widehat{\mQ}_y \mT_y \widehat{\mQ}^{-1}_y$
where $\|\mT_x\| = \|\mT_y\| < 1$. Furthermore, 
the following results hold:

$\bullet$ \textbf{ED.} Setting $\mA_x = \mW_{x}, \mC_x = \mathrm{I}_{Kd_1}, \mB_x =(\mathrm{I}_{Kd_1} - \mW_x)^{1/2}$,
$\mA_y = \mW_{y}, \mC_y = \mathrm{I}_{Kd_2}, \mB_y =(\mathrm{I}_{Kd_2} - \mW_y)^{1/2}$, assume $W \ge 0$, then  

\begin{align}
v^2_1&=\|\widehat{\mQ}_x\|^2=\|\widehat{\mQ}_y\|^2 \le 4, \quad v^2_2=\|\widehat{\mQ}^{-1}_x\|^2= \|\widehat{\mQ}^{-1}_y\|^2 \le \frac{2}{\underline{\lambda}}, \\
\rho &= \rho_x = \rho_y = \sqrt{\lambda}, \quad \lambda_a = \lambda_{a_x} = \lambda_{a_y} = \lambda, \quad \underline{\lambda_b}=\underline{\lambda_{b_x}} = \underline{\lambda_{b_y}} = \sqrt{1-\lambda}.
\end{align}
where $\lambda = \max_{i=2,\dots, K} \lambda_i$ and $\underline{\lambda}$
is the minimum non-zero eigenvalue of $W$.

$\bullet$ \textbf{EXTRA.}
$\mA_{x} = \mathrm{I}_{Kd_1}, \mC_x = \mW_x, \mB_x = (\mathrm{I}_{Kd_1} - \mW_x)^{1/2}$,
$\mA_{y} = \mathrm{I}_{Kd_2}, \mC_y = \mW_y, \mB_y = (\mathrm{I}_{Kd_2} - \mW_y)^{1/2}$, then 
\begin{align}
v^2_1&=\|\widehat{\mQ}_x\|^2=\|\widehat{\mQ}_y\|^2 \le 4, \quad v^2_2=\|\widehat{\mQ}^{-1}_x\|^2= \|\widehat{\mQ}^{-1}_y\|^2 \le \frac{2}{\underline{\lambda}}, \\
\rho &= \rho_x= \rho_y = \sqrt{\lambda}, \quad \lambda_a = \lambda_{a_x} = \lambda_{a_y} = 1, \quad \underline{\lambda_b}=\underline{\lambda_{b_x}} = \underline{\lambda_{b_y}} = \sqrt{1-\lambda}
.\end{align}

$\bullet$ \textbf{ATC-GT.}
 $\mA_x = \mW^2_x,\mC_x = \mathrm{I}_{Kd_1}, \mB_x = (\mathrm{I}_{Kd_1} - \mW_x)$
and  $\mA_y = \mW^2_y,\mC_y = \mathrm{I}_{Kd_2}, \mB_y = (\mathrm{I}_{Kd_2} - \mW_y)$, then 
\begin{align}
v^2_1&=\|\widehat{\mQ}_x\|^2=\|\widehat{\mQ}_y\|^2 \le 3, \quad v^2_2=\|\widehat{\mQ}^{-1}_x\|^2= \|\widehat{\mQ}^{-1}_y\|^2 \le 9, \\
\rho &= \rho_x= \rho_y \le \frac{1+\lambda}{2}, \quad \lambda_a = \lambda_{a_x} = \lambda_{a_y} = \lambda^2, \quad \underline{\lambda_b}=\underline{\lambda_{b_x}} = \underline{\lambda_{b_y}} = 1-\lambda
.\end{align}

\end{Lemma}

\subsection{Performance bound of STORM-based algorithms}
\label{appendix:corollary:storm}
By setting $\gamma_1=1, \gamma_2=0, p=0$, and $\beta_x=\beta_y\neq 0$, we deactivate large-batch computations, and \textbf{DAMA} reduces to the STORM-based algorithm. The following results hold for online scenarios.

\begin{Corollary}[\textbf{STORM+ED}]
\label{appendix:corollary:sd+storm}
Under Theorem \ref{main:theorem:main}
and assuming $W$ is positive semi-definite, if we consider the matrix choice of ED and set the hyperparameters for \textbf{DAMA} as follows
\begin{align}
\mu_x = \mathcal{O}\Big(\frac{K^{2/3}}{T^{1/3}\kappa^2}\Big), \mu_y = \mathcal{O} \Big( \frac{K^{2/3}}{T^{1/3}}\Big), \beta_x = \beta_y = \mathcal{O}
\Big(\frac{K^{1/3}}{T^{2/3}}\Big), p =0,  b_0 = \mathcal{O}\Big(\frac{T^{1/3}}{K^{2/3}}\Big), b=\mathcal{O}(1).
\label{appendix:condition:STORM+ED}
\end{align}
It can be verified that the hyperparameter conditions 
in Theorem \ref{main:theorem:main}
are satisfied for sufficiently large $T$, i.e.,
\begin{subequations}
  \begin{align}
\mu_x &= \mathcal{O}\Big(\frac{K^{2/3}}{T^{1/3}\kappa^2}\Big) \le \min \Big\{
\frac{1}{32L}, \mathcal{O}
\Big(
\frac{K^{2/3}}{\kappa^2T^{1/3}} \Big),
 \mathcal{O}
\Big(
\frac{K^{2/3}}{\kappa T^{1/3}}
\Big) 
\Big\}, \\
\mu_y &= 
\mathcal{O}
\Big(
\frac{K^{2/3}}{T^{1/3}}\Big) \le
\min \Bigg\{ 
\frac{1}{\nu},
\frac{1}{2L_f},
\mathcal{O}\Big(
\frac{K^{2/3}}{T^{1/3}}
\Big), \mathcal{O}\Big(
\frac{(1-\sqrt{\lambda})\sqrt{1-\lambda}}{L_fv_1v_2\lambda} \Big), \mathcal{O}
\Big(
\frac{(1-\sqrt{\lambda})\sqrt{1-\lambda}T^{1/3}}{L_fv_1v_2\lambda K^{1/6}}\Big), \notag\\
&\qquad \mathcal{O}\Big(
\frac{(1-\sqrt{\lambda})^{\frac{2}{3}}(1-\lambda)^{\frac{1}{3}}K^{4/9}}{L_f\kappa^{\frac{1}{3}}v^{\frac{2}{3}}_1v^{\frac{2}{3}}_2\lambda^{\frac{2}{3}} T^{2/9}}
\Big)
\Bigg\}, \\
\beta &= \beta_x =\beta_y= \mathcal{O}\Big(\frac{K^{1/3}}{T^{2/3}}\Big)
\le \min\Big\{
\mathcal{O}\Big(\frac{K^{2/3}}{T^{1/3}}\Big), 1, \frac{1}{K}
\Big\}. \label{appendix:stability:STORM+ED}
\end{align}  
\end{subequations}
Then the performance bound is given by
\begin{align}
&\frac{1}{T}
\sum_{i=0}^{T-1}
(\mE\|\nabla_x J(\bx_{c,i}, \by_{y,i})\|^2 + \mE\|\nabla_y J(\bx_{c,i}, \by_{c,i})\|^2) \notag \\
&\le \mathcal{O}
\Big( 
\frac{
\mE G_{p,0} \kappa^2
}{(TK)^{2/3}} +  \frac{\kappa^2 \zeta^2_0}{T(1-\lambda)^2}
+ \frac{\kappa^2\lambda^2K\sigma^2}{T^2(1-\lambda)^3}
+ \frac{\kappa^2\lambda^2K^{5/3}\zeta^2_0}{T^{7/3}(1-\lambda)^4}+ \frac{\kappa^2\lambda^2K^{2/3}\sigma^2}{T^{4/3}(1-\lambda)^3}
 \notag \\
&\quad 
+ \frac{\kappa^2\mE\Delta_{c,0}}{(KT)^{2/3}}
+ \frac{\kappa^2\sigma^2}{(KT)^{2/3}} + \frac{\kappa^2\zeta^2_0}{T(1-\lambda)} + \frac{\kappa^2\sigma^2}{(KT)^{2/3}}
\Big)
\notag \\
&\le 
 \mathcal{O}
\Big( 
\frac{\kappa^2
\mE G_{p,0} 
}{(TK)^{2/3}} +  \frac{\kappa^2\mE\Delta_{c,0}}{(KT)^{2/3}} + \frac{\kappa^2\sigma^2}{(KT)^{2/3}} + \frac{\kappa^2 \zeta^2_0}{T(1-\lambda)^2}
+ \frac{\kappa^2\lambda^2K^{2/3}\sigma^2}{T^{4/3}(1-\lambda)^3}+ \frac{\kappa^2\lambda^2K\sigma^2}{T^2(1-\lambda)^3} \notag \\
&\quad 
+ \frac{\kappa^2\lambda^2K^{5/3}\zeta^2_0}{T^{7/3}(1-\lambda)^4} 
\Big).
\end{align}
Note that the above convergence rate is dominated by 
$\mathcal{O}\Big(\frac{\kappa^2}{(TK)^{2/3}} + \frac{\kappa^2}{T(1-\lambda)^2}\Big)$. In addition, the transient time in achieving linear speedup is given by the choice of $T$ such that 
\begin{align}
&\mathcal{O}\Big(\frac{\kappa^2}{(KT)^{2/3}}\Big) \ge \mathcal{O}\Big(\frac{\kappa^2}{T(1-\lambda)^2} \Big)\Longrightarrow T \ge \mathcal{O}\Big( \frac{K^2}{(1-\lambda)^6}\Big), \notag \\
&\mathcal{O}\Big(
\frac{\kappa^2}{(KT)^{2/3}} \Big)\ge \mathcal{O}\Big(
\frac{\lambda^2\kappa^2K^{2/3}}{T^{4/3}(1-\lambda)^3}
\Big)\Longrightarrow T \ge \mathcal{O}\Big(\frac{\lambda^3K^2}{(1-\lambda)^{4.5}}\Big), \notag \\
&\mathcal{O}\Big(
\frac{\kappa^2}{(KT)^{2/3}} \Big) \ge \mathcal{O}\Big(
\frac{\lambda^2\kappa^2K}{T^{2}(1-\lambda)^3}
\Big)\Longrightarrow T \ge \mathcal{O}\Big(\frac{\lambda^{3/2} K^{5/4}}{(1-\lambda)^{9/4}}\Big), \notag \\
&\mathcal{O}\Big(
\frac{\kappa^2}{(KT)^{2/3}} \ge \mathcal{O}\Big(
\frac{\lambda^2\kappa^2K^{5/3}}{T^{7/3}(1-\lambda)^4}
\Big)\Longrightarrow T \ge \mathcal{O}\Big(\frac{\lambda^{6/5}K^{7/5}}{(1-\lambda)^{2.4}}\Big).
\end{align}
Therefore, the transient time $T$ is given by
\begin{align}
T = \max \Big\{\mathcal{O}\Big(\frac{\lambda^{6/5}K^{7/5}}{(1-\lambda)^{2.4}}\Big), \mathcal{O}\Big(\frac{\lambda^{3/2}K^{5/4}}{(1-\lambda)^{9/4}}\Big), \mathcal{O}\Big(\frac{\lambda^{3}K^2}{(1-\lambda)^{4.5}}\Big), \mathcal{O}\Big( \frac{K^2}{(1-\lambda)^6}\Big)\Big\}.
\end{align}
Furthermore, when $T$ is sufficiently large, the sample complexity is given by 
\begin{align}
SC=\mathcal{O}\Big(
    \frac{\kappa^3\varepsilon^{-3}}{K}+\frac{\kappa^2\varepsilon^{-2}}{(1-\lambda)^2} +\frac{\kappa^{1.5} \lambda^{1.5}K^{0.5}\varepsilon^{-1.5}}{(1-\lambda)^{9/4}}
\Big).
\end{align}
Note that the first term dominates the sample complexity since $\varepsilon$ is usually chosen small. 
\end{Corollary}
\begin{Corollary}[\textbf{STORM+EXTRA}]
\label{appendix:corollary:extra+storm}
 Under Theorem \ref{main:theorem:main} and assuming $W$ is positive semi-definite, if we consider the matrix choice of EXTRA and set the hyperparameters for \textbf{DAMA} similar to 
Corollary \ref{appendix:corollary:sd+storm},
we get 
\begin{align}
&\frac{1}{T}
\sum_{i=0}^{T-1}
(\mE\|\nabla_x J(\bx_{c,i}, \by_{y,i})\|^2 + \mE\|\nabla_y J(\bx_{c,i}, \by_{c,i})\|^2) \notag \\
&\le  \mathcal{O}
\Big( 
\frac{\kappa^2
\mE G_{p,0} 
}{(TK)^{2/3}} +  \frac{\kappa^2\mE\Delta_{c,0}}{(KT)^{2/3}} + \frac{\kappa^2\sigma^2}{(KT)^{2/3}} + \frac{\kappa^2 \zeta^2_0}{T(1-\lambda)^2}
+ \frac{\kappa^2K^{2/3}\sigma^2}{T^{4/3}(1-\lambda)^3}+ \frac{\kappa^2K\sigma^2}{T^2(1-\lambda)^3} \notag \\
&\quad 
+ \frac{\kappa^2K^{5/3}\zeta^2_0}{T^{7/3}(1-\lambda)^4} 
\Big).
\end{align}
The transient time of STORM+EXTRA is given by 
\begin{align}
T = \max \Big\{\mathcal{O}\Big(\frac{K^{7/5}}{(1-\lambda)^{2.4}}\Big), \mathcal{O}\Big(\frac{K^{5/4}}{(1-\lambda)^{9/4}}\Big), \mathcal{O}\Big(\frac{K^2}{(1-\lambda)^{4.5}}\Big), \mathcal{O}\Big( \frac{K^2}{(1-\lambda)^6}\Big)\Big\}.
\end{align}
The sample complexity of STORM+EXTRA is given by 
\begin{align}
SC=\mathcal{O}\Big(
    \frac{\kappa^3\varepsilon^{-3}}{K}+\frac{\kappa^2\varepsilon^{-2}}{(1-\lambda)^2} +\frac{\kappa^{1.5} K^{0.5}\varepsilon^{-1.5}}{(1-\lambda)^{9/4}}
\Big).
\end{align}
\end{Corollary}
\begin{Corollary}[\textbf{STORM+ATC-GT}]
\label{appendix:corollary:storm+atc-gt}
Under Theorem \ref{main:theorem:main} and assuming $W$ is positive semi-definite, if we consider the matrix choice of ATC-GT and set the hyperparameters for \textbf{DAMA} similar to 
Corollary \ref{appendix:corollary:sd+storm},
it is not difficult to verify that the hyperparameter conditions in \eqref{appendix:stepsize:summarya}---\eqref{appendix:stepsize:summaryc} are satisfied for sufficiently large $T$.
Then, we get 
\begin{align}
&\frac{1}{T}
\sum_{i=0}^{T-1}
(\mE\|\nabla_x J(\bx_{c,i}, \by_{y,i})\|^2 + \mE\|\nabla_y J(\bx_{c,i}, \by_{c,i})\|^2) \notag \\
&\le \mathcal{O}
\Big( 
\frac{
\mE G_{p,0} \kappa^2
}{(TK)^{2/3}} +  \frac{\kappa^2 \zeta^2_0}{T(1-\lambda)^3}
+ \frac{\kappa^2\lambda^4K\sigma^2}{T^2(1-\lambda)^4}
+ \frac{\kappa^2\lambda^4K^{5/3}\zeta^2_0}{T^{7/3}(1-\lambda)^6}+ \frac{\kappa^2\lambda^4K^{2/3}\sigma^2}{T^{4/3}(1-\lambda)^4}
 \notag \\
&\quad 
+ \frac{\kappa^2\mE\Delta_{c,0}}{(KT)^{2/3}}
+ \frac{\kappa^2\sigma^2}{(KT)^{2/3}} + \frac{\kappa^2\zeta^2_0}{T(1-\lambda)^2} + \frac{\kappa^2\sigma^2}{(KT)^{2/3}}
\Big)  \notag \\
&\le 
 \mathcal{O}
\Big( 
\frac{\kappa^2
\mE G_{p,0} 
}{(TK)^{2/3}} +  \frac{\kappa^2\mE\Delta_{c,0}}{(KT)^{2/3}} + \frac{\kappa^2\sigma^2}{(KT)^{2/3}} + \frac{\kappa^2 \zeta^2_0}{T(1-\lambda)^3}
+ \frac{\kappa^2\lambda^4K^{2/3}\sigma^2}{T^{4/3}(1-\lambda)^4}+ \frac{\kappa^2\lambda^4K\sigma^2}{T^2(1-\lambda)^4} \notag \\
&\quad + \frac{\kappa^2\lambda^4K^{5/3}\zeta^2_0}{T^{7/3}(1-\lambda)^6} 
\Big).
\end{align}
Note that the above convergence rate is dominated by 
$\mathcal{O}\Big(\frac{\kappa^2}{(TK)^{2/3}} + \frac{\kappa^2}{T(1-\lambda)^3}\Big)$. The transient time in achieving linear speedup is given by the choice of $T$ such that 
\begin{align}
&\mathcal{O}\Big(\frac{\kappa^2}{(KT)^{2/3}}\Big) \ge \mathcal{O}\Big(\frac{\kappa^2}{T(1-\lambda)^3} \Big)\Longrightarrow T \ge \mathcal{O}\Big( \frac{K^2}{(1-\lambda)^9}\Big), \notag \\
&\mathcal{O}\Big(
\frac{\kappa^2}{(KT)^{2/3}} \Big)\ge \mathcal{O}\Big(
\frac{\kappa^2\lambda^4K^{2/3}}{T^{4/3}(1-\lambda)^4}
\Big)\Longrightarrow T \ge \mathcal{O}\Big(\frac{\lambda^6K^2}{(1-\lambda)^{6}}\Big), \notag \\
&\mathcal{O}\Big(
\frac{\kappa^2}{(KT)^{2/3}} \Big) \ge \mathcal{O}\Big(
\frac{\kappa^2\lambda^4K}{T^{2}(1-\lambda)^4}
\Big)\Longrightarrow T \ge \mathcal{O}\Big(\frac{\lambda^3K^{5/4}}{(1-\lambda)^{3}}\Big), \notag \\
&\mathcal{O}\Big(
\frac{\kappa^2}{(KT)^{2/3}} \ge \mathcal{O}\Big(
\frac{\kappa^2\lambda^4K^{5/3}}{T^{7/3}(1-\lambda)^6}
\Big)\Longrightarrow T \ge \mathcal{O}\Big(\frac{\lambda^{12/5}K^{7/5}}{(1-\lambda)^{3.6}}\Big).
\end{align}
Therefore, the transient time  is given by
\begin{align}
T = \max \Big\{\mathcal{O}\Big(\frac{\lambda^{12/5}K^{7/5}}{(1-\lambda)^{3.6}}\Big), \mathcal{O}\Big(\frac{\lambda^3K^{5/4}}{(1-\lambda)^{3}}\Big), \mathcal{O}\Big(\frac{\lambda^6K^2}{(1-\lambda)^{6}}\Big), \mathcal{O}\Big( \frac{K^2}{(1-\lambda)^9}\Big)\Big\}.
\end{align}
Furthermore, when $T$ is sufficiently large, the sample complexity is given by 
\begin{align}
SC=\mathcal{O} \Big(\frac{\kappa^3\varepsilon^{-3}}{K}+\frac{\kappa^2\varepsilon^{-2}}{(1-\lambda)^3}+\frac{\kappa^{1.5}\lambda^3K^{0.5}\varepsilon^{-1.5}}{(1-\lambda)^{3}} 
\Big).
\end{align}
\end{Corollary}
\subsection{Performance bound of PAGE-based algorithms in online/offline scenarios}
By setting $p \not =0$, $\beta_x= \beta_y = 0$, and setting $b$ in a specific order,
we activate the strategy of PAGE.
Below, we discuss the performance bound of PAGE-based algorithms in both online/offline scenarios.

\subsubsection{Offline scenario}
\begin{Corollary}[\textbf{PAGE+ED}] 
\label{appendix:corollary:ed+page:offline}
Under Theorem \ref{main:theorem:main},  assuming $W$ is positive semi-definite and the local sample size $N_k \equiv N$ for simplicity, if we consider the matrix choice of ED and set $B=N$ and $b = b_0=\sqrt{\frac{N}{K}}, p = \frac{1}{\sqrt{KN}}, \beta_x=\beta_y=0$  under a large $N$ and a sparsely connected network, i.e., $\lambda \rightarrow 1$,  the hyperparameter conditions  \eqref{appendix:stepsize:summarya}---\eqref{appendix:stepsize:summaryc} become 
\begin{align}
\mu_x &\le \min \Big\{
\frac{1}{32L}, \frac{\mu_y}{16\kappa^2}, \frac{1}{24\sqrt{3}\kappa L_f} \Big\} , \\
\mu_y &\le 
\min \Big\{ 
\frac{1}{\nu}, \frac{1}{2L_f}, \frac{1}{12L_f}, \frac{(1-\sqrt{\lambda})\sqrt{1-\lambda}}{\sqrt{620}L_f v_1v_2 \lambda},  \frac{(1-\sqrt{\lambda})\sqrt{1-\lambda}}{12L_f v_1v_2\lambda},\frac{(1-\sqrt{\lambda})\sqrt{1-\lambda}}{24L_f v_1v_2\lambda}  \Big(\frac{N}{K}\Big)^{1/4}, \notag \\
&\quad 
\frac{(1-\sqrt{\lambda})^{\frac{2}{3}}(1-\lambda)^{1/3}}{90L_f \kappa^{\frac{1}{3}}v^{\frac{2}{3}}_1 v^{\frac{2}{3}}_2 \lambda^{\frac{2}{3}}}
\Big\}, \notag \\
p &\le \frac{\nu \mu_y}{2}, bp\le \frac{1}{K}, p \le 1.
\end{align}
We can then choose $\mu_x = \Theta \Big(\frac{(1-\lambda)^{1.5}}{\kappa^2}\Big)$ and $\mu_y = \Theta((1-\lambda)^{1.5})$, and the performance bound is given by
\begin{align}
&\frac{1}{T}
\sum_{i=0}^{T-1}
(\mE\|\nabla_x J(\bx_{c,i}, \by_{y,i})\|^2 + \mE\|\nabla_y J(\bx_{c,i}, \by_{c,i})\|^2) \notag \\
&\le \mathcal{O}
\Big( 
\frac{\kappa^2 \mE G_{p,0}}{T(1-\lambda)^{1.5}}
+\frac{\kappa^2\zeta^2_0(1-\lambda)}{T}
+ \frac{\kappa^2\lambda^2\sqrt{K}\sigma^2}{\sqrt{N}T}
+ \frac{\kappa^2\lambda^2(1-\lambda)^2\sqrt{K}\zeta^2_0}{\sqrt{N}T} + \frac{\kappa^2\mE\Delta_{c,0}}{T(1-\lambda)^{1.5}}
+\frac{\kappa^2\sigma^2}{T}
\notag \\
&\quad 
+\frac{\kappa^2(1-\lambda)^2\zeta^2_0}{T}
\Big) \notag \\
&\le  \mathcal{O}
\Big( 
\frac{\kappa^2 \mE G_{p,0}}{T(1-\lambda)^{1.5}}+ \frac{\kappa^2\mE\Delta_{c,0}}{T(1-\lambda)^{1.5}}
+\frac{\kappa^2\zeta^2_0(1-\lambda)}{T}+\frac{\kappa^2\sigma^2}{T}
+ \frac{\kappa^2\lambda^2\sqrt{K}\sigma^2}{\sqrt{N}T}
+ \frac{\kappa^2\lambda^2(1-\lambda)^2\sqrt{K}\zeta^2_0}{\sqrt{N}T} 
\Big).
\end{align}
Therefore, the communication complexity in achieving an $\varepsilon$-stationary point is given by 
\begin{align}
    T = \mathcal{O}\Big(
    \frac{\kappa^2\varepsilon^{-2}}{(1-\lambda)^{1.5}}
    + \kappa^2\varepsilon^{-2} +  \kappa^2\varepsilon^{-2}(1-\lambda)\Big) \approx \mathcal{O}\Big(
    \frac{\kappa^2\varepsilon^{-2}}{(1-\lambda)^{1.5}}
    \Big).
\end{align}
Hence,
the sample complexity is given by 
\begin{align}
SC &\approx \mathcal{O}\Big( \frac{\kappa^2\varepsilon^{-2}}{(1-\lambda)^{1.5}}
\times (1-p) \times b +  \frac{\kappa^2\varepsilon^{-2}}{(1-\lambda)^{1.5}} \times p \times N  + b_0\Big)\notag \\
&=\mathcal{O}\Big(
\frac{\kappa^2\sqrt{N}\varepsilon^{-2}}{\sqrt{K}(1-\lambda)^{1.5}}
+\frac{\kappa^2\sqrt{N}\varepsilon^{-2}}{\sqrt{K}(1-\lambda)^{1.5}} + \sqrt{\frac{N}{K}}
\Big).
\end{align}
\end{Corollary}
\begin{Corollary}[\textbf{PAGE+EXTRA}]
\label{appendix:corollary:extra+page:offline}
Under Theorem \ref{main:theorem:main}, assuming $W$ is positive semi-definite and  the local sample size $N_k \equiv N$ for simplicity, if we consider the matrix choice of EXTRA and set $B=N$ and $b = b_0=\sqrt{\frac{N}{K}}, p = \frac{1}{\sqrt{KN}}$, under sparsely connected network, i.e., $\lambda \rightarrow 1$,
the step size condition is similar to Corollary \ref{appendix:corollary:ed+page:offline}.
Choosing $\mu_x = \Theta \Big(\frac{(1-\lambda)^{1.5}}{\kappa^2}\Big)$ and $\mu_y = \Theta((1-\lambda)^{1.5})$,
the performance bound is given by 
\begin{align}
&\frac{1}{T}
\sum_{i=0}^{T-1}
(\mE\|\nabla_x J(\bx_{c,i}, \by_{y,i})\|^2 + \mE\|\nabla_y J(\bx_{c,i}, \by_{c,i})\|^2) \notag \\
&\le  \mathcal{O}
\Big( 
\frac{\kappa^2 \mE G_{p,0}}{T(1-\lambda)^{1.5}}+ \frac{\kappa^2\mE\Delta_{c,0}}{T(1-\lambda)^{1.5}}
+\frac{\kappa^2\zeta^2_0(1-\lambda)}{T}+\frac{\kappa^2\sigma^2}{T}
+ \frac{\kappa^2\sqrt{K}\sigma^2}{\sqrt{N}T}
+ \frac{\kappa^2(1-\lambda)^2\sqrt{K}\zeta^2_0}{\sqrt{N}T} 
\Big) .
\end{align}
Then, the communication complexity and sample complexity are given by 
\begin{align}
  T &= \mathcal{O}\Big(
   \frac{\kappa^2\varepsilon^{-2}}{(1-\lambda)^{1.5}}
    + \kappa^2\varepsilon^{-2} +  \kappa^2\varepsilon^{-2}(1-\lambda)\Big) \approx \mathcal{O}\Big(
   \frac{\kappa^2\varepsilon^{-2}}{(1-\lambda)^{1.5}}
   \Big), \\
   SC &\approx \mathcal{O} \Big(\frac{\kappa^2\sqrt{N}\varepsilon^{-2}}{\sqrt{K}(1-\lambda)^{1.5}} + \sqrt{\frac{N}{K}}\Big).
\end{align}
\end{Corollary}
\begin{Corollary}[\textbf{PAGE+ATC-GT}]
\label{appendix:corollary:gt-atc+page:offline}
Under Theorem \ref{main:theorem:main}, assuming $W$ is positive semi-definite and the local sample size $N_k \equiv N$ for simplicity, if we consider the matrix choice of ATC-GT and set $B=N$ and $b = b_0=\sqrt{\frac{N}{K}}, p = \frac{1}{\sqrt{KN}}$  for a large $N$, under sparsely connected network, i.e., $\lambda \rightarrow 1$,  the hyperparameter condition \eqref{appendix:stepsize:summarya}---\eqref{appendix:stepsize:summaryc} become 
\begin{align}
\mu_x &\le \min \Big\{
\frac{1}{32L}, \frac{\mu_y}{16\kappa^2}, \frac{1}{24\sqrt{3}\kappa L_f}\Big\} , \\
\mu_y &\le 
\min \Big\{ 
\frac{1}{\nu}, \frac{1}{2L_f}, \frac{1}{12L_f}, \frac{(1-\lambda)^2}{2\sqrt{620}L_f v_1v_2 \lambda^2},  \frac{(1-\lambda)^2}{24L_f v_1v_2\lambda^2},\frac{(1-\lambda)^2}{48L_f v_1v_2\lambda^2}\Big(\frac{N}{K}\Big)^{1/4}, \frac{(\frac{1-\lambda}{2})^{\frac{2}{3}}(1-\lambda)^{2/3}}{90L_f \kappa^{\frac{1}{3}}v^{\frac{2}{3}}_1 v^{\frac{2}{3}}_2 \lambda^{\frac{4}{3}}}
\Big\}, \notag \\
p &\le \frac{\nu \mu_y}{2}, bp\le \frac{1}{K}, p \le 1.
\end{align}
Choosing $\mu_x = \Theta \Big(\frac{(1-\lambda)^{2}}{\kappa^2}\Big)$ and $\mu_y = \Theta((1-\lambda)^{2})$,
the performance bound is given by 
\begin{align}
&\frac{1}{T}
\sum_{i=0}^{T-1}
(\mE\|\nabla_x J(\bx_{c,i}, \by_{y,i})\|^2 + \mE\|\nabla_y J(\bx_{c,i}, \by_{c,i})\|^2) \notag \\
&\le \mathcal{O}
\Big( 
\frac{\kappa^2 \mE G_{p,0}}{T(1-\lambda)^{2}}
+\frac{\kappa^2\zeta^2_0(1-\lambda)}{T}
+ \frac{\kappa^2\lambda^4\sqrt{K}\sigma^2}{\sqrt{N}T}
+ \frac{\kappa^2\lambda^4(1-\lambda)^2\sqrt{K}\zeta^2_0}{\sqrt{N}T} + \frac{\kappa^2\mE\Delta_{c,0}}{T(1-\lambda)^{2}}
+\frac{\kappa^2\sigma^2}{T}
\notag \\
&\quad 
+\frac{\kappa^2(1-\lambda)^2\zeta^2_0}{T}
\Big) \notag \\
&\le  \mathcal{O}
\Big( 
\frac{\kappa^2 \mE G_{p,0}}{T(1-\lambda)^{2}}+ \frac{\kappa^2\mE\Delta_{c,0}}{T(1-\lambda)^{2}}
+\frac{\kappa^2\zeta^2_0(1-\lambda)}{T}+\frac{\kappa^2\sigma^2}{T}
+ \frac{\kappa^2\lambda^4\sqrt{K}\sigma^2}{\sqrt{N}T}
+ \frac{\kappa^2\lambda^4(1-\lambda)^2\sqrt{K}\zeta^2_0}{\sqrt{N}T} 
\Big) .
\end{align}
The communication complexity and sample complexity are given by 
\begin{align}
  T = \mathcal{O}\Big(
    \frac{\kappa^2\varepsilon^{-2}}{(1-\lambda)^{2}}
    + \kappa^2\varepsilon^{-2} +  \kappa^2\varepsilon^{-2}(1-\lambda)\Big) \approx  \mathcal{O}\Big(
    \frac{\kappa^2\varepsilon^{-2}}{(1-\lambda)^{2}}\Big), \quad SC \approx \mathcal{O}\Big(\frac{\kappa^2\sqrt{N}\varepsilon^{-2}}{\sqrt{K}(1-\lambda)^{2}}\Big).
\end{align}
\end{Corollary}

$\bullet$ \textbf{Online Scenario}

\begin{Corollary}[\textbf{PAGE+ED}]
\label{appendix:corollary:ed+page:online}
Under Theorem \ref{main:theorem:main}
and assuming $W$ is positive semi-definite in the online scenario, 
if we consider the matrix choice of ED and set $B = \mathcal{O}\Big(\frac{(1-\lambda)^{1.5}T}{K}\Big), b =b_0 = \mathcal{O} \Big(
\frac{(1-\lambda)^{3/4}T^{1/2}}{K}
\Big)$, $p =\mathcal{O}\Big(
\frac{1}{(1-\lambda)^{3/4}T^{1/2}}
\Big)$, 
and $\mu_x = \Theta \Big(\frac{(1-\lambda)^{1.5}}{\kappa^2}\Big), \mu_y = \Theta((1-\lambda)^{1.5})$, 
then it is not difficult to verify that the conditions 
in Theorem \ref{main:theorem:main}
are satisfied for sufficiently large $T$ and small $1-\lambda$. The performance bound is given by 
\begin{align}
&\frac{1}{T}
\sum_{i=0}^{T-1}
(\mE\|\nabla_x J(\bx_{c,i}, \by_{y,i})\|^2 + \mE\|\nabla_y J(\bx_{c,i}, \by_{c,i})\|^2) \notag \\
&\le \mathcal{O}
\Big( 
\frac{\kappa^2 \mE G_{p,0}}{T(1-\lambda)^{1.5}}
+\frac{\kappa^2\zeta^2_0(1-\lambda)}{T}
+ \frac{\kappa^2\lambda^2K\sigma^2}{(1-\lambda)^{3/4}T^{3/2}}
+ \frac{\kappa^2\lambda^2K(1-\lambda)^{5/4}\zeta^2_0}{T^{3/2}} + \frac{\kappa^2\lambda^2K\sigma^2}{(1-\lambda)^{9/4}T^{3/2}}
\notag \\
&\quad +\frac{\kappa^2\mE\Delta_{c,0}}{T(1-\lambda)^{1.5}}
+\frac{\kappa^2\sigma^2}{T}
+\frac{\kappa^2(1-\lambda)^2\zeta^2_0}{T} 
+\frac{\kappa^2\sigma^2}{T(1-\lambda)^{1.5}}\Big)\notag \\
&\le  \mathcal{O}
\Big( 
\frac{\kappa^2 \mE G_{p,0}}{T(1-\lambda)^{1.5}}+ \frac{\kappa^2\mE\Delta_{c,0}}{T(1-\lambda)^{1.5}}
+\frac{\kappa^2\sigma^2}{T(1-\lambda)^{1.5}}
+\frac{\kappa^2\zeta^2_0(1-\lambda)}{T}+\frac{\kappa^2\sigma^2}{T} \notag\\
&\quad 
+ \frac{\kappa^2\lambda^2K\sigma^2}{(1-\lambda)^{9/4}T^{3/2}}+ \frac{\kappa^2\lambda^2K(1-\lambda)^{5/4}\zeta^2_0}{T^{3/2}}
\Big) .
\end{align}
For sufficiently large $T$, the communication complexity and sample complexity are given by 
\begin{align}
   T &= \mathcal{O}\Big(
    \frac{\kappa^2\varepsilon^{-2}}{(1-\lambda)^{1.5}}
    + \kappa^2\varepsilon^{-2} +  \kappa^2\varepsilon^{-2}(1-\lambda)\Big) \approx \mathcal{O}\Big(
    \frac{\kappa^2\varepsilon^{-2}}{(1-\lambda)^{1.5}}\Big) ,  \\
    SC 
    &\approx  \mathcal{O}\Big( \frac{\kappa^2\varepsilon^{-2}}{(1-\lambda)^{1.5}}
\times (1-p) \times b +  \frac{\kappa^2\varepsilon^{-2}}{(1-\lambda)^{1.5}} \times p \times B  + b_0\Big)\notag \\
 &=
  \mathcal{O}\Big( 
 \frac{\kappa^3\varepsilon^{-3}}{(1-\lambda)^{1.5}K}+
 \frac{\kappa^3\varepsilon^{-3}}{(1-\lambda)^{1.5}K} + \frac{\kappa\varepsilon^{-1}}{K}
 \Big).
\end{align}
\end{Corollary}
\begin{Corollary}[\textbf{PAGE+EXTRA}]
\label{appendix:corollary:extra+pag:online}
Under Theorem \ref{main:theorem:main}
and assuming $W$ is positive semi-definite in the online scenario,
if we consider the matrix choice of EXTRA and set $B = \mathcal{O}\Big(\frac{(1-\lambda)^{1.5}T}{K}\Big), b =b_0 = \mathcal{O} \Big(
\frac{(1-\lambda)^{3/4}T^{1/2}}{K}
\Big)$, $p =\mathcal{O}\Big(
\frac{1}{(1-\lambda)^{3/4}T^{1/2}}
\Big)$, 
and $\mu_x = \Theta \Big(\frac{(1-\lambda)^{1.5}}{\kappa^2}\Big), \mu_y = \Theta((1-\lambda)^{1.5})$,
it is then easy to verify that the conditions 
in Theorem \ref{main:theorem:main}
are satisfied for sufficiently large $T$ and small $1-\lambda$. The performance bound is given by 
\begin{align}
&\frac{1}{T}
\sum_{i=0}^{T-1}
(\mE\|\nabla_x J(\bx_{c,i}, \by_{y,i})\|^2 + \mE\|\nabla_y J(\bx_{c,i}, \by_{c,i})\|^2) \notag \\
&\le \mathcal{O}
\Big( 
\frac{\kappa^2 \mE G_{p,0}}{T(1-\lambda)^{1.5}}
+\frac{\kappa^2\zeta^2_0(1-\lambda)}{T}
+ \frac{\kappa^2K\sigma^2}{(1-\lambda)^{3/4}T^{3/2}}
+ \frac{\kappa^2K(1-\lambda)^{5/4}\zeta^2_0}{T^{3/2}} 
\notag \\
&\quad + \frac{\kappa^2K\sigma^2}{(1-\lambda)^{9/4}T^{3/2}}+\frac{\kappa^2\mE\Delta_{c,0}}{T(1-\lambda)^{1.5}}
+\frac{\kappa^2\sigma^2}{T}
+\frac{\kappa^2(1-\lambda)^2\zeta^2_0}{T} 
+\frac{\kappa^2\sigma^2}{T(1-\lambda)^{1.5}}\Big)\notag \\
&\le  \mathcal{O}
\Big( 
\frac{\kappa^2 \mE G_{p,0}}{T(1-\lambda)^{1.5}}+ \frac{\kappa^2\mE\Delta_{c,0}}{T(1-\lambda)^{1.5}}
+\frac{\kappa^2\sigma^2}{T(1-\lambda)^{1.5}}
+\frac{\kappa^2\zeta^2_0(1-\lambda)}{T}+\frac{\kappa^2\sigma^2}{T}  + \frac{\kappa^2K\sigma^2}{(1-\lambda)^{9/4}T^{3/2}}\notag \\
&\quad 
+ \frac{\kappa^2K(1-\lambda)^{5/4}\zeta^2_0}{T^{3/2}}
\Big) .
\end{align}
The communication complexity and sample complexity are given by 
\begin{align}
   T &= \mathcal{O}\Big(
    \frac{\kappa^2\varepsilon^{-2}}{(1-\lambda)^{1.5}}
    + \kappa^2\varepsilon^{-2} +  \kappa^2\varepsilon^{-2}(1-\lambda)\Big) \approx \mathcal{O}\Big(
    \frac{\kappa^2\varepsilon^{-2}}{(1-\lambda)^{1.5}}\Big),  \\
    SC 
    &\approx 
  \mathcal{O}\Big( 
 \frac{\kappa^3\varepsilon^{-3}}{(1-\lambda)^{1.5}K}+
 \frac{\kappa^3\varepsilon^{-3}}{(1-\lambda)^{1.5}K} + \frac{\kappa\varepsilon^{-1}}{K}
 \Big).
\end{align}
\end{Corollary}
\begin{Corollary}[\textbf{PAGE+ATC-GT}]
\label{appendix:corollary:atc-gt+page:online}
Under Theorem \ref{main:theorem:main}
and assuming $W$ is positive semi-definite in the online scenario,
if we consider the matrix choice of ATC-GT and set $B = \mathcal{O}\Big(\frac{(1-\lambda)^{2}T}{K}\Big), b =b_0 = \mathcal{O} \Big(
\frac{(1-\lambda)T^{1/2}}{K}
\Big)$, $p =\mathcal{O}\Big(
\frac{1}{(1-\lambda)T^{1/2}}
\Big)$, 
and $\mu_x = \Theta \Big(\frac{(1-\lambda)^{2}}{\kappa^2}\Big), \mu_y = \Theta((1-\lambda)^{2})$, it is then easy to verify that the conditions 
in Theorem \ref{main:theorem:main}
are satisfied for sufficiently large $T$ and small $1-\lambda$. The performance bound is given by 
\begin{align}
&\frac{1}{T}
\sum_{i=0}^{T-1}
(\mE\|\nabla_x J(\bx_{c,i}, \by_{y,i})\|^2 + \mE\|\nabla_y J(\bx_{c,i}, \by_{c,i})\|^2) \notag \\
&\le \mathcal{O}
\Big( 
\frac{\kappa^2 \mE G_{p,0}}{T(1-\lambda)^{2}}
+\frac{\kappa^2\zeta^2_0(1-\lambda)}{T}
+ \frac{\kappa^2\lambda^4K\sigma^2}{(1-\lambda)T^{3/2}}
+ \frac{\kappa^2\lambda^4K(1-\lambda)\zeta^2_0}{T^{3/2}} + \frac{\kappa^2 \lambda^4 K\sigma^2}{(1-\lambda)^{3}T^{3/2}}+\frac{\kappa^2\mE\Delta_{c,0}}{T(1-\lambda)^{2}}
\notag \\
&\quad 
+\frac{\kappa^2\sigma^2}{T}
+\frac{\kappa^2(1-\lambda)^2\zeta^2_0}{T} 
+\frac{\kappa^2\sigma^2}{T(1-\lambda)^{2}}\Big)\notag \\
&\le  \mathcal{O}
\Big( 
\frac{\kappa^2 \mE G_{p,0}}{T(1-\lambda)^{2}}+ \frac{\kappa^2\mE\Delta_{c,0}}{T(1-\lambda)^{2}}
+\frac{\kappa^2\sigma^2}{T(1-\lambda)^{2}}
+\frac{\kappa^2\zeta^2_0(1-\lambda)}{T}+\frac{\kappa^2\sigma^2}{T}  + \frac{\kappa^2\lambda^4K\sigma^2}{(1-\lambda)^{3}T^{3/2}} \notag\\
&\quad 
+ \frac{\kappa^2\lambda^4K(1-\lambda)\zeta^2_0}{T^{3/2}}
\Big) .
\end{align}
For sufficiently large $T$, the communication complexity and sample complexity are given by 
\begin{align}
   T &= \mathcal{O}\Big(
    \frac{\kappa^2\varepsilon^{-2}}{(1-\lambda)^{2}}
    +\frac{\kappa^2\varepsilon^{-2}}{(1-\lambda)^{2}}
    + \kappa^2\varepsilon^{-2} +  \kappa^2\varepsilon^{-2}(1-\lambda)\Big) \approx \mathcal{O}\Big(
    \frac{\kappa^2\varepsilon^{-2}}{(1-\lambda)^{2}}\Big),  \\
    SC 
    &\approx 
  \mathcal{O}\Big( 
 \frac{\kappa^3\varepsilon^{-3}}{(1-\lambda)^{2}K}+
 \frac{\kappa^3\varepsilon^{-3}}{(1-\lambda)^{2}K} + \frac{\kappa\varepsilon^{-1}}{K}
 \Big).
\end{align}
\end{Corollary}
\subsection{Performance bound of Loopless SARAH-based algorithms}
\label{appendix:corollary:L2S}
It is observed that setting $b$ to a large value  
can significantly improve the communication complexity of PAGE over STORM.
However, increasing the batch size can also be memory-consuming, which is not always feasible in practice.
In the following, we demonstrate the results of 
Loopless SARAH-based algorithm by adopting a minibatch size $b= \mathcal{O}(1)$.
For simplicity, we focus on the offline scenario, and the extension to the online scenario can follow the derivation of PAGE, which is straightforward. To achieve optimal sample complexity, we assume $B= \frac{N}{K}$. In this case, when $N$ is sufficiently large, we can still verify that the terms involving $\mathbb{I}_{\rm online}$ are nearly zero and thus negligible.

\begin{Corollary}[\textbf{Loopless SARAH+ED}]
\label{appendix:corollary:lopplesssarah+ed}
Under theorem \ref{main:theorem:main}
and assume the local sample size $N_k \equiv N$ in an offline scenario,
if we consider the matrix choice of ED and set the following hyperparameters for \textbf{DAMA}
\begin{align}
 b = \mathcal{O}(1), b_0 = 
\mathcal{O}\Big(\frac{\sqrt{N}}{K}\Big), p = \mathcal{O}\Big(\frac{K}{N}\Big), B = \frac{N}{K}.
\end{align}
under a large $N$,  the hyperparamter conditions \eqref{appendix:stepsize:summarya}---\eqref{appendix:stepsize:summaryc} become 
\begin{align}
\mu_x &\le \min \Big\{
\frac{1}{32L}, \frac{\mu_y}{16\kappa^2}, \frac{K}{24\sqrt{3}\kappa L_f\sqrt{N}}\Big\} , \\
\mu_y &\le 
\min \Big\{ 
\frac{1}{\nu}, \frac{1}{2L_f}, \mathcal{O}\Big(\frac{K}{12L_f\sqrt{N}}\Big), \frac{(1-\sqrt{\lambda})\sqrt{1-\lambda}}{\sqrt{620}L_f v_1v_2 \lambda},  \frac{(1-\sqrt{\lambda})\sqrt{1-\lambda}}{12L_f v_1v_2\lambda},\mathcal{O}\Big(\frac{(1-\sqrt{\lambda})\sqrt{1-\lambda}}{24L_f v_1v_2\lambda} \Big), \notag \\
&\mathcal{O}\Big(\frac{(1-\sqrt{\lambda})^{\frac{2}{3}}(1-\lambda)^{1/3}K^{\frac{2}{3}}}{90L_f \kappa^{\frac{1}{3}}v^{\frac{2}{3}}_1 v^{\frac{2}{3}}_2 \lambda^{\frac{2}{3}}N^{\frac{1}{3}}}\Big)
\Big\}, \notag \\
p &\le \frac{\nu \mu_y}{2}, bp\le \frac{1}{K}, p \le 1.
\end{align}
It is easy to verify that the step size choices $\mu_x= \mathcal{O}\Big(\frac{K}{\sqrt{N}\kappa^2}\Big),
\mu_y = \mathcal{O}\Big(
\frac{K}{\sqrt{N}}\Big)$ satisfy the above condition under a large $N$. Then,
the performance bound is given by 
\begin{align}
&\frac{1}{T}
\sum_{i=0}^{T-1}
(\mE\|\nabla_x J(\bx_{c,i}, \by_{y,i})\|^2 + \mE\|\nabla_y J(\bx_{c,i}, \by_{c,i})\|^2) \notag \\
&\le 
\mathcal{O}
\Big(
\frac{\kappa^2\sqrt{N}\mE G_{p,0}}{KT}
+\frac{\kappa^2 \zeta^2_0}{T(1-\lambda)^2}
+ \frac{\kappa^2\lambda^2 K\sigma^2}{\sqrt{N}T(1-\lambda)^3}
+ \frac{\kappa^2\lambda^2K^2\zeta^2_0}{NT(1-\lambda)^4} \notag \\
&\quad 
+\frac{\kappa^2 \sqrt{N}\mE \Delta_{c,0}}{KT}
+ \frac{\kappa^2\sqrt{N}\sigma^2}{KT}
+ \frac{\kappa^2\zeta^2_0}{T(1-\lambda)}
\Big).
\end{align}
When the sample size $N$ is large enough, the communication complexity and sample complexity are approximately given by 
\begin{align}
T &= \mathcal{O}\Big(
\frac{\kappa^2\sqrt{N}\varepsilon^{-2}}{K} +\frac{\kappa^2\varepsilon^{-2}}{(1-\lambda)^2}
\Big)  \approx \mathcal{O}\Big(
\frac{\kappa^2\sqrt{N}\varepsilon^{-2}}{K}
\Big),\\
SC &\approx 
T\times p \times B + T\times (1-p) \times b
+b_0
\approx  \mathcal{O}\Big(\frac{\kappa^2 \sqrt{N} \varepsilon^{-2}}{K}
+\frac{\kappa^2\varepsilon^{-2}}{(1-\lambda)^2}
+\frac{\sqrt{N}}{K}
\Big).
\end{align}
\end{Corollary}

\begin{Corollary}[\textbf{\textbf{Loopless SARAH+EXTRA}}]
\label{appendix:corollary:looplesssarch+extra}
Under theorem \ref{main:theorem:main}
and assume the local sample size $N_k \equiv N$ in an offline scenario, 
if we consider the matrix choice of EXTRA and set the hyperparameters for \textbf{DAMA} as follows
\begin{align}
\mu_x= \mathcal{O}\Big( \frac{K}{\sqrt{N}\kappa^2} \Big),
\mu_y = \mathcal{O}\Big(
\frac{K}{\sqrt{N}} \Big), b = \mathcal{O}(1), b_0 = 
\mathcal{O}\Big(\frac{\sqrt{N}}{K}\Big) , p =  \mathcal{O}\Big(\frac{K}{N}\Big), B = \frac{N}{K}.
\end{align}
It is not difficult to verify that the hyperparameter conditions in Theorem 
\ref{main:theorem:main}
are satisfied for sufficiently large $N$.
The performance bound is given by 
\begin{align}
&\frac{1}{T}
\sum_{i=0}^{T-1}
(\mE\|\nabla_x J(\bx_{c,i}, \by_{y,i})\|^2 + \mE\|\nabla_y J(\bx_{c,i}, \by_{c,i})\|^2) \notag \\
&\le 
\mathcal{O}
\Big(
\frac{\kappa^2\sqrt{N}\mE G_{p,0}}{KT}
+\frac{\kappa^2  \zeta^2_0}{T(1-\lambda)^2}
+ \frac{\kappa^2 K\sigma^2}{\sqrt{N}T(1-\lambda)^3}
+ \frac{\kappa^2K^2\zeta^2_0}{NT(1-\lambda)^4} \notag \\
&\quad 
+\frac{\kappa^2 \sqrt{N}\mE \Delta_{c,0}}{KT}
+ \frac{\kappa^2\sqrt{N}\sigma^2}{KT}
+ \frac{\kappa^2\zeta^2_0}{T(1-\lambda)}
\Big).
\end{align}
when the sample size $N$ is large enough, the communication complexity and sample complexity are approximated given by 
\begin{align}
T &= \mathcal{O}\Big(
\frac{\kappa^2\sqrt{N}\varepsilon^{-2}}{K} +\frac{\kappa^2\varepsilon^{-2}}{(1-\lambda)^2}
\Big) \approx \mathcal{O}\Big(
\frac{\kappa^2\sqrt{N}\varepsilon^{-2}}{K}
\Big), \\
SC &\approx 
T\times p \times B + T\times (1-p) \times b
+b_0
\approx  \mathcal{O}\Big(\frac{\kappa^2 \sqrt{N} \varepsilon^{-2}}{K}
+\frac{\kappa^2\varepsilon^{-2}}{(1-\lambda)^2}
+ \frac{\sqrt{N}}{K}
\Big).
\end{align}
\end{Corollary}

\begin{Corollary}[\textbf{\textbf{Loopless SARAH+ATC-GT}}]
\label{appendix:corollary:looplesssarch+atc-gt}
Under theorem \ref{main:theorem:main}
and assume the local sample size $N_k \equiv N$ in an offline scenario, consider the matrix choice for ATC-GT and 
if we set the hyperparameters for \textbf{DAMA} as follows
\begin{align}
\mu_x= \frac{K}{\sqrt{N}\kappa^2},
\mu_y = 
\frac{K}{\sqrt{N}}, b = 1,b_0 = 
\mathcal{O}\Big(\frac{\sqrt{N}}{K}\Big), p = \frac{K}{N}, B = \frac{N}{K}.
\end{align}
It is not difficult to verify that the conditions in Theorem 
\ref{main:theorem:main}
are satisfied for sufficiently large $N$.
The performance bound is given by 
\begin{align}
&\frac{1}{T}
\sum_{i=0}^{T-1}
(\mE\|\nabla_x J(\bx_{c,i}, \by_{y,i})\|^2 + \mE\|\nabla_y J(\bx_{c,i}, \by_{c,i})\|^2) \notag \\
&\le 
\mathcal{O}
\Big(
\frac{\kappa^2\sqrt{N}\mE G_{p,0}}{KT}
+\frac{\kappa^2 \zeta^2_0}{T(1-\lambda)^3}
+ \frac{\kappa^2 \lambda^4 K\sigma^2}{\sqrt{N}T(1-\lambda)^4}
+ \frac{\kappa^2\lambda^4K^2\zeta^2_0}{NT(1-\lambda)^6} \notag \\
&\quad 
+\frac{\kappa^2 \sqrt{N}\mE \Delta_{c,0}}{KT}
+ \frac{\kappa^2\sqrt{N}\sigma^2}{KT}
+ \frac{\kappa^2\zeta^2_0}{T(1-\lambda)^2}
\Big).
\end{align}
when the sample size $N$ is large enough, the communication complexity and sample complexity are approximated given by 
\begin{align}
T &= \mathcal{O}\Big(
\frac{\kappa^2\sqrt{N}\varepsilon^{-2}}{K} + \frac{\kappa^2\varepsilon^{-2}}{(1-\lambda)^3}
\Big)  \approx \mathcal{O}\Big(
\frac{\kappa^2\sqrt{N}\varepsilon^{-2}}{K} 
\Big) ,\\
SC &\approx 
T\times p \times B + T\times (1-p) \times b
+b_0
\approx  \mathcal{O}\Big(\frac{\kappa^2 \sqrt{N} \varepsilon^{-2}}{K}
+ \frac{\kappa^2\varepsilon^{-2}}{(1-\lambda)^3}
+\frac{\sqrt{N}}{K}
\Big).
\end{align}
\end{Corollary}

\end{document}